\documentclass[11pt,a4paper]{article}

\usepackage{comment}
\usepackage{graphicx}
\usepackage{amsmath}
\usepackage{amsfonts}
\usepackage{amssymb}
\usepackage{enumerate}
\usepackage{lscape}
\usepackage{longtable}
\usepackage{rotating}
\usepackage{color}
\usepackage{url}
\usepackage{rotating}
\usepackage{algpseudocode}
\usepackage[ruled]{algorithm}
\usepackage{pgf,tikz}
\usepackage{caption}
\usepackage[labelfont=sf]{subcaption}
\usepackage[titletoc]{appendix}
\usepackage{wrapfig}
\usepackage{adjustbox}
\usepackage{mathrsfs}
\usepackage{pgfplots}
\usetikzlibrary{arrows}
\usepackage{color}

\usepackage{bbm}

\usepackage{booktabs}
\usepackage{multirow}
\usepackage{mathtools}
\usepackage{makecell}

\usepackage[sort,numbers]{natbib}
\usepackage{eqparbox}%

\numberwithin{equation}{section}%

\usepackage{bm}

\newtheorem{definition}{Definition}
\newtheorem{example}{Example}[section]

\newtheorem{remark}{Remark}[section]

\newtheorem{assumption}{Assumption}[section]

\DeclareMathOperator{\MCI}{\tcb{MCI}}
\DeclareMathOperator{\MCF}{MCF}
\DeclareMathOperator{\MCFg}{gMCF}

\newcommand{\beqn}[1]{\begin{equation}\label{#1}}
\newcommand{\eeqn}{\end{equation}}

\definecolor{darkgreen}{rgb}{0,0.6,0}
\definecolor{aau2}{rgb}{0.0, 0.5, 0.69}
\definecolor{aau3}{rgb}{0.0, 0.53, 0.74}
\definecolor{aau4}{rgb}{0.0, 0.48, 0.65}
\definecolor{aau5}{rgb}{0.0, 0.45, 0.73}
\definecolor{rsap}{RGB}{130, 36, 51}
\definecolor{gsap}{RGB}{112, 164, 137}
\definecolor{tud}{rgb}{0.43,0.73,0.11}
\definecolor{verde}{rgb}{0.33,0.53,0.11}
\definecolor{ttffqq}{rgb}{0.0, 0.48, 0.65} 
\definecolor{ffqqqq}{rgb}{0.0, 0.5, 0.69} 

\newcommand{\tcb}{}

\usetikzlibrary{arrows.meta, calc, backgrounds}

\tikzstyle{decision} = [diamond, draw, fill=blue!20,
text width=4.5em, text badly centered, node distance=3cm, inner sep=0pt]
\tikzstyle{block} = [rectangle, draw, fill=blue!20,
text centered, rounded corners, minimum height=4em]
\tikzstyle{line} = [draw, -latex']
\tikzstyle{cloud} = [draw, ellipse,fill=red!20, node distance=3cm,
minimum height=2em]
\tikzstyle{cloud2} = [draw, ellipse,fill=green!20, node distance=3cm,
minimum height=2em]

\begin{document}

	\title{Pareto sensitivity, most-changing sub-fronts, \\and knee solutions}
	
	\author{
		T. Giovannelli\thanks{Department of Mechanical, Materials, and Industrial Engineering, University of Cincinnati, Cincinnati, OH 45221, USA ({\tt giovanto@ucmail.uc.edu}).}
		\and
		M. M. Raimundo\thanks{Departamento de Sistemas de Informa\c{c}\~ao, UNICAMP - Universidade Estadual de Campinas, Campinas, São Paulo, Brasil ({\tt mraimundo@ic.unicamp.br}).}
		\and
		L. N. Vicente\thanks{Department of Industrial and Systems Engineering, Lehigh University, Bethlehem, PA 18015-1582, USA ({\tt lnv@lehigh.edu}).}
	}
	
	\maketitle
 
	\begin{abstract}
              When dealing with a multi-objective optimization problem, \tcb{identifying Pareto optimal solutions that capture relevant trade-offs may not be straightforward. Popular choices include} the so-called Pareto knee solutions, which correspond to nondominated points on the Pareto front where a small improvement in any objective leads to a large deterioration in at least one other objective.
      In this paper, using Pareto sensitivity, \tcb{we relate Pareto knee solutions to Geoffrion's concept of proper efficiency and show how to compute them} according to their verbal (informal) definition of least maximal change. We refer to the resulting approach as the sensitivity knee (snee) approach, and we apply it to unconstrained and constrained problems. Pareto sensitivity can also be used to compute local most-changing Pareto sub-fronts around a nondominated point, where points on the sub-fronts are distributed along directions of maximum change.
      Our approach is restricted to scalarized methods, in particular to the weighted-sum or epsilon-constrained methods, and requires the computation or approximations of first- and second-order derivatives. We include numerical results from synthetic problems that \tcb{illustrate our approach}.	 
	\end{abstract}

\section{Introduction}\label{sec:intro} 

In this paper, we focus on the following multi-objective optimization~(MOO) problem	
\begin{equation}\label{prob:multiobj}
    \begin{aligned}
        \min_{x \in X \subseteq \mathbb{R}^n} \ & F(x) \; = \; \left(f_{1}(x), \ldots, f_{q}(x)\right), 
    \end{aligned}
\end{equation}
where each~$f_i: \mathbb{R}^n \to \mathbb{R}$ is a twice continuously differentiable function\tcb{,} for all~$i \in \{1,\ldots,q\}$, with~$q \ge 2$, and~$X$ represents the feasible set.
\tcb{We consider both unconstrained problems, with~$X = \mathbb{R}^n$, and constrained problems, with~$X \subset \mathbb{R}^n$. Additional regularity conditions for the unconstrained and constrained cases are introduced later in Assumption~\ref{ass:IFT} (unconstrained
case, Section~\ref{sec:unconstr}) and
Assumption~\ref{ass:LL_assumptions_constr} (constrained case,
Section~\ref{sec:constr}), respectively, together with Assumption~\ref{ass:convexity} (Section~\ref{subsec:point_based}). In particular, strong convexity of the objective functions~$f_1, \ldots, f_q$ is a sufficient condition for the assumptions used in the unconstrained case, whereas in the constrained case we assume a convex feasible set and impose second-order sufficient optimality conditions, the linear independence constraint qualification~(LICQ), and strict complementarity.}

Problems with multiple objectives\tcb{, such as problem~\eqref{prob:multiobj},} find applications across a wide range of domains. 
In engineering design, multiple objectives correspond to various performance measures (e.g., efficiency, reliability, and safety) that must be maximized simultaneously~\citep{PJFleming_etal_2005,GChiandussi_etal_2012,JHakanen_etal_2021,RHStewart_etal_2021,MJameel_etal_2026}. In finance, MOO is at the heart of portfolio optimization problems, where trade-offs between risk and return need to be balanced to achieve the best investment strategy~\citep{MCSteinbach_2001,YChen_AZhou_2022,ASulas_etal_2025}.
 In the healthcare sector, MOO is used to strike a balance between providing equitable or rapid service to patients and managing operational costs effectively~\citep{WZhang_etal_2016,MBoresta_etal_2024}. In transportation systems, MOO helps optimize routing and speed decisions to minimize environmental impact while maximizing efficiency for service providers and ensuring high service quality for users~\citep{EDemir_etal_2014,TGiovannelli_LNVicente_2023}. In the field of machine learning, MOO problems are particularly relevant in scenarios like multi-task learning~\citep{LXi_2019,LChen_HFernando_YYing_TChen_2023}, where a model must perform well on several related conflicting tasks simultaneously, and in learning problems that address fairness and privacy concerns~\citep{ANavon_etal_2020,SLiu_LNVicente_2020}.

 \tcb{To formalize the discussion, we briefly introduce the terminology used throughout Subsections~\ref{subsec:intro_pareto} and~\ref{subsec:intro_pareto_sensitivity},} and formally defined in Definitions~\ref{def:pareto} and~\ref{def:weakpareto} of Section~\ref{sec:basic_results}. Specifically, \tcb{a feasible solution to problem~\eqref{prob:multiobj} is called Pareto optimal (or efficient) if there is no feasible solution that improves at least one objective without worsening any of the others. The set of all Pareto optimal solutions is called the efficient set.} 
 The Pareto front of problem~\eqref{prob:multiobj}, also referred to as the nondominated set or Pareto set, is the set obtained by mapping the Pareto optimal solutions from the decision space~$\mathbb{R}^n$ to the objective space~$\mathbb{R}^q$. Each point on the Pareto front is referred to as a nondominated point.
 \tcb{Weak Pareto optimality (or weak efficiency) relaxes Pareto optimality by requiring only that there be no feasible solution that strictly improves all objective functions simultaneously. The corresponding points in the objective space are called weakly nondominated points.}

\subsection{Finding relevant Pareto optimal solutions}
\label{subsec:intro_pareto}

Determining a large number of Pareto optimal solutions can be computationally expensive, particularly with a high number of objectives~\citep{HWang_etal_2017,XLin_YLiu_etal_2024}. While all Pareto optimal solutions are mathematically valid solutions, it is generally preferable to identify only a few \tcb{relevant} solutions~\citep{ETriantaphyllou_2000,NBoehmer_etal_2025}. 
Well-known multi-criteria decision-making techniques have been proposed for selecting a single Pareto optimal solution. 
\textit{A priori methods}~\citep{MEhrgott_2005,KMiettinen_1999} consist of reducing an~MOO problem into a single-objective one by incorporating user preferences, which can be achieved by either weighting the objective functions into a single function with user-defined weights (weighted-sum method), minimizing one objective function and considering the other objectives as constraints with user-defined right-hand sides~($\varepsilon$-constraint method), minimizing a user-defined utility function (utility-based method), or minimizing the deviation of the objectives from user-defined target values (also known as goal programming). The so-called \textit{knee solutions}\footnote{Throughout this paper, we interpret a \textit{(Pareto) knee solution} as a Pareto optimal solution in the decision space~$\mathbb{R}^n$, and refer to the corresponding nondominated point in the objective space~$\mathbb{R}^q$ as a \textit{knee nondominated point}.}~\citep{Burke_etal_1988,IDas_1999,KMiettinen_1999,KDeb_2003} are Pareto optimal solutions that correspond to nondominated points on the Pareto front where a small improvement in any objective would lead to a large deterioration in at least one other objective, which is a commonly-used \textbf{verbal (informal) definition}.
\textit{Sharpe solutions}~\citep{Sharpe_1966,SLiu_LNVicente_2021} have only been proposed for bi-objective problems and are defined as Pareto optimal solutions with the smallest function value in the first objective per unit of function value in the second objective. 

In addition to a priori methods, other widely studied classes of multi-criteria decision-making approaches include \textit{a posteriori} and \textit{interactive methods}~\citep{MEhrgott_2005,KMiettinen_1999,BAfsar_KMiettinen_FRuiz_2021}. 
A posteriori methods aim to approximate \tcb{the entire nondominated set in the objective space,} after which a decision maker selects a preferred \tcb{nondominated point} in a post-processing step. 
Interactive methods iteratively generate \tcb{Pareto optimal solutions or nondominated points} while progressively incorporating decision-maker preferences throughout the optimization process\tcb{, with the goal of identifying a preferred solution or point}. 

In this work, we focus on knee solutions, which we identify using a sensitivity-based approach that does not require explicit user preferences, applies to an arbitrary number of objective functions, avoids computing an approximation of the \tcb{nondominated} set, and does not involve iterative interaction with a decision maker.
Knee solutions have gained significant popularity and have been extensively studied in the literature. A review of methods to compute knee solutions~(and the corresponding knee nondominated points) is provided in~\cite{JBranke_etal_2004} and~\cite{KDeb_SGupta_2010}. However, \tcb{although several definitions of knee solutions have been proposed, no widely accepted quantitative definition exists,} and identifying knee solutions in high-dimensional objective spaces remains challenging. Furthermore, as discussed in~\citet[Subsection~2.2]{KDeb_SGupta_2010}, most current methods identify nondominated points that do not necessarily correspond to knee
solutions based on the verbal definition, but are simply intermediate points on the Pareto front, where the trade-off between the improvement in any objective and the corresponding deterioration in at least one other objective is not necessarily significant.
Throughout this paper, following \tcb{the} terminology adopted in~\cite{KDeb_SGupta_2010}, \tcb{when discussing bi-objective and tri-objective problems, we informally} refer to the intermediate region of the Pareto front \tcb{as the middle portion of the front, away from its two endpoints~(when~$q=2$) or its edge curves (when~$q=3$)}. 

Existing methods for determining knee solutions generally fall into three categories: angle-based, utility-based, and normal boundary intersection-based approaches. Such methods are described below, highlighting their motivations and limitations. Our sensitivity-based approach, which is alternative to these three categories, is described in Subsection~\ref{subsec:intro_pareto_sensitivity}.  

Angle-based methods, which are applicable only to problems with two objective functions, identify knee nondominated points by measuring reflex angles~\citep{JBranke_etal_2004} or bend angles~\citep{KDeb_SGupta_2010} using a predetermined set of nondominated points. Specifically, a reflex angle is formed between a nondominated point and two chosen neighboring nondominated points (one to the left and one to the right). Since a reflex angle is determined using neighboring solutions, it reflects a local characteristic of the Pareto front, describing its behavior in the vicinity of a specific nondominated point. A bend angle is similar to a reflex angle but is obtained by replacing neighboring nondominated points with \tcb{single-objective optimal nondominated points}\footnote{In general, given~$q$ objectives, the corresponding Pareto front has at most~$q$ {\it \tcb{single-objective optimal nondominated points}} in the objective space, given by~$F(x_*^i)$, where~$x_*^i$ is a minimizer of~$f_i$, with~$i \in \{1,\ldots,q\}$.}, which are nondominated points that correspond to the minimizers of each objective function individually. Since a bend angle is calculated using \tcb{such nondominated points}, it is better able to capture the global behavior of the Pareto front. In general, reflex and bend angles approximate the shape of the Pareto front at a given nondominated point, and the knee solution is the Pareto optimal solution corresponding to the nondominated point with the maximum reflex or bend angle. In~\cite{KDeb_SGupta_2010}, the authors propose using a problem-specific, user-defined threshold to determine whether a bend angle is sharp enough to classify the corresponding solution as a knee solution. If the threshold is set too high, a knee may not exist, and nondominated points that are practically significant may be excluded without proper fine-tuning of the threshold. 


According to~\citet{JBranke_etal_2004}, utility-based methods identify knee solutions using a linear utility function of the form~$U(F(x),\mu)=\sum_{i=1}^{q} \mu_i f_i(x)$, where~$\mu$ lies in the simplex set.
 To identify a knee solution, one can assume a uniform distribution over the values of~$\mu$ and maximize the expected value of a {\it marginal} utility function over a predetermined set of Pareto optimal solutions. 
 For a given value of~$\mu$, the marginal utility of a Pareto optimal solution quantifies the loss in utility that would occur if such a solution were unavailable and the decision maker had to select the best remaining solution. In this sense, Pareto optimal solutions with high marginal utility are those whose removal would lead to a significant degradation in achievable utility, and are therefore interpreted as knee solutions; we refer to~\citet[Section 3.2]{JBranke_etal_2004} for the details.
 A variant of such a method is provided in~\citet[Definition~3]{KDeb_SGupta_2010}, which defines a knee solution as a Pareto optimal solution that minimizes~$U(F(\cdot),\mu)$ over a predetermined set of Pareto optimal solutions for the largest number of weight vectors~$\mu$ chosen from the simplex set. Consequently, one must solve a minimization problem for each selected~$\mu$ and count how often each Pareto optimal solution is obtained. The solution with the highest count is then defined as the knee solution. 

Normal boundary intersection-based methods, applicable to an arbitrary number of objective functions, search for knee nondominated points on the Pareto front by maximizing their distance to the convex hull of the \tcb{single-objective optimal nondominated points}, referred to as the boundary line. Such knee nondominated points are typically found in the intermediate region of the Pareto front, especially when the front is convex~\citep{KDeb_SGupta_2010}.  
To identify knee nondominated points and the corresponding knee solutions, one can employ a nonlinear constrained optimization approach, as introduced in~\cite{IDas_1999} and \cite{IDas_JEDennis_1998}, and further refined in~\cite{PKShukla_2013}. 
Evolutionary algorithms were used to search for normal boundary intersection-based knee solutions in~\cite{XZhang_etal_2015} and \cite{SBechikh_etal_2010}.

\subsection{Computing knee solutions through Pareto sensitivity}
\label{subsec:intro_pareto_sensitivity}

As an advancement over existing multi-criteria decision-making techniques, the main contribution of our paper is the introduction of a novel approach for determining knee solutions and knee nondominated points through Pareto sensitivity. By Pareto sensitivity, we mean the application of sensitivity analysis principles to obtain a quantitative measure of how sensitive a Pareto optimal solution or a nondominated point is to variations in problem parameters. In our paper, such parameters correspond to the weights arising from the weighted-sum method (or to the right-hand side~$\varepsilon$ when using the~$\varepsilon$-constraint method). We note that the study of sensitivity analysis in~MOO dates back to the late 1980s~\citep{TTanino_1988,DRiosInsua_1990}, where it was investigated through set-valued perturbation mappings and variational analysis frameworks, and has continued through more recent work~\citep{KAAndersen_TKBoomsma_etal_2024}. 
Related sensitivity techniques have recently gained attention in the bilevel and multilevel optimization literature, where they are used to compute the gradient of the objective function in the single-level reformulation of bilevel or trilevel problems~\citep{SGhadimi_MWang_2018,KJi_JYang_YLiang_2020,TChen_YSun_WYin_2021,TGiovannelli_GKent_LNVicente_2025_BSG,TGiovannelli_GKent_LNVicente_2024,TGiovannelli_GDKent_LNVicente_2025_trilevel}. 

 Consider the weighted-sum function associated with problem~\eqref{prob:multiobj},~$\sum_{i=1}^{q} \lambda_i f_i(x)$, where $\lambda$ takes values in the simplex set~$\Lambda$. Note that such a weighted-sum function can be interpreted as the linear utility function~$U(F(x),\lambda)$ introduced in Subsection~\ref{subsec:intro_pareto}, with~$\lambda$ in place of~$\mu$. 
\tcb{For each~$\lambda \in \Lambda$, any global minimizer~$x(\lambda)$ of the weighted-sum function over~$x \in X$ is a weakly Pareto optimal solution.} \tcb{Under the convexity condition of Assumption~\ref{ass:convexity} (Section~\ref{subsec:point_based}), $x(\lambda)$ is Pareto optimal, rather than only weakly Pareto optimal. Under the assumptions needed for the application of the implicit function theorem (collected in Assumption~\ref{ass:IFT} in Subsection~\ref{subsec:pareto_sensitivity_calculation} for the unconstrained case~$X = \mathbb{R}^n$, and in Assumption~\ref{ass:LL_assumptions_constr} in Subsection~\ref{subsec:pareto_sensitivity_calculation_constr} for the constrained case~$X \subset \mathbb{R}^n$), such an~$x(\lambda)$ is unique.
To define and determine a Pareto knee solution that follows the spirit of the verbal definition, we connect knee solutions to the classical notion of Geoffrion proper efficiency (see Definition~\ref{def:knee} in Section~\ref{subsec:point_based}). A Pareto optimal solution is properly efficient, in the sense of Geoffrion, when the trade-offs required to improve one objective at the expense of another remain bounded. 
Building on this notion, for a given~$\lambda \in \Lambda$, we define a regularized Geoffrion trade-off constant, referred to as the \textit{maximal-change} function~(MCF), which measures the largest local trade-off ratio attainable around the nondominated point corresponding to~$x(\lambda)$. We identify knee solutions by minimizing the~MCF over~$\lambda \in \Lambda$, thereby selecting Pareto optimal solutions at which the maximal local trade-off is minimized.}

\tcb{Evaluating the~MCF at a given~$\lambda$ requires, in general, solving an auxiliary optimization problem over a continuum of directions, which can be computationally expensive. When this computational cost is a concern, we therefore introduce a computationally tractable proxy for the~MCF, which we refer to as the \textit{gradient maximal-change} function~(gMCF). The~gMCF depends only on the magnitudes of the objective function gradients, admits a closed-form expression, and is considerably cheaper to evaluate and minimize than the~MCF. Minimizing the~gMCF typically identifies knee solutions similar to those obtained by minimizing~MCF directly, although the two functions can occasionally lead to different solutions, as demonstrated by the computational experiments. When a knee solution is selected using either the~MCF or the~gMCF,} the decision maker is guaranteed that trading among objectives is optimally balanced, minimizing the extent of compromise required for one objective relative to the others. We refer to the resulting approach as the \textit{sensitivity knee} (snee) approach.
We note that the techniques used in our~\textit{snee} approach require the computation or approximations of first- and second-order derivatives and are restricted to scalarized methods, in particular to the weighted-sum method considered in this paper.

Since the \tcb{the~MCF and the~gMCF}
may be non-convex (regardless of the convexity of the Pareto front)
and thus may have more than one local minimizer, \tcb{the knee solutions determined through our \textit{snee} approach} may have a local or global nature. Moreover, \tcb{such} knee solutions may not necessarily correspond to nondominated points in the intermediate region of the Pareto front (we again emphasize that by intermediate region, we mean \tcb{the middle portion of the Pareto front}). The concept of knee nondominated points outside the intermediate region was analyzed in~\citet[Section~5]{KDeb_SGupta_2010}, where such points, located near or at the \tcb{single-objective optimal nondominated points} of the Pareto front, are referred to as {\it edge-knee}. However, the definition of an edge-knee nondominated point applies only to~MOO problems with~2~objectives and depends on a problem-specific, user-defined parameter, whereas our approach overcomes such limitations.

Given~$\lambda \in \Lambda$, let~$x(\lambda)$ be the corresponding solution obtained by minimizing the weighted-sum function~$\sum_{i=1}^{q}\lambda_i f_i(x)$ over~$x \in X$. 
We will also see in Section~\ref{sec:neighborhoods} of this paper that the Jacobian matrix of the vector function~$F(x(\lambda))$ with respect to~$\lambda$ can be used to define a neighborhood of weight vectors around~$\lambda$ (analogously, a neighborhood of Pareto optimal solutions around~$x(\lambda)$), where the objective functions~$f_i$ exhibit the greatest variation per unit change of the weights. 
Local most-changing Pareto sub-fronts can then be computed by mapping the weight vectors in such a neighborhood through~$F(x(\lambda))$ to the objective space~$\mathbb{R}^q$. Nondominated points on these sub-fronts are distributed along directions of maximum change on the Pareto front. 
Such sub-fronts can assist decision makers in locally exploring different trade-offs.

To derive Pareto sensitivity in a convenient way and calculate the gradient of~$f_i(x(\lambda))$ with respect to~$\lambda$, required by our~\textit{snee} approach, with~$i \in \{1, \ldots, q\}$, we focus on the first-order necessary conditions for minimizing the weighted-sum function~$\sum_{i=1}^q \lambda_i f_i(x)$ over~$x \in X$, given~$\lambda \in \Lambda$. By applying the chain rule to \tcb{such} conditions, we will first characterize the Jacobian, or sensitivity matrix,~$\nabla x(\lambda)^\top$, and then calculate the gradient of~$f_i(x(\lambda))$ with respect to~$\lambda$. 
All these steps will be formally described in Sections~\ref{sec:unconstr} and~\ref{sec:constr}.

Several advantages are offered by our \textit{snee} approach 
compared to existing techniques:
\begin{enumerate}
    \item Change on the Pareto front is quantified using the most universal concept of rate of change, i.e., first-order derivatives, eliminating the need for problem-specific, user-defined parameters.
    \item It offers the flexibility to compute knee solutions according to the verbal definition, but in a quantifiable and optimal form, searching the Pareto front for local or global knee solutions, without relying on an approximation of the front in advance.
    \item It can be applied to any number~$q$ of objective functions.
    \item It detects a knee \tcb{nondominated point} in the intermediate region of the Pareto front only if the trade-off between improvement and deterioration in the objectives is significant enough\tcb{, rather than forcing it into the intermediate region regardless of whether such a trade-off exists.}
     \item The approach also allows for the identification of most-changing Pareto sub-fronts of predefined size, rather than focusing on a single \tcb{nondominated point}.
\end{enumerate}

Our paper is organized as follows. In Section~\ref{sec:basic_results}, we present basic~MOO definitions and facts that will be used throughout the paper. In~Section~\ref{subsec:point_based}, we \tcb{define knee solutions and introduce the \textit{snee} approach \tcb{by developing optimization formulations for identifying them}}. 
In~Section~\ref{sec:neighborhoods}, we address how to compute most-changing Pareto sub-fronts around a nondominated point.
In Sections~\ref{sec:unconstr} and~\ref{sec:constr}, we introduce general assumptions for problem~\eqref{prob:multiobj} and derive Pareto sensitivity formulas for the unconstrained case~($X = \mathbb{R}^n$) and the constrained case~($X \subset \mathbb{R}^n$), respectively, including numerical results that \tcb{illustrate our approach}. Finally, we draw some concluding remarks and ideas for future work in Section~\ref{sec:future_work}. 

All code was implemented in Python~3.12.\tcb{3}. The experimental results 
were obtained using a \tcb{Lenovo laptop (model~21FBS17600) with~128\,GB of~RAM and an Intel(R) Core(TM) i9-13950HX processor}.
Throughout the paper, $\|\cdot\|$ denotes the Euclidean norm~$\|\cdot\|_2$. 

\section{MOO background}\label{sec:basic_results} 

Consistent with common conventions in the~MOO literature~\citep[Table 1.2]{MEhrgott_2005}, we adopt the following notation. Let~$y^1, y^2 \in \mathbb{R}^q$. 
We write~$y^1 \leq y^2$ if $y^1_i \leq y^2_i$ for all~$i \in \{1, \ldots, q\}$, with~$y^1 \neq y^2$ (componentwise order), and~$y^1 < y^2$ if~$y^1_i < y^2_i$ for all~$i \in \{1, \ldots, q\}$ (strict componentwise order). The image of the feasible set~$X$ under the vector-valued function~$F$ is denoted as~$Y = F(X)$ \tcb{and called the image set}.
Pareto optimal solutions (or efficient solutions) of problem~\eqref{prob:multiobj} and the corresponding weakly versions are introduced in Definitions~\ref{def:pareto} and~\ref{def:weakpareto} below~\citep[p. 24, Definition~2.1, and p. 38, Definition~2.24, respectively]{MEhrgott_2005}.

\begin{definition}\label{def:pareto}
   A feasible \tcb{solution}~$\hat{x} \in X$ is called Pareto optimal or efficient if there is no other~$x \in X$ such that~$F(x) \le F(\hat{x})$. If~$\hat{x}$ is Pareto optimal, $F(\hat{x})$ is called nondominated point. If~$x^1, x^2 \in X$ and~$F(x^1) \le F(x^2)$, we say~$x^1$ dominates~$x^2$ and~$F(x^1)$ dominates~$F(x^2)$. The set of all Pareto optimal solutions~$\hat{x} \in X$ is denoted~$X_E$ and called the efficient set. The set of all nondominated points~$\hat{y} = F(\hat{x}) \in Y$, where~$\hat{x} \in X_E$, is denoted~$Y_N$ and called the nondominated set.
\end{definition}

\begin{definition}\label{def:weakpareto}
    A feasible \tcb{solution}~$\hat{X} \in X$ is called weakly Pareto optimal or weakly efficient if there is no~$x \in X$ such that~$F(x) < F(\hat{x})$. The point~$\hat{y} = F(\hat{x})$ is then called weakly nondominated. The weakly efficient and nondominated sets are denoted~$X_{wE}$ and~$Y_{wN}$, respectively.
\end{definition}

To be consistent with a significant portion of the multi-objective optimization literature~\citep{KDeb_2003,KDeb_1999_MOO,ASinha_KDeb_2009,KDeb_etal_2002,KDeb_SGupta_2010,SFallah_etal_2024,VJBowman_1976,MDeSantis_etal_2020}, we also refer to~$Y_N$ as the Pareto front, and we note that alternative terminology, such as Pareto set, is also used~\citep{LWeerasena_etal_2017}. Note that from Definitions~\ref{def:pareto} and~\ref{def:weakpareto}, it follows that~$Y_N \tcb{\subseteq} Y_{wN}$ and~$X_E \tcb{\subseteq} X_{wE}$.

\tcb{A refinement of efficiency that will play a central role in this paper is Geoffrion's notion of proper efficiency, recalled in Definition~\ref{def:geoffrion} below~\citep[p. 51, Definition~2.39]{MEhrgott_2005}.
\begin{definition}\label{def:geoffrion}
A feasible solution~$\hat{x} \in X$ is called properly efficient if it is efficient and if there is a real number~$M > 0$ such that for all~$i \in \{1, \ldots, q\}$ and~$x \in X$ satisfying~$f_i(x) < f_i(\hat{x})$, there exists an index~$j \in \{1, \ldots, q\}$ such that~$f_j(\hat{x}) < f_j(x)$ and
\begin{equation}\label{eq:geoffrion_def}
\frac{f_i(\hat{x}) - f_i(x)}{f_j(x) - f_j(\hat{x})} \; \leq \; M.
\end{equation}
The corresponding point~$\hat{y} = F(\hat{x})$ is called properly nondominated. The properly efficient and nondominated sets are denoted~$X_{pE}$ and~$Y_{pN}$, respectively.
\end{definition}}

Definition~\ref{def:ideal_nadir} below, based on~\citet[Definition~2.22]{MEhrgott_2005}, defines ideal and nadir points, which are the points in the objective space associated with the best and worst possible values for all objectives, respectively.
Such points can be used to define the denominator of the \tcb{indicator} that we will introduce to assess how well a local Pareto sub-front aligns with the direction of maximum change on the Pareto front (see~\eqref{eq:metric} in Section~\ref{sec:neighborhoods}).  
\begin{definition}[Ideal and nadir points]\label{def:ideal_nadir}
   The ideal point~$F^I \in \mathbb{R}^q$ is the vector whose $i$-th component is given by~$F_i^I = \min_{x \in X} f_i(x) = \min_{y \in Y} y_i$, for all~$i \in \{1,\ldots,q\}$. 
   The nadir point~$F^N \in \mathbb{R}^q$ is the vector whose $i$-th component is given by~$F_i^N = \max_{x \in X_E} f_i(x) = \max_{y \in Y_N} y_i$, for all~$i \in \{1,\ldots,q\}$.
\end{definition}

\tcb{Definition~\ref{def:geoffrion} states that~$\hat{x}$ is properly efficient when the trade-off required to improve one objective at the expense of another remains bounded by some constant~$M$. Figure~\ref{fig:geoffrion}
illustrates this definition in the bi-objective case. Since every properly efficient solution is Pareto optimal, we have~$X_{pE} \subseteq X_E$ and~$Y_{pN} \subseteq Y_N$.}

To compute \tcb{weakly} Pareto optimal solutions, one can use scalarization techniques to reduce a multi-objective problem into a single-objective one and then apply classical optimization methods. One popular scalarization technique is the weighted-sum method~\citep{MEhrgott_2005,KMiettinen_1999}, which consists of weighting the objective functions into a single objective~$\sum_{i=1}^{q} \lambda_i f_i(x)$, where~$\lambda_i$ are non-negative weights.
The resulting optimization problem is
\begin{equation}\label{prob:multiobj_weightedsum}
    \begin{aligned}
        \min_{x \in X} \ & \sum_{i=1}^{q} \lambda_i f_i(x). 
    \end{aligned}
\end{equation}

We now present some properties of the weighted-sum method. Let $\Lambda$ denote the simplex set, i.e.,
\begin{equation*}\label{eq:simplex_set}
\Lambda = \{\lambda \in \mathbb{R}^q ~|~ \sum_{i=1}^{q} \lambda_i = 1, \, \lambda_i \ge 0, \ \forall i \in \{1, \ldots, q\} \}.
\end{equation*}
It is well known
that when $X$ is a convex set and the objective functions~$f_1, \ldots, f_q$ are convex, then~$x_* \in X_{wE}$ if and only if there exists some weight vector~$\lambda \in \Lambda$ such that~$x_*$ is an optimal solution to problem~\eqref{prob:multiobj_weightedsum}~\citep[Propositions 3.9 and 3.10]{MEhrgott_2005}. \tcb{Under the same convexity assumptions, $x_* \in X_{pE}$ if and only if there exists some weight vector~$\lambda \in \Lambda$, with~$\lambda_i > 0$ for all~$i \in \{1,\ldots,q\}$, such that~$x_*$ is an optimal solution to problem~\eqref{prob:multiobj_weightedsum}~\citep[Theorem 3.15]{MEhrgott_2005}.} Moreover, if~$X$ is a convex set and the objective functions~$f_1, \ldots, f_q$ are strictly convex, then~$X_{wE} = X_E$ (the proof follows from basic facts and is omitted for brevity).

\vspace{-0.1cm}
\begin{wrapfigure}{r}{0.43\textwidth}
\centering
\begin{tikzpicture}[>=Stealth, thick, font=\small, scale=0.4]

  \draw[->] (0,0) -- (7.6,0) node[right] {$f_1$};
  \draw[->] (0,0) -- (0,6.2) node[left] {$f_2$};

  \draw[darkgreen!65!black, line width=1.6pt]
    (0.7,5.5) .. controls (0.7,2.2) and (1.6,0.6) .. (6.8,0.6);

  \node[darkgreen!45!black] at (5.2,1.3) {$Y_N$};

  \filldraw[fill=white, draw=darkgreen!65!black, line width=1.1pt]
    (0.7,5.5) circle (3pt);
  \filldraw[fill=white, draw=darkgreen!65!black, line width=1.1pt]
    (6.8,0.6) circle (3pt);

  \coordinate (Px)  at (1.18,2.66);   
  \coordinate (Pxh) at (2.41,1.37);   
  \filldraw (Px)  circle (2.4pt);
  \filldraw (Pxh) circle (2.4pt);
  \node[above right=1pt] at (Px)  {$F(x) = \left(f_1(x),f_2(x)\right)$};
  \node[above right=1pt] at (Pxh) {$F(\hat{x})=\left(f_1(\hat{x}),f_2(\hat{x})\right)$};

  \draw[gray!60, dashed, thin] (Px)  -- (0,2.66);
  \draw[gray!60, dashed, thin] (Px)  -- (1.18,0);
  \draw[gray!60, dashed, thin] (Pxh) -- (0,1.37);
  \draw[gray!60, dashed, thin] (Pxh) -- (2.41,0);

  \draw[darkgreen!60!black, decorate,
        decoration={brace, amplitude=5pt}]
    (-0.15,1.37) -- (-0.15,2.66);
  \node[darkgreen!60!black, left=7pt, align=center] at (-0.15,2.015)
    {$f_2(x)-f_2(\hat{x})$};

  \draw[blue!70!black, decorate,
        decoration={brace, amplitude=5pt, mirror}]
    (1.18,-0.15) -- (2.41,-0.15);
  \node[blue!70!black, below=7pt] at (1.795,-0.15)
    {$f_1(\hat{x})-f_1(x)$};

\end{tikzpicture}
\caption{\tcb{Illustration of Geoffrion's definition of proper efficiency in the
bi-objective case. Moving from~$F(\hat{x})$ to~$F(x)$
improves~$f_1$ by the amount~$f_1(\hat{x})-f_1(x)$ at the cost of
deteriorating~$f_2$ by the amount~$f_2(x)-f_2(\hat{x})$.}}
\label{fig:geoffrion}
\end{wrapfigure}

Finally, we briefly introduce the geometric concepts that we will use to describe the Pareto front~$Y_N$. \tcb{The Pareto front can be viewed as the portion of the boundary of the image set~$Y=F(X)$ consisting of the nondominated points of~$Y$.} Under the settings considered in this paper, in the bi-objective case ($q=2$), the front is a curve in~$\mathbb{R}^2$, \tcb{with its two endpoints given by the single-objective optimal nondominated points}. In the tri-objective case ($q=3$), the front is a two-dimensional surface in~$\mathbb{R}^3$ \tcb{whose edge curves meet at the single-objective optimal nondominated points. Figure~\ref{fig:pareto} illustrates these two cases.} 
Following the terminology of~\citet{KDeb_SGupta_2010}, \tcb{when discussing bi-objective and tri-objective problems, we informally refer to the intermediate region of the Pareto front as the middle portion of the front, away from its two endpoints (when~$q=2$) or its edge curves (when~$q=3$).}

\begin{figure}[ht]
\centering
\begin{minipage}[c]{0.50\textwidth}
  \centering
  \scalebox{0.645}{%
\begin{tikzpicture}[>=Stealth, thick, font=\Large]
  \draw[->] (0,0) -- (7.6,0) node[right] {$f_1$};
  \draw[->] (0,0) -- (0,6.2) node[above] {$f_2$};
  \def\Ax{0.55} \def\Ay{5.8}   
  \def\Bx{7.0}  \def\By{0.45}  
  \def\CAx{0.55} \def\CAy{0.9}
  \def\CBx{2.8}  \def\CBy{0.45}
  \begin{scope}[on background layer]
    \fill[gray!13]
      (\Ax,\Ay)
        .. controls (\CAx,\CAy) and (\CBx,\CBy) ..
      (\Bx,\By)
      -- (7.1,\By) -- (7.1,5.8) -- (0,5.8) -- (0,\Ay) -- cycle;
  \end{scope}
  \draw[darkgreen!65!black, line width=2pt]
    (\Ax,\Ay) .. controls (\CAx,\CAy) and (\CBx,\CBy) .. (\Bx,\By);
  \filldraw[red!75!black] (\Ax,\Ay) circle (3.8pt);
  \filldraw[red!75!black] (\Bx,\By) circle (3.8pt);
  \node[right=5pt, red!75!black, align=left, font=\large] at (\Ax,\Ay)
    {Single-objective optimal nondominated point};
  \node[above=5pt, red!75!black, align=center,font=\large] at (\Bx,\By)
    {Single-objective optimal nondominated point};
  \node[darkgreen!45!black, align=center,font=\large] at (4.2,2.6)
    {Intermediate region};
  \draw[->, darkgreen!45!black, thin] (2.8,2.4) -- (1.6,1.95);
  \node[black!65!black, rotate=-72] at (1.12,3.5) {$Y_N$};
  \node[gray!55!black] at (5.5,5.0) {$Y = F(X)$};
\end{tikzpicture}
}
  \subcaption{Example of a Pareto front for a bi-objective problem ($q=2$)}
  \label{fig:pareto2d}
\end{minipage}%
\hfill
\begin{minipage}[c]{0.50\textwidth}
  \centering
  \scalebox{0.515}{%
\begin{tikzpicture}[>=Stealth, font=\Large,
    x={( 0.7cm, 0.25cm)},  
    y={(-0.7cm, 0.25cm)},  
    z={( 0.0cm,-1.00cm)}]  

  \draw[->] (0,0,0) -- (6.3,0,0) node[right]      {$f_1$};
  \draw[->] (0,0,0) -- (0,6.0,0) node[left]       {$f_2$};
  \draw[->] (0,0,0) -- (0,0,5.8) node[below]      {$f_3$};


  \coordinate (S00) at (0.000,0.000,5.000);
  \coordinate (S01) at (0.000,0.000,5.000);
  \coordinate (S02) at (0.000,0.000,5.000);
  \coordinate (S03) at (0.000,0.000,5.000);
  \coordinate (S04) at (0.000,0.000,5.000);
  \coordinate (S10) at (1.913,0.000,4.619);
  \coordinate (S11) at (1.768,0.732,4.619);
  \coordinate (S12) at (1.353,1.353,4.619);
  \coordinate (S13) at (0.732,1.768,4.619);
  \coordinate (S14) at (0.000,1.913,4.619);
  \coordinate (S20) at (3.536,0.000,3.536);
  \coordinate (S21) at (3.266,1.353,3.536);
  \coordinate (S22) at (2.500,2.500,3.536);
  \coordinate (S23) at (1.353,3.266,3.536);
  \coordinate (S24) at (0.000,3.536,3.536);
  \coordinate (S30) at (4.619,0.000,1.913);
  \coordinate (S31) at (4.268,1.768,1.913);
  \coordinate (S32) at (3.266,3.266,1.913);
  \coordinate (S33) at (1.768,4.268,1.913);
  \coordinate (S34) at (0.000,4.619,1.913);
  \coordinate (S40) at (5.000,0.000,0.000);
  \coordinate (S41) at (4.619,1.913,0.000);
  \coordinate (S42) at (3.536,3.536,0.000);
  \coordinate (S43) at (1.913,4.619,0.000);
  \coordinate (S44) at (0.000,5.000,0.000);

  \fill[darkgreen!12!white,opacity=0.90] (S30)--(S40)--(S41)--(S31)--cycle;
  \fill[darkgreen!11!white,opacity=0.90] (S31)--(S41)--(S42)--(S32)--cycle;
  \fill[darkgreen!12!white,opacity=0.90] (S32)--(S42)--(S43)--(S33)--cycle;
  \fill[darkgreen!11!white,opacity=0.90] (S33)--(S43)--(S44)--(S34)--cycle;
  \fill[darkgreen!15!white,opacity=0.90] (S20)--(S30)--(S31)--(S21)--cycle;
  \fill[darkgreen!14!white,opacity=0.90] (S21)--(S31)--(S32)--(S22)--cycle;
  \fill[darkgreen!15!white,opacity=0.90] (S22)--(S32)--(S33)--(S23)--cycle;
  \fill[darkgreen!14!white,opacity=0.90] (S23)--(S33)--(S34)--(S24)--cycle;
  \fill[darkgreen!18!white,opacity=0.90] (S10)--(S20)--(S21)--(S11)--cycle;
  \fill[darkgreen!17!white,opacity=0.90] (S11)--(S21)--(S22)--(S12)--cycle;
  \fill[darkgreen!18!white,opacity=0.90] (S12)--(S22)--(S23)--(S13)--cycle;
  \fill[darkgreen!17!white,opacity=0.90] (S13)--(S23)--(S24)--(S14)--cycle;
  \fill[darkgreen!20!white,opacity=0.90] (S00)--(S10)--(S11)--cycle;
  \fill[darkgreen!20!white,opacity=0.90] (S00)--(S11)--(S12)--cycle;
  \fill[darkgreen!20!white,opacity=0.90] (S00)--(S12)--(S13)--cycle;
  \fill[darkgreen!20!white,opacity=0.90] (S00)--(S13)--(S14)--cycle;

  \draw[red!75!black, line width=2pt]
    (5.000,0.000,0.000) .. controls (5.000,0.657,0.000) and (4.871,1.307,0.000) .. (4.619,1.913,0.000) .. controls (4.368,2.520,0.000) and (4.000,3.071,0.000) .. (3.536,3.536,0.000) .. controls (3.071,4.000,0.000) and (2.520,4.368,0.000) .. (1.913,4.619,0.000) .. controls (1.307,4.871,0.000) and (0.657,5.000,0.000) .. (0.000,5.000,0.000);
  \draw[red!75!black, line width=2pt]
    (0.000,0.000,5.000) .. controls (0.657,0.000,5.000) and (1.307,0.000,4.871) .. (1.913,0.000,4.619) .. controls (2.520,0.000,4.368) and (3.071,0.000,4.000) .. (3.536,0.000,3.536) .. controls (4.000,0.000,3.071) and (4.368,0.000,2.520) .. (4.619,0.000,1.913) .. controls (4.871,0.000,1.307) and (5.000,0.000,0.657) .. (5.000,0.000,0.000);
  \draw[red!75!black, line width=2pt]
    (0.000,0.000,5.000) .. controls (0.000,0.657,5.000) and (0.000,1.307,4.871) .. (0.000,1.913,4.619) .. controls (0.000,2.520,4.368) and (0.000,3.071,4.000) .. (0.000,3.536,3.536) .. controls (0.000,4.000,3.071) and (0.000,4.368,2.520) .. (0.000,4.619,1.913) .. controls (0.000,4.871,1.307) and (0.000,5.000,0.657) .. (0.000,5.000,0.000);

  \filldraw[red!75!black] (S40) circle (2.5pt);  
  \filldraw[red!75!black] (S44) circle (2.5pt);  
  \filldraw[red!75!black] (S00) circle (2.5pt);  

  \node[red!75!black, right=5pt, below=2pt, align=center,font=\large] at (S40)
    {Single-objective optimal \\ nondominated point};
  \node[red!75!black, left=5pt, below=2pt, align=center,font=\large] at (S44)
    {Single-objective optimal \\nondominated point};
  \node[red!75!black, below=0pt, font=\large] at (S00)
    {Single-objective optimal nondominated point};
  \node[red!70!black, font=\large, above=3pt] at (S42)
    {Upper edge curve};
  \node[red!70!black, font=\large, right=3pt] at (S20)
    {Right edge curve};
  \node[red!70!black, font=\large, left=3pt] at (S24)
    {Left edge curve};
  \node[darkgreen!45!black, align=center, font=\large] at (5.0,0.3,2.15)
    {Intermediate region};
  \draw[->, darkgreen!45!black, thin] (4.5,0.5,2.2) -- (1.5,1.2,2.0);
  \node[black!65!black, font=\Large] at (2.8,2.0,3.8) {$Y_N$};

\end{tikzpicture}
}
  \subcaption{Example of a Pareto front for a tri-objective problem ($q=3$)}
  \label{fig:pareto3d}
\end{minipage}
\caption{%
  Pareto fronts~$Y_N$ for a minimization problem in the bi-objective
  ($q=2$, left) and tri-objective ($q=3$, right) cases.}
\label{fig:pareto}
\end{figure}

\section{Finding knee solutions using Pareto sensitivity}\label{subsec:point_based} 

\tcb{In this section, we first define knee solutions and knee nondominated points, and we then develop optimization formulations and computational procedures for identifying them. Our definition follows the verbal notion of a knee (i.e., a nondominated point where a small improvement in any objective leads to a large deterioration in at least one other objective) by quantifying the trade-offs between objectives incurred when moving away from a given nondominated point. Specifically, to quantify such trade-offs, we resort to Definition~\ref{def:geoffrion} in Section~\ref{sec:basic_results}, which states that a Pareto optimal solution is properly efficient in the sense of Geoffrion if there exists a constant~$M > 0$ such that the trade-off ratio in~\eqref{eq:geoffrion_def} remains uniformly bounded by~$M$. Proper efficiency requires only that these trade-offs remain bounded, whereas Definition~\ref{def:knee} below identifies a knee solution as an efficient solution for which the largest such trade-off across all pairs of objectives is minimized\footnote{\tcb{Geoffrion's definition requires that, for every~$i$ and every~$x \in X$ satisfying~$f_i(x) < f_i(\hat{x})$, there exists an index~$j$ such that~$f_j(x) > f_j(\hat{x})$ and the corresponding trade-off ratio is bounded by~$M$. Thus, the trade-off ratio is required to be bounded for at least one~$j$, rather than for all~$j$. In contrast, our definition requires the trade-off ratio to be bounded for every pair~$i,j$ with~$i \neq j$ and every~$x \in X$ satisfying~$f_i(x) < f_i(\hat{x})$ and~$f_j(x) > f_j(\hat{x})$. Our definition is therefore stronger, but it leads to tractable formulations. We do not distinguish between weaker and stronger notions of knee solutions and focus exclusively on the stronger notion introduced here.}}.
\begin{definition}\label{def:knee}
For an efficient solution~$\hat{x} \in X_E$, for each~$x \in X$ and ordered pair~$(i,j)$ with~$i \neq j$, define
\[
g_x^{ij}(\hat{x})
=
\begin{cases}
\dfrac{f_i(\hat{x})-f_i(x)}
{f_j(x)-f_j(\hat{x})},
&
f_i(x)<f_i(\hat{x}),\;
f_j(x)>f_j(\hat{x}),\\[2mm]
0, & \text{otherwise}.
\end{cases}
\]
We define~$
M(\hat{x})
=
\sup_{x\in X}
\sup_{i\neq j}
g_x^{ij}(\hat{x})$.
An efficient solution~$\hat{x}_*$ minimizing~$M(\hat{x})$ over the set of efficient solutions~$X_E$ is called a \emph{knee solution}. The associated nondominated point~$F(\hat{x}_*) \in Y_N$ is called a \emph{knee nondominated point}.
\end{definition}}
\tcb{Since~$M(\hat{x})=+\infty$ for efficient solutions that are not properly efficient, any knee solution is properly efficient whenever some efficient solution has~$M(\hat{x})<\infty$. Remark~\ref{remark:existence_knee} establishes existence of a knee solution.
\begin{remark}\label{remark:existence_knee}
The strict-convexity condition, together with convexity of~$X$ in the constrained case, in Assumption~\ref{ass:convexity} below implies that the weakly efficient set~$X_{wE}$ and the efficient set~$X_E$ coincide, i.e.,~$X_{wE} = X_E$, and the weighted-sum characterization of weakly efficient solutions yields~$
X_E=\{x(\lambda):\lambda\in\Lambda\}$ (see Section~\ref{sec:basic_results}). Under Assumption~\ref{ass:IFT} in the unconstrained case (Subsection~\ref{subsec:pareto_sensitivity_calculation}) or Assumption~\ref{ass:LL_assumptions_constr} in the constrained case (Subsection~\ref{subsec:pareto_sensitivity_calculation_constr}), the solution mapping~$x(\cdot)$ is continuous. Since~$\Lambda$ is compact, it follows that~$X_E$ is compact.
Moreover,~$M(\cdot)$ is lower semicontinuous on~$X_E$. Indeed, for every fixed~$x\in X$ and ordered pair~$(i,j)$ with~$i\neq j$, the function~$g_x^{ij}(\cdot)$ is lower semicontinuous. In particular, when the denominator approaches zero from the positive side while the numerator remains positive, the ratio tends to~$+\infty$. Since~$M(\hat{x})$ is the pointwise supremum of the functions~$g_x^{ij}(\hat{x})$ over all~$x\in X$ and all ordered pairs~$(i,j)$ with~$i\neq j$,~$M(\cdot)$ is lower semicontinuous. Therefore, it attains its minimum over the compact set~$X_E$, and a knee solution exists.
\hfill$\square$
\end{remark}}

\tcb{Definition~\ref{def:knee} is not directly amenable to computation, since evaluating~$M(\hat{x})$ requires considering all feasible solutions~$x \in X$ for each~$\hat{x} \in X_{E}$. We therefore derive a local directional version of the trade-off ratio in~\eqref{eq:geoffrion_def}, restricted to directions tangent to the efficient set. Let~$\hat{x}$ be an efficient solution and let~$T_{X_E}(\hat{x})$ denote the tangent cone to~$X_E$ at~$\hat{x}$, which consists of all directions that can be approached by sequences of efficient solutions converging to~$\hat{x}$, as defined in~\citet[Definition~12.2]{JNocedal_SJWright_2006}.
For any distinct pair of objectives~$\{i,j\}\subseteq\{1,\ldots,q\}$, consider the set
\[
P_\varepsilon^{ij}(\hat{x})
\; = \;
\left\{
p\in T_{X_E}(\hat{x})~:~
\|p\|=1,\;
\nabla f_i(\hat{x})^\top p\leq-\varepsilon,\;
\nabla f_j(\hat{x})^\top p\geq\varepsilon
\right\},
\]
where~$\varepsilon>0$ ensures that the magnitudes of both directional changes are bounded away from zero.
For any~$p\in T_{X_E}(\hat{x})$, there exist sequences~$\{x_k\}\subseteq X_E$ and~$\{t_k\}$ with~$t_k\downarrow0$ such that~$x_k=\hat{x}+t_kp+o(t_k)$. A first-order Taylor expansion therefore gives~$f_i(\hat{x})-f_i(x_k)=-t_k\nabla f_i(\hat{x})^\top p+o(t_k)$ and~$f_j(x_k)-f_j(\hat{x})=t_k\nabla f_j(\hat{x})^\top p+o(t_k)$. Thus, as~$k\to\infty$, the ratio in~\eqref{eq:geoffrion_def} converges to the local trade-off ratio~$(-\nabla f_i(\hat{x})^\top p)/(\nabla f_j(\hat{x})^\top p).$
We can then define the \emph{maximal-change function}~(MCF) at~$\hat{x}$ as
\begin{equation}\label{eq:geoffrion_aux}
\MCF(\hat{x}) \; = \; \max_{i\neq j}
\sup_{p\in P_\varepsilon^{ij}(\hat{x})}
\frac{
-\nabla f_i(\hat{x})^\top p
}{
\nabla f_j(\hat{x})^\top p
},
\end{equation}
which measures the largest local trade-off between any pair of objectives over directions that improve one objective and deteriorate the other. 
Motivated by Definition~\ref{def:knee}, we therefore seek a knee solution by minimizing~\eqref{eq:geoffrion_aux} over the efficient set~$X_E$.
Directly solving this problem over~$X_E$ would require characterizing the efficient set, and hence the Pareto front, in advance. 
We instead represent~\eqref{eq:geoffrion_aux} in terms of the weight vector~$\lambda \in \Lambda$ used in the weighted-sum problem~\eqref{prob:multiobj_weightedsum}, where each~$\lambda$ yields a weakly efficient solution~$x(\lambda)$.
To ensure that~$x(\lambda)$ is efficient, rather than only weakly efficient, we require Assumption~\ref{ass:convexity}.
\begin{assumption}\label{ass:convexity}
The objective functions~$f_1,\ldots,f_q$ are strictly convex. In the constrained case~($X \subset \mathbb{R}^n$), the feasible set~$X$ is convex.
\end{assumption}
Under Assumption~\ref{ass:convexity}, we have~$X_{wE} = X_E$, and the weighted-sum characterization yields~$X_E=\{x(\lambda):\lambda\in\Lambda\}$ (see Section~\ref{sec:basic_results}).
Let~$T_\Lambda(\lambda)
=
\{
d\in\mathbb{R}^q:
\sum_{i=1}^q d_i=0,\;
d_i\geq0
\ \text{for all }i\in\{1, \ldots,q\}$ such that $\lambda_i=0
\}$ denote the tangent cone to~$\Lambda$ at~$\lambda$. A direction~$d\in T_\Lambda(\lambda)$ induces the direction~$\nabla x(\lambda)^{\top} d$ in the decision space. Therefore, expressing~\eqref{eq:geoffrion_aux} in terms of variations in~$\lambda$ generally restricts the directional search to the directions generated by perturbations of the scalarization weights. 
When
\begin{equation}\label{eq:tangent_coverage}
\nabla x(\lambda)^{\top} T_\Lambda(\lambda)
\; = \;
T_{X_E}(x(\lambda)),
\end{equation}
this restriction is without loss of generality, and the scalarization-based formulation is equivalent to~\eqref{eq:geoffrion_aux}.
Condition~\eqref{eq:tangent_coverage} states that perturbations of the scalarization weights generate all tangent directions to the efficient set at~$x(\lambda)$.
Under the assumptions needed for the application of the implicit function theorem, collected in Assumption~\ref{ass:IFT} in Subsection~\ref{subsec:pareto_sensitivity_calculation} for the unconstrained case and in Assumption~\ref{ass:LL_assumptions_constr} in Subsection~\ref{subsec:pareto_sensitivity_calculation_constr} for the constrained case, the solution mapping~$x(\lambda)$ is unique and differentiable, so that~$\nabla x(\lambda)$ and the sensitivities~$\nabla_{\lambda}(f_i(x(\lambda)))$,~$i\in\{1,\ldots,q\}$, are well defined.} \tcb{The following example illustrates a setting in which condition~\eqref{eq:tangent_coverage} holds.} 
\tcb{\begin{example}\label{example:tangent_space_condition}
Consider the bi-objective unconstrained problem~$f_1(x)=x^2$, $f_2(x)=(x-1)^2$, $x\in\mathbb{R}$.
For~$\lambda=(\lambda_1,\lambda_2)\in\Lambda$, the weighted-sum problem is~$\min_{x\in\mathbb{R}}
\lambda_1x^2+\lambda_2(x-1)^2$, whose unique solution is~$x(\lambda) = \lambda_2/(\lambda_1+\lambda_2)$.
Since~$\lambda_1+\lambda_2=1$ on~$\Lambda$, we have~$x(\lambda)=\lambda_2$, and the efficient set is~$X_E=[0,1]$. Consider, for instance,~$\lambda=(0.4,0.6)$, for which~$x(\lambda)=0.6$. The Jacobian of the solution mapping with respect to the weight vector is
\[
\nabla x(\lambda)^{\top}
\; = \;
\left(
-\frac{\lambda_2}{(\lambda_1+\lambda_2)^2},
\frac{\lambda_1}{(\lambda_1+\lambda_2)^2}
\right)
\; = \;
(-0.6,0.4).
\]
Since~$\lambda$ belongs to the relative interior of~$\Lambda$, we have~$T_\Lambda(\lambda) = \{(a,-a):a\in\mathbb{R}\}$.
Hence, for any~$(a,-a)\in T_\Lambda(\lambda)$, we have~$\nabla x(\lambda)^{\top} (a,-a)^{\top} = -a$,
and therefore~$\nabla x(\lambda)^{\top}T_\Lambda(\lambda)=\mathbb{R}$.
Since~$x(\lambda)=0.6$ is an interior point of~$X_E=[0,1]$, we also have~$T_{X_E}(x(\lambda))=\mathbb{R}$.
Thus,~\eqref{eq:tangent_coverage} holds at~$\lambda=(0.4,0.6)$.
\hfill$\square$
\end{example}}

\tcb{We now reformulate the~MCF~\eqref{eq:geoffrion_aux} in the weight vector space. To this end, for any distinct pair of objectives~$\{i,j\}$, we define
\[
D_\varepsilon^{ij}(\lambda)
\; = \;
\left\{
\begin{array}{l|l}
d\in T_\Lambda(\lambda)
&
\begin{array}{l}
\|d\|=1,\quad \nabla x(\lambda)^{\top}d\neq0,\\[1mm]
\nabla_\lambda(f_i(x(\lambda)))^\top d
\leq
-\varepsilon\|\nabla x(\lambda)^{\top}d\|,\\[1mm]
\nabla_\lambda(f_j(x(\lambda)))^\top d
\geq
\varepsilon\|\nabla x(\lambda)^{\top}d\|
\end{array}
\end{array}
\right\}.
\]
We first show that every direction~$d\in D_\varepsilon^{ij}(\lambda)$ induces a direction in~$P_\varepsilon^{ij}(x(\lambda))$. For any~$d\in D_\varepsilon^{ij}(\lambda)$, define the corresponding unit direction by~$p=(\nabla x(\lambda)^{\top}d)/\|\nabla x(\lambda)^{\top}d\|$. Since~$x(\lambda)\in X_E$ for every~$\lambda\in\Lambda$, any perturbation of~$\lambda$ along a direction~$d\in T_\Lambda(\lambda)$ generates points in~$X_E$. By differentiability of~$x(\cdot)$, the corresponding first-order direction is~$\nabla x(\lambda)^{\top}d$, and hence~$\nabla x(\lambda)^{\top}d\in T_{X_E}(x(\lambda))$. Thus,~$p\in T_{X_E}(x(\lambda))$ and~$\|p\|=1$. Moreover, dividing the last two inequalities in the definition of~$D_\varepsilon^{ij}(\lambda)$ by~$\|\nabla x(\lambda)^{\top}d\|$ and using the chain rule~$\nabla_\lambda(f_i(x(\lambda)))^\top d=\nabla f_i(x(\lambda))^\top\nabla x(\lambda)^{\top}d$ yields~$\nabla f_i(x(\lambda))^\top p\leq-\varepsilon$ and~$\nabla f_j(x(\lambda))^\top p\geq\varepsilon$. Therefore,~$p\in P_\varepsilon^{ij}(x(\lambda))$. The chain rule also gives
\[
\frac{
-\nabla_\lambda(f_i(x(\lambda)))^\top d
}{
\nabla_\lambda(f_j(x(\lambda)))^\top d
}
=
\frac{
-\nabla f_i(x(\lambda))^\top p
}{
\nabla f_j(x(\lambda))^\top p
},
\]
so the corresponding trade-off ratio in~\eqref{eq:geoffrion_aux} is preserved.
When condition~\eqref{eq:tangent_coverage} holds, every direction in~$P_\varepsilon^{ij}(x(\lambda))$ can also be generated by a direction in~$D_\varepsilon^{ij}(\lambda)$. Indeed, let~$\tilde p\in P_\varepsilon^{ij}(x(\lambda))$. Since~$\tilde p\in T_{X_E}(x(\lambda))$, condition~\eqref{eq:tangent_coverage} guarantees the existence of~$\bar d\in T_\Lambda(\lambda)$ such that~$\nabla x(\lambda)^{\top}\bar d=\tilde p$. Setting~$d=\bar d/\|\bar d\|$ gives~$d\in T_\Lambda(\lambda)$,~$\|d\|=1$, and~$\frac{\nabla x(\lambda)^{\top}d}
{\|\nabla x(\lambda)^{\top}d\|}
= \tilde p$.
Since~$\tilde p\in P_\varepsilon^{ij}(x(\lambda))$, multiplying the inequalities~$\nabla f_i(x(\lambda))^\top\tilde p\leq-\varepsilon$ and~$\nabla f_j(x(\lambda))^\top\tilde p\geq\varepsilon$ by~$\|\nabla x(\lambda)^{\top}d\|$ and using the chain rule yields~$
\nabla_\lambda(f_i(x(\lambda)))^\top d
\leq
-\varepsilon\|\nabla x(\lambda)^{\top}d\|$ and~$\nabla_\lambda(f_j(x(\lambda)))^\top d
\geq
\varepsilon\|\nabla x(\lambda)^{\top}d\|$,
so that~$d\in D_\varepsilon^{ij}(\lambda)$. Therefore, under~\eqref{eq:tangent_coverage},
\[
P_\varepsilon^{ij}(x(\lambda))
=
\left\{
\frac{\nabla x(\lambda)^{\top}d}
{\|\nabla x(\lambda)^{\top}d\|}
:
d\in D_\varepsilon^{ij}(\lambda)
\right\},
\]
and the corresponding trade-off ratios in~\eqref{eq:geoffrion_aux} are preserved.}

\tcb{We can now define the~MCF in the weight vector space as
\begin{equation}\label{eq:MCF}
\MCF_{\varepsilon}(\lambda)
\; = \;
\max_{i\neq j}
\sup_{d \in D_\varepsilon^{ij}(\lambda)}
\frac{
-\nabla_\lambda(f_i(x(\lambda)))^\top d
}{
\nabla_\lambda(f_j(x(\lambda)))^\top d
}.
\end{equation}
\tcb{When condition~\eqref{eq:tangent_coverage} holds, the scalarization-based directions cover all tangent directions to the efficient set, and therefore~$
\MCF_{\varepsilon}(\lambda)
\equiv
\MCF_{\varepsilon}(x(\lambda))$.}
\tcb{A weight vector corresponding to a knee solution can then be identified by solving the following minimization problem:}
\begin{equation}\label{prob:point_based_formulation_exact}
\begin{alignedat}{2}
    \min_{\lambda \in \Lambda} \ 
    \MCF_{\varepsilon}(\lambda).
\end{alignedat}
\end{equation}
}

\tcb{Although~$\MCF_{\varepsilon}(\lambda)$ has a direct interpretation in terms of local trade-offs, its evaluation requires solving a directional optimization problem over the tangent cone at each~$\lambda$. This motivates the search for a computationally simpler proxy, which we now introduce by first considering the bi-objective case.}
When~$q=2$, a Pareto front can be \tcb{represented} in the objective space as either the curve~$f_2 = f_2(f_1)$ or the curve~$f_1 = f_1(f_2)$. Let~$df_2/df_1$ and~$df_1/df_2$ denote the derivatives of~$f_2$ with respect to~$f_1$ and vice versa. These derivatives represent the \tcb{slopes} of the tangent line to the corresponding curve at a given nondominated point. 
\tcb{Recall that we look for a knee solution among properly efficient points in the sense of Geoffrion, at which the trade-offs between objectives remain bounded. We assume that the derivatives~$df_2/df_1$ and~$df_1/df_2$ at the corresponding nondominated points exist, so that the Pareto front curve is locally differentiable without kinks or other singularities.
In~\citet[Section~6]{KDeb_SGupta_2010}, the authors note that, according to the verbal definition of a knee solution, knee solutions correspond to nondominated points at which the derivatives~$df_2/df_1$ and~$df_1/df_2$ are equal in magnitude. We therefore observe that, in the bi-objective case, a knee solution minimizes~$\max\{\left|df_2/df_1\right|,\left|df_1/df_2\right|\}$.
For~$q>2$, the Pareto front is a $(q-1)$-dimensional surface, and each objective~$f_i$ can be represented locally as a function of the remaining objectives. The role of~$df_2/df_1$ is then played by the partial derivatives~$\partial f_i/\partial f_j$, with~$j\neq i$. This suggests extending the bi-objective formulation by minimizing the maximum~$\left|\partial f_i/\partial f_j\right|$ over all pairs~$\{i,j\}\subseteq\{1,\ldots,q\}$, with~$i\neq j$, thereby seeking to levelize the magnitudes of the pairwise partial derivatives. However, this approach is computationally impractical because it requires constructing a representation of the entire Pareto front and then modeling each objective as a function of the remaining objectives in order to evaluate the partial derivatives~$\partial f_i/\partial f_j$.}

\tcb{To replicate the role of~$|\partial f_i/\partial f_j|$ without computing the entire Pareto front in advance, our \textit{snee} approach uses the sensitivities~$\nabla_\lambda (f_i(x(\lambda)))$ and~$\nabla_\lambda (f_j(x(\lambda)))$ of the $i$-th and $j$-th objectives to variations in~$\lambda$, for all pairs~$\{i,j\} \subseteq \{1,\ldots,q\}$ with~$i \neq j$, resulting in the \emph{gradient maximal-change function}~($\MCFg$):
\begin{equation}\label{eq:maximal_change_function}
\MCFg(\lambda) \; = \; \max_{\substack{\{i,j\} \subseteq \{1,\ldots,q\} \\ i \neq j}} 
\; \frac{\lVert \nabla_{\lambda} (f_i(x(\lambda))) \rVert}{\lVert \nabla_{\lambda} (f_j(x(\lambda))) \rVert}
\; = \; \frac{g_{\max}(\lambda)}{g_{\min}(\lambda)},
\end{equation}
where $g_{\max}(\lambda) := \max_{i \in \{1,\ldots,q\}} \|\nabla_{\lambda} (f_i(x(\lambda)))\|$ and $g_{\min}(\lambda) := \min_{i \in \{1,\ldots,q\}} \|\nabla_{\lambda} (f_i(x(\lambda)))\|$.
Unlike~$\MCF_{\varepsilon}(\lambda)$, the quantity~$\MCFg(\lambda)$ depends only on the norms of the sensitivities and does not require the tangent cone.}
\tcb{Therefore, it is considerably simpler to evaluate than~$\MCF_{\varepsilon}(\lambda)$. Furthermore, because~$\MCFg(\lambda)$ admits a closed-form expression, minimizing it over~$\Lambda$ is significantly easier than minimizing~$\MCF_{\varepsilon}(\lambda)$, which requires repeatedly solving a directional optimization problem. Thus, $\MCFg(\lambda)$ provides a computationally attractive indicator of the local trade-off constant in Definition~\ref{def:knee}.
To determine the weight vector corresponding to a knee solution, we can therefore solve the following minimization problem:
\begin{equation}\label{prob:point_based_formulation}
\begin{alignedat}{2}
    \min_{\lambda \in \Lambda} \ 
    \MCFg(\lambda).
\end{alignedat}
\end{equation}}
Intuitively, the minimization problem in~\eqref{prob:point_based_formulation} aims to make the norms of $\nabla_{\lambda} (f_i (x(\lambda)))$ and $\nabla_{\lambda} (f_j (x(\lambda)))$ as close as possible.
Numerically, to avoid division by zero, we recommend replacing
$\lVert \nabla_{\lambda} (f_j (x(\lambda))) \rVert$ in~\eqref{eq:maximal_change_function}
by $\max \{ \lVert \nabla_{\lambda} (f_j (x(\lambda))) \rVert, \text{eps}\}$,
where $\text{eps}$ represents the machine precision.  

\tcb{Our~\textit{snee} approach identifies a Pareto knee solution and its corresponding knee nondominated point through the procedure summarized in Algorithm~\ref{alg:snee} \tcb{when problem~\eqref{prob:point_based_formulation_exact} is solved, and identifies them approximately when problem~\eqref{prob:point_based_formulation} is used instead}. 
Problem~\eqref{prob:point_based_formulation_exact} requires the solution of a directional optimization problem over the tangent cone for each candidate~$\lambda$, whereas problem~\eqref{prob:point_based_formulation}} involves the minimization of a non-smooth and non-convex function whose subgradients involve third-order derivatives. \tcb{Since the dimension of both problems is equal to the number~$q$ of objective functions, which is typically small, they} can be efficiently solved by a derivative-free optimization~(DFO) algorithm~\citep{ARConn_KScheinberg_LNVicente_2009,CAudet_WHare_2017,JLarson_MMenickelly_SMWild_2019,ALCustodio_KScheinberg_LNVicente_2017}.
For instance, one may use popular~DFO algorithms\tcb{, such as NOMAD\footnote{\tcb{The NOMAD solver implements the Mesh Adaptive Direct Search~(MADS) algorithm, a derivative-free optimization method designed to handle nonlinear and nonsmooth constraints. It does not require derivatives and can incorporate constraints directly into the optimization problem, making it particularly suitable when additional constraints need to be imposed.}}~\citep{CAudet_SLDigabel_VRochonMontplaisir_CTribes_2022},} Nelder-Mead\footnote{The Nelder-Mead~(NM) algorithm uses a simplex of points to navigate the search space, adjusting the simplex through operations such as reflection, expansion, contraction, and shrinkage. It requires an initial starting point, but it typically involves fewer function evaluations than global methods like DIRECT, making it suitable for problems where computational efficiency is a priority.}~\citep{JANelder_RMead_1965}, \tcb{or} DIRECT\footnote{The DIRECT algorithm is a global optimization method that systematically divides the search space into hyper-rectangles, calculating objective function values to find potential global minima without requiring derivative information. It does not require an initial starting point and it may necessitate a substantial number of function evaluations.}~\citep{DRJones_etal_1993}. Note that DFO is a well-studied class of optimization methods. In the literature, it is also referred to as zeroth-order optimization~\citep{YNesterov_VSpokoiny_2017} or black-box optimization~\citep{CAudet_WHare_2017}.

     \begin{algorithm}[htbp]
	\caption{Pareto Knee Solution Search}\label{alg:snee}
	\begin{algorithmic}[1]
		\medskip
		\item[] {\bf Input:} Initial weight vector~$\lambda_0$.
		\smallskip
        
    \item[] {\bf Main step.} Starting from~$\lambda_0$, solve \tcb{either problem~\eqref{prob:point_based_formulation_exact} or~\eqref{prob:point_based_formulation}} using a~DFO algorithm to obtain~$\lambda_*$. 
    \tcb{For a given~$\lambda \in \Lambda$, evaluate the objective function of the selected formulation as follows:}
    
    \item[] \quad\quad {\bf Step~A.} Given~$\lambda \in \Lambda$, apply the weighted-sum method with weight vector~$\lambda$ to obtain~$x(\lambda)$
    \item[] \quad\quad by solving problem~\eqref{prob:multiobj_weightedsum}.
    
    \item[] \quad\quad {\bf Step~B.} Given~$x(\lambda)$ from \textbf{Step~A}, compute the matrix~$\nabla_{\lambda}(F(x(\lambda)))$ using~\eqref{eq:adjoint} in Subsection~\ref{subsec:pareto_sensitivity_calculation}
    \item[] \quad\quad for the unconstrained case~($X = \mathbb{R}^n$) or~\eqref{eq:adjoint_2} in Subsection~\ref{subsec:pareto_sensitivity_calculation_constr} for the constrained case~($X \subset \mathbb{R}^n$).
    
    \item[] \quad\quad {\bf Step~C.} Evaluate \tcb{the~MCF~\eqref{eq:MCF} if problem~\eqref{prob:point_based_formulation_exact} is selected, or the~gMCF~\eqref{eq:maximal_change_function} if problem~\eqref{prob:point_based_formulation} is selected,}  
    \item[] \quad\quad where the gradients~$\nabla_{\lambda} (f_i (x(\lambda)))$ and~$\nabla_{\lambda} (f_j (x(\lambda)))$ correspond to the~$i$-th and~$j$-th columns,
    \item[] \quad\quad respectively, of matrix~$\nabla_{\lambda}(F(x(\lambda)))$ from~\textbf{Step~B}.

    \smallskip
    \item[] {\bf Output:} Optimal weight vector~$\lambda_*$ and corresponding solution~$x(\lambda_*)$.

		\par\bigskip\noindent
    	\end{algorithmic}
    \end{algorithm}

To evaluate the \tcb{MCF~\eqref{eq:MCF} or the~gMCF~\eqref{eq:maximal_change_function}} at each iteration of the selected~DFO algorithm for a given~$\lambda \in \Lambda$, one first applies the weighted-sum method to obtain~$x(\lambda)$ (\textbf{Step~A}), and then computes the matrix~$\nabla_{\lambda} (F(x(\lambda)))$ (\textbf{Step~B}). Such a matrix is formally derived in~\eqref{eq:adjoint} in Subsection~\ref{subsec:pareto_sensitivity_calculation} for the unconstrained case~($X = \mathbb{R}^n$) and~\eqref{eq:adjoint_2} in Subsection~\ref{subsec:pareto_sensitivity_calculation_constr} for the constrained case~($X \subset \mathbb{R}^n$).
\tcb{The~MCF~\eqref{eq:MCF}, when problem~\eqref{prob:point_based_formulation_exact} is selected, and the~gMCF~\eqref{eq:maximal_change_function}, when problem~\eqref{prob:point_based_formulation} is selected, are} then evaluated by noting that the gradients~$\nabla_{\lambda} (f_i (x(\lambda)))$ and~$\nabla_{\lambda} (f_j (x(\lambda)))$ correspond to the columns of such a matrix (\textbf{Step~C}).

\tcb{Although using a~DFO algorithm avoids the need for derivatives of the objective functions when solving problems~\eqref{prob:point_based_formulation_exact} and~\eqref{prob:point_based_formulation}, the evaluation of both the~MCF and the~gMCF still requires computing the gradients in~\eqref{eq:MCF} and~\eqref{eq:maximal_change_function}.} For large-scale machine learning applications, this may require using stochastic gradients (see Remark~\ref{remark:dfo} below).

\begin{remark}\label{remark:dfo}
Even when using a~DFO algorithm, \tcb{we need to efficiently compute the gradients required in~\eqref{eq:MCF} and~\eqref{eq:maximal_change_function} to evaluate the~MCF and~gMCF, respectively.} 
In machine learning applications, where computing a gradient involves summing the loss over the entire dataset, researchers typically use a stochastic gradient computed from randomly sampled \tcb{subsets} of data~\citep{KLChung_1954,HRobbins_SMonro_1951,JSacks_1958,LBottou_FECurtis_JNocedal_2018}. 
When the gradients~$\nabla_{\lambda} (f_i (x(\lambda)))$ and~$\nabla_{\lambda} (f_j (x(\lambda)))$ \tcb{in~\eqref{eq:MCF} and~\eqref{eq:maximal_change_function} are approximated using stochastic gradients, the resulting evaluations of the~MCF and~gMCF become stochastic. One possible approach is to formulate the objectives of problems~\eqref{prob:point_based_formulation_exact} and~\eqref{prob:point_based_formulation} in terms of expectations with respect to the probability distribution used to sample the data, which introduces} additional challenges that are beyond the scope of this paper. 
Moreover, note that many popular loss functions in machine learning are nonconvex (see~\citet[Subsection~4.3]{LBottou_FECurtis_JNocedal_2018}) and therefore do not satisfy Assumptions~\ref{ass:IFT} and~\ref{ass:LL_assumptions_constr}, which will be introduced in Subsections~\ref{subsec:pareto_sensitivity_calculation} and~\ref{subsec:pareto_sensitivity_calculation_constr}. Finally, we note that our approach requires solving the weighted-sum problem~\eqref{prob:multiobj_weightedsum} for each~$\lambda$ generated by a~DFO algorithm (see~\textbf{Step~A} of Algorithm~\ref{alg:snee}). This may be computationally expensive when the weighted-sum problem is difficult to solve. In such cases, the normal boundary intersection method can be preferable because it requires solving only a single nonlinear constrained optimization problem, provided that \tcb{single-objective optimal nondominated points, or estimates thereof,} are available. 
\hfill$\square$
\end{remark}

\section{Most-changing Pareto sub-fronts around nondominated points}\label{sec:neighborhoods} 


\tcb{Having identified a \textit{snee} solution~$x(\lambda_*)$ via Algorithm~\ref{alg:snee} in Section~\ref{subsec:point_based}, a natural question is how the nondominated set behaves locally around the corresponding nondominated point~$F(x(\lambda_*))$. In this section, we introduce the notion of most-changing Pareto sub-fronts, which capture the directions of greatest variation in the objectives around a given nondominated point. Specifically, given a Pareto optimal solution~$x_c$ of interest, we seek a neighborhood of weight vectors~$\lambda$ centered at~$\lambda_c$, with~$x_c = x(\lambda_c)$, where~$F(x(\lambda))$ changes the most. The resulting Pareto sub-front consists of the nondominated points corresponding to the Pareto optimal solutions identified by weight vectors in that neighborhood. While our primary motivation is to apply this construction to the \textit{snee} solution, by setting~$\lambda_c = \lambda_*$, so as to visually highlight the directions of maximum change around the knee nondominated point~$F(x(\lambda_*))$, the construction applies to any nondominated point~$F(x(\lambda_c))$ for a given~$\lambda_c \in \Lambda$.} Next, we detail how to construct such a Pareto sub-front.

Given a neighborhood size~$\alpha > 0$, a sub-front of most change is identified using the following neighborhood of weight vectors
\begin{equation} \label{eq:Ninverse}
\mathcal{E}_\alpha (\lambda_c) \; = \; \{ \lambda \in \mathbb{R}^q \; | \;
\| \nabla_{\lambda} (F(x(\lambda_c)))^\dagger (\lambda - \lambda_c) \| \leq \alpha \},
\end{equation}
where~$\nabla_{\lambda} (F(x(\lambda_c)))$ is the transpose of the Jacobian matrix of~$F(x(\cdot))$ evaluated at~$\lambda_c$,
and $\nabla_{\lambda} (F(x(\lambda_c)))^\dagger$ is its pseudo-inverse.  
Neighborhood~(\ref{eq:Ninverse}) can be motivated in two different, related ways. First note that when~$n=q=1$, we obtain~$| \lambda - \lambda_c| \leq \alpha \cdot \frac{d}{d\lambda} (f(x(\lambda_c)))$, where~$\frac{d}{d\lambda} (f(x(\lambda_c)))$ denotes the derivative of~$f(x(\cdot))$ evaluated at~$\lambda_c$, and we can see that we are basically considering larger intervals centered at~$\lambda_c$ with amplitudes increasing with the slope. This observation provides a geometric interpretation of~\eqref{eq:Ninverse}, as the resulting neighborhood is an ellipsoid with its major axis aligned with the direction of maximum change. For the remainder of the paper, we refer to neighborhood~(\ref{eq:Ninverse}) as the ellipsoidal neighborhood.

One can also introduce a motivating argument based on the general notion of steepest ascent/descent. We are interested in the weight vectors~$\lambda$ that lead to a meaningful change in~$\| F(x(\lambda)) -  F(x(\lambda_c)) \|$.
Given the expansion
\[
\| F(x(\lambda)) -  F(x(\lambda_c)) \| \; = \;  \| \nabla_{\lambda} (F(x(\lambda_c))) \Delta \lambda \| + \mathcal{O}(\| \Delta \lambda \|^2),
\]
we ask for a quadratic decrease in~$\| \lambda - \lambda_c \|$ and introduce the neighboorhood
\begin{equation} \label{eq:N}
\mathcal{E}_\beta (\lambda_c) \; = \; \{ \lambda \in \mathbb{R}^q \; | \;
 \| \nabla_{\lambda} (F(x(\lambda_c))) (\lambda - \lambda_c) \| \geq \beta  \|\lambda - \lambda_c \|^2 \},
\end{equation}
for some~$\beta > 0$.
Note that~$ \| \nabla_{\lambda} (F(x(\lambda_c))) (\lambda - \lambda_c) \| \geq \beta  \|\lambda - \lambda_c \|^2$ implies 
\[
\frac{\kappa(\nabla_{\lambda} (F(x(\lambda_c))))}{\beta} \; \geq \;  \| \nabla_{\lambda} (F(x(\lambda_c)))^\dagger (\lambda - \lambda_c) \|,
\]
where~$\kappa(\nabla_{\lambda} (F(x(\lambda_c))))= \| \nabla_{\lambda} (F(x(\lambda_c))) \| 
\| \nabla_{\lambda} (F(x(\lambda_c)))^\dagger \|$ measures the conditioning of \linebreak $\nabla_{\lambda} (F(x(\lambda_c)))$ beyond singularity. We conclude that for moderate~$\kappa(\nabla_{\lambda} (F(x(\lambda_c))))$, we recover~(\ref{eq:Ninverse}) \linebreak from~(\ref{eq:N}) by setting~$\alpha = \kappa(\nabla_{\lambda} (F(x(\lambda_c))))/\beta$. Geometrically, neighborhood~\eqref{eq:N} has a shape similar to a Cassini oval~\cite[p.~153]{JDLawrence_1972}, which is visible in some of the plots of this paper. For reference, a Cassini oval is the set of points in the plane such that the product of the distances to two fixed points (foci) is constant. Neighborhood~\eqref{eq:N} only roughly resembles this shape. 

Given a weight vector~$\lambda_c$, to evaluate the effectiveness of~a Pareto neighborhood~$\mathcal{E}(\lambda_c)$—either $\mathcal{E}_{\alpha}(\lambda_c)$ in~\eqref{eq:Ninverse} or~$\mathcal{E}_{\beta}(\lambda_c)$ in~\eqref{eq:N}—in identifying 
the greatest variation across the objectives on the Pareto front through the corresponding Pareto sub-front~$\{F(x(\lambda)) ~|~ \lambda \in \mathcal{E}(\lambda_c)\}$, we will use the most-changing \tcb{indicator}~\eqref{eq:metric} below. Such \tcb{an indicator} measures the extent of the Pareto front covered by a set of nondominated points relative to the full extent of the front. In the literature, such \tcb{an indicator} is referred to as the overall Pareto spread~\citep{JWu_SAzarm_2000,CAudet_JBigeon_etal_2021}. 
In our paper, we prefer not to use the term {\it spread} because it is not relevant to our context.
 The most-changing \tcb{indicator}~\eqref{eq:metric} evaluates the extent of the Pareto front covered by the nondominated points corresponding to the Pareto optimal solutions obtained by applying the weighted-sum method with weight vectors in the intersection between a Pareto neighborhood and the simplex set. In particular, the most-changing \tcb{indicator~(MCI)} is given by
\begin{equation}\label{eq:metric}
\MCI ~=~ \MCI( \mathcal{E}(\lambda_c) ) ~=~
\prod_{i=1}^{q} 
 \frac{|\max_{\lambda \in \mathcal{E}(\lambda_c) \cap \Lambda} f_i(x(\lambda)) - \min_{\lambda \in \mathcal{E}(\lambda_c) \cap \Lambda} f_i(x(\lambda))|}
 {|\max_{\lambda \in \Lambda} f_i(x(\lambda)) - \min_{\lambda \in \Lambda} f_i(x(\lambda))|}. 
\end{equation} 

The numerator in~\eqref{eq:metric} measures the maximum variation of the objective function~$f_i$ over the set of weights in~$\mathcal{E} (\lambda_c) \cap \Lambda$. The denominator measures the maximum variation of~$f_i$ over the set of weights in~$\Lambda$. The denominator can also be written as~$|F_i^I - F_i^N|$, where~$F^I$ and~$F^N$ are the ideal and nadir points introduced in Definition~\ref{def:ideal_nadir} in Section~\ref{sec:basic_results}. Note that metric~\eqref{eq:metric} lies between~0 and~1. Since such \tcb{an indicator} is computed relative to the full Pareto front, it is not affected by differences in the orders of magnitude of the objective functions. 
A higher value of \tcb{indicator}~\eqref{eq:metric} indicates that the corresponding Pareto sub-front exhibits greater variation across the objectives.

To determine the values of~$\alpha$ and~$\beta$ in neighborhoods~\eqref{eq:Ninverse} and~\eqref{eq:N}, one can use either a fixed stepsize or an adaptive rule. Using a fixed value requires fine-tuning and may be unreliable, as observed in the numerical experiments (see Remark~\ref{rem:MCFvsMCM} in Subsection~\ref{subsec:results_unc}). To apply the proposed adaptive rule, we consider a finite set of equidistant weight vectors corresponding to a fine-scale discretization of~$\Lambda$, i.e.,~$\Lambda_m = \{\lambda_1, \ldots, \lambda_m\} \subset \Lambda$. Specifically, for neighborhood~\eqref{eq:Ninverse}, we set~$\alpha$ to~$\gamma(1/m)\sum_{i=1}^{m} \Vert \nabla_{\lambda} (F(x( \lambda_c)))^\dagger (\lambda_i - \lambda_c)\Vert$, where~$\gamma > 0$. The value of~$\gamma$ should be chosen so that a reasonable proportion of the weight vectors in~$\Lambda_m$ lie within each neighborhood~$\mathcal{E}(\lambda_c)$. A similar adaptive rule can be applied for neighborhood~\eqref{eq:N}.

The next remark gives insights into how to determine a most-changing Pareto sub-front when multiple nondominated points are available.

\begin{remark}
If one wants to determine a most-changing Pareto sub-front around multiple nondominated points, there are several approaches one can use based on the above procedure for a single point. A first approach consists of computing a most-changing neighborhood using the centroid of the weight vectors that correspond to the given nondominated points.
A second approach involves determining the centroid of the given nondominated points and computing a most-changing neighborhood using a matrix that could be either the Jacobian at the centroid or a linear combination of the Jacobians at the given nondominated points. A third approach consists of computing a most-changing Pareto sub-front for all the given nondominated points, and then taking the union of such sub-fronts.
\hfill$\square$
\end{remark}

\section{The unconstrained case}\label{sec:unconstr}

Throughout this section, we focus on the case~$X = \mathbb{R}^n$.

\subsection{Pareto sensitivity calculation in the unconstrained case}\label{subsec:pareto_sensitivity_calculation}

In this subsection, given~$\lambda \in \Lambda$, we first characterize the Jacobian, or sensitivity matrix,~$\nabla x(\lambda)^\top$, and then calculate the gradient of~$f_i(x(\lambda))$ with respect to~$\lambda$, for all~$i \in \{1, \ldots, q\}$. Such a gradient allows us to derive the Jacobian matrix of~$F(x(\lambda)) = (f_1(x(\lambda)), \ldots, f_q(x(\lambda)))^\top$ with respect to~$\lambda$, which plays a crucial role in determining Pareto knee solutions in Section~\ref{subsec:point_based} and in computing Pareto sub-fronts in Section~\ref{sec:neighborhoods}. 

For any given~$\lambda \in \Lambda$, let us assume the existence of $x(\lambda)$ satisfying the first-order necessary optimality conditions
\begin{equation} \label{eq:1storderKKT}
\sum_{i=1}^{q} \lambda_i \nabla_x f_i(x(\lambda)) \; = \; 0.
\end{equation}
We will later collect all the assumptions required for the application of
the implicit function theorem at~$x(\lambda)$ (i.e., Assumption~\ref{ass:IFT} below), which ensures the uniqueness of $x(\lambda)$ and its continuous differentiability~\citep[p.~224, Theorem~9.28]{Ruding_1953}. To calculate the Jacobian~$\nabla x(\lambda)^{\top} \in \mathbb{R}^{n \times q}$, we take derivatives with respect to~$\lambda$ on both sides of the first-order necessary optimality conditions~\eqref{eq:1storderKKT}, obtaining
\begin{equation*}  
    (\nabla_x f_1(x(\lambda)), \ldots, \nabla_x f_q(x(\lambda)))^\top + \sum_{i=1}^{q} \lambda_i \nabla x(\lambda) \nabla_{xx}^2 f_i(x(\lambda)) = 0,
\end{equation*}
which, under \tcb{Assumption~\ref{ass:IFT} below}, yields
\begin{equation}\label{eq:04}
   \nabla x(\lambda) = - (\nabla_x f_1(x(\lambda)), \ldots, \nabla_x f_q(x(\lambda)))^\top \left(\sum_{i=1}^{q} \lambda_i \nabla_{xx}^2 f_i(x(\lambda))\right)^{-1}.
\end{equation}

By applying the chain rule and using~\eqref{eq:04}, the transpose of the Jacobian matrix of~$F(x(\cdot))$ at~$\lambda$ is given by
\begin{equation}\label{eq:adjoint}
\begin{alignedat}{2}
    \nabla_{\lambda} ( F(x(\lambda))) \; &= \; \nabla x(\lambda) \nabla_x F(x(\lambda)) \\ \; &= \;  - (\nabla_x f_1, \ldots, \nabla_x f_q)^\top \left(\sum_{i=1}^{q} \lambda_i \nabla_{xx}^2 f_i\right)^{-1} (\nabla_x f_1, \ldots, \nabla_x f_q),
\end{alignedat}
\end{equation}
where all the gradients~$\nabla_x f_i$ and Hessians~$\nabla_{xx}^2 f_i$ are evaluated at $x(\lambda)$, with~$i \in \{1, \ldots, q\}$. 
Recalling that~$F(x(\lambda)) = (f_1(x(\lambda)), \ldots, f_q(x(\lambda)))^\top$, the columns of~$\nabla_{\lambda}(F(x(\lambda)))$ correspond to the gradients of the functions~$f_i(x(\lambda))$, i.e.,~$\nabla_{\lambda} (F(x(\lambda))) = (\nabla_{\lambda}(f_1(x(\lambda))), \ldots, \nabla_{\lambda}(f_q(x(\lambda))))^\top$. Note that~$\nabla_{\lambda}(F(x(\lambda)))$ is symmetric. It is also singular due to the fact that~(\ref{eq:1storderKKT}) implies the linear dependence of the gradients at $x(\lambda_c)$.

To apply this derivation when~$X = \mathbb{R}^n$, we \tcb{require} the matrix resulting from the convex linear combination of the individual Hessian matrices to be non-singular at~$x(\lambda)$, as evident from the expression for~$\nabla_{\lambda}( F(x(\lambda)))$ in~\eqref{eq:adjoint}.
We collect all the required working assumptions in Assumption~\ref{ass:IFT} below. Note that a sufficient condition for Assumption~\ref{ass:IFT} to hold is that the objective functions~$f_1, \ldots, f_q$ are strongly convex~\cite[p.~74, Definition~2.1.3]{YNesterov_2018}. In this case, for any $\lambda \in \Lambda$, the weighted-sum function~$\sum_{i=1}^{q} \lambda_i f_i(x)$ has a unique minimizer~$x(\lambda)$ and the Hessian~$\sum_{i=1}^{q} \lambda_i \nabla_{xx}^2 f_i(x(\lambda))$ is positive definite and therefore nonsingular~\cite[p.~75, Theorem~2.1.11]{YNesterov_2018}.

\begin{assumption}[Existence of Pareto optimal solutions (LL unconstrained case)]\label{ass:IFT}
Let $X = \mathbb{R}^n$. The objective functions~$f_1, \ldots, f_q$ are twice continuously differentiable. For any~$\lambda \in \Lambda$, there exists a point~$x(\lambda)$ such that $\sum_{i=1}^{q} \lambda_i$  $\nabla_x f_i(x(\lambda)) = 0$ and the convex linear combination~$\sum_{i=1}^{q} \lambda_i \nabla_{xx}^2 f_i(x(\lambda))$ of the individual Hessian matrices is non-singular.
\end{assumption}

Note that if we are given a Pareto optimal solution~$x_c$, and the objective functions meet the conditions of Assumption~\ref{ass:IFT}, one can calculate a~$\lambda_c$ such that the first-order necessary conditions~\eqref{eq:1storderKKT} are satisfied, and then proceed from there to calculate a sub-front of most change. Note also that such a~$\lambda_c$ may not be unique, which could potentially lead to different sub-fronts of most change (see Example~\ref{ex:nonunique_lambda} below).

\begin{example}\label{ex:nonunique_lambda}
Consider $n=1$, $q=3$ with $f_1(x)=(x-1)^2$, $f_2(x)=x^2$,
$f_3(x)=(x+1)^2$.
Each~$f_i$ has Hessian~$f_i''(x)=2>0$, so all three objectives are
strongly convex and Assumption~1 holds.
The first-order necessary optimal conditions~\eqref{eq:1storderKKT} reduce to
$\sum_{i=1}^3 \lambda_i f_i'(x)=0$, which gives $x(\lambda)=\lambda_1-\lambda_3$.
Now, choose arbitrarily~$x_c=0.2$. The condition~$\lambda_1-\lambda_3=0.2$ together with~$\lambda_1+\lambda_2+\lambda_3=1$ and~$\lambda_i\ge 0$ admits
the one-parameter family
$\lambda_c(t)=(0.2+t,\;0.8-2t,\;t)^\top$, $t\in[0,0.4]$,
so infinitely many weight vectors correspond to the same~$x_c$.
Three members of this family are
$\lambda_c^{(1)}=(0.2,\,0.8,\,0)^\top$,
$\lambda_c^{(2)}=(0.4,\,0.4,\,0.2)^\top$, and~$\lambda_c^{(3)}=(0.6,\,0,\,0.4)^\top$.
This illustrates the non-uniqueness of $\lambda_c$ when~$n<q-1$
(here $n=1$ and $q-1=2$).
Although all three weight vectors yield $x(\lambda_c)=0.2$ and produce the same $\nabla_{\lambda}(F(x(\lambda_c)))$ (since it depends on~$\lambda$ only through~$x(\lambda_c)$), neighborhoods~\eqref{eq:Ninverse} and~\eqref{eq:N} are centered at
different vectors, so the corresponding Pareto sub-fronts will generally differ across choices of~$\lambda_c$.
\hfill$\square$
\end{example}


\subsection{MOO unconstrained test problems}\label{subsec:moo_test_problems}

Table~\ref{tab:test_prob} specifies the unconstrained test problems used for the experiments in Subsections~\ref{subsec:results_unc} and~\ref{subsec:results_unconstr}, along with the number of variables and objectives (i.e.,~$n$ and~$q$, respectively).  
All of the problems have strictly convex objective functions.
\begin{table}
    \footnotesize
    \centering
    \begin{tabular}{ c|c|c|c|l } 
    Problem & $n$ & $q$ & Ref. & \multicolumn{1}{c}{Objective Functions}\\[4pt]
     \hline
    \multirow{3}{*}{ZLT1} & \multirow{3}{*}{3} & \multirow{3}{*}{3} & \multirow{3}{*}{\cite{SHuband_etal_2006}} & $f_1(x_1, x_2, x_3) = (x_1 - 1)^2 + x_2^2 + x_3^2$ \\
     &  &  & &  $f_2(x_1, x_2, x_3) = x_1^2 + (x_2 - 1)^2 + x_3^2$\\
     &  &  & &  $f_3(x_1, x_2, x_3) = x_1^2 + x_2^2 + (x_3 - 1)^2$\\[4pt]
      \hline
    \multirow{3}{*}{GRV1} & \multirow{3}{*}{$2\bar{n}$} & \multirow{3}{*}{3} &  & $f_1(x_1, x_2) = \frac{1}{2} x_1^{\top} H_1 x_1 + \frac{1}{2} x_2^{\top} H_2 x_2 + x_1^\top H_3 x_2 + a_1^\top x_1 + a_2^\top x_2$ \\
     &  &  & &  $f_2(x_1, x_2) = \frac{1}{2} x_1^{\top} H_3 x_1 + \frac{1}{2} x_2^{\top} H_4 x_2 + x_1^\top H_5 x_2 + a_3^\top x_1 + a_4^\top x_2$\\
     &  &  & &  $f_3(x_1, x_2) = \frac{1}{2} x_1^{\top} H_7 x_1 + \frac{1}{2} x_2^{\top} H_8 x_2 + x_1^\top H_9 x_2 + a_5^\top x_1 + a_6^\top x_2$\\[4pt]
     \hline
    \multirow{3}{*}{VFM1} & \multirow{3}{*}{$2$} & \multirow{3}{*}{3} & \multirow{3}{*}{\cite{SHuband_etal_2006}} & $f_1(x_1, x_2) = x_1^2 + (x_2 - 1)^2$ \\
     &  &  & &  $f_2(x_1, x_2) = x_1^2 + (x_2 + 1)^2 + 1$\\
     &  &  & &  $f_3(x_1, x_2) = (x_1 - 1)^2 + x_2^2 + 2$\\[4pt]
      \hline
    ZLT1$q$ & $\bar{n}$ & $\bar{q}$ & \cite{SHuband_etal_2006} & $f_j(x) = (x_j - 1)^2 + \sum_{1\le i \le \bar{n}, i \ne j} x_i^2, \quad j \in \{1,\ldots,\bar{q}\}$ \\[4pt]
    \hline
    \multirow{2}{*}{GRV2} & \multirow{2}{*}{$\bar{n}$} & \multirow{2}{*}{$2$} &  & $f_1(x) = \frac{1}{\bar{n}} \sum_{i=1}^{\bar{n}} x_i^2 + \frac{1}{2} \sum_{i=1}^{\bar{n}} x_i^4$\\
    &  &  & & $f_2(x) = \frac{1}{\bar{n}} \sum_{i=1}^{\bar{n}} (x_i - 2)^2 + \frac{1}{2} \sum_{i=1}^{\bar{n}} (x_i - 2)^4$\\
    \end{tabular}
    \caption{Unconstrained test problems~($\bar{n}$ and~$\bar{q}$ are arbitrary positive scalars).}\label{tab:test_prob}
\end{table}
In Problem~GRV1, $a_1$, $a_2$, and~$a_5$ were randomly generated according to a uniform distribution over~$[-5,0)$, resulting in~$a_1=-1.87$, $a_2=-4.75$, and~$a_5=-0.78$. Similarly, $a_3$, $a_4$, and~$a_6$ were randomly generated according to a uniform distribution over~$[0,5)$, resulting in~$a_3=3.66$, $a_4=2.99$, and~$a_6=0.78$. Additionally, we randomly generated three $\bar{n}\times\bar{n}$ symmetric positive definite matrices~$H^{(1)}$, $H^{(2)}$, and~$H^{(3)}$, and we set the matrices in the objective functions~$f_i$, with~$i \in \{1,2,3\}$, as follows:
\[
H^{(1)} = \begin{bmatrix}
  H_{1} & H_{3} \\
  H_{3}^{\top} & H_{2}
\end{bmatrix} = \begin{bmatrix}
  50.82 & -0.23 \\
  -0.23 & 10.57
\end{bmatrix}, \quad
H^{(2)} = \begin{bmatrix}
  H_{4} & H_{6} \\
  H_{6}^{\top} & H_{5}
\end{bmatrix} = \begin{bmatrix}
  38.25 & 12.19 \\
  12.19 & 6.53
\end{bmatrix},
\]
\[
H^{(3)} = \begin{bmatrix}
  H_{7} & H_{9} \\
  H_{9}^{\top} & H_{8}
\end{bmatrix} = \begin{bmatrix}
  45.10 & -9.55 \\
  -9.55 & 9.91
\end{bmatrix}.
\]

\subsection{Results for knee solutions using Pareto sensitivity in the unconstrained case}\label{subsec:results_unc} 

We consider the unconstrained problems from Table~\ref{tab:test_prob} in Subsection~\ref{subsec:moo_test_problems}. 
For each problem with~$n \le 3$ and~$q \le 3$, we graphically represent the weight vector, objective, and decision spaces.
To solve problem~\tcb{\eqref{prob:point_based_formulation_exact} and~\eqref{prob:point_based_formulation}} in the main step of Algorithm~\ref{alg:snee}, 
\tcb{we employ NOMAD's Mesh Adaptive Direct Search~(MADS) algorithm via its Python interface, PyNomadBBO~\citep{CAudet_SLDigabel_VRochonMontplaisir_CTribes_2022}, introduced in Section~\ref{subsec:point_based}, run with default parameters.
We also solved problem~\eqref{prob:point_based_formulation_exact} and~\eqref{prob:point_based_formulation} with the NM and DIRECT algorithms from the Python SciPy library~\citep{PVirtanen_etal_2020_Scipy}. Since the resulting knee solutions are nearly identical to those obtained with MADS (there are no noticeable differences in the plots), we report them in the Appendix rather than here.
At each iteration of all three methods, feasibility is ensured by computing orthogonal projections onto~$\Lambda$.}
To perform \textbf{Step~A} and obtain the Pareto optimal solution associated with a given weight vector, we solve the weighted-sum problem~\eqref{prob:multiobj_weightedsum} using the BFGS algorithm implementation available in the Python SciPy library~\citep{RFletcher_1987,PVirtanen_etal_2020_Scipy}, again with default parameters.

\tcb{When solving problem~\eqref{prob:point_based_formulation_exact}, the parameter~$\varepsilon$ in the~MCF~\eqref{eq:MCF} is set adaptively at each candidate~$\lambda$ as~$\varepsilon(\lambda) = \eta g_{\min}(\lambda)$, where~$\eta \in (0,1]$ is a fixed constant (we use~$\eta = 0.5$ for all problems in Table~\ref{tab:test_prob}, except~$\eta = 0.4$ for ZLT1) and~$g_{\min}(\lambda)$ denotes the norm of the minimum-norm gradient across all objectives.
Since both~$\MCF_\varepsilon(\lambda)$ and~$\varepsilon(\lambda)$ depend on the gradients~$\nabla_\lambda(f_k(x(\lambda)))$, $k \in \{1,\ldots,q\}$, we discard any candidate~$\lambda$ for which~$g_{\min}(\lambda) < \tau_{\mathrm{norm}} = 10^{-6}$ or~$g_{\max}(\lambda)/g_{\min}(\lambda) > \tau_{\mathrm{cond}} = 10^{6}$ (where~$g_{\max}$ denotes the norm of the maximum-norm gradient across all objectives) to avoid numerical instability, assigning it a large fixed penalty of~$10^{8}$. This steers the outer~DFO search away from such~$\lambda$.}

\begin{figure}[t]
    \centering
    \includegraphics[width=0.95\textwidth]{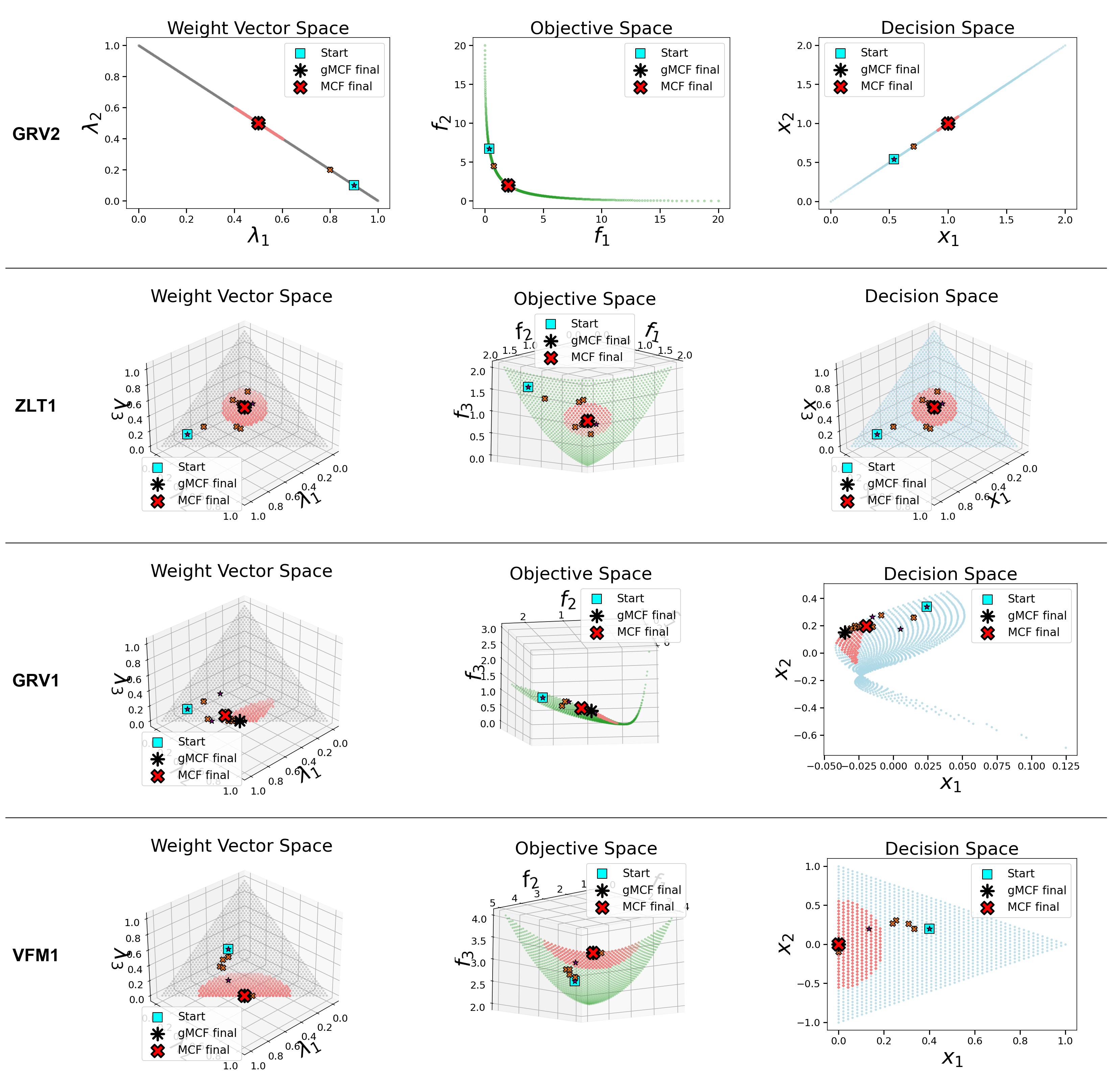}
    \caption{\tcb{Weight vector, objective, and decision spaces obtained when applying the MADS algorithm in the main step of Algorithm~\ref{alg:snee}, for problems GRV2, ZLT1, GRV1, and VFM1 (one row per problem). The cyan square marks the initial weight vector~$\lambda_0$; stars mark the MADS iterates. The black asterisk marks the knee obtained by minimizing~$\MCFg$~\eqref{eq:maximal_change_function} (``gMCF final''); the red cross marks the knee obtained by minimizing~$\MCF$~\eqref{eq:MCF} (``MCF final''). The red-shaded region shows the ellipsoidal neighborhood centered at the gMCF-based knee weight vector, with the corresponding Pareto sub-front and neighborhood of Pareto optimal solutions in the objective and decision spaces, respectively.}}
    \label{fig:WeightObjDecision-MADS}
\end{figure}

\begin{figure}[t]
    \centering
    \includegraphics[width=0.95\textwidth]{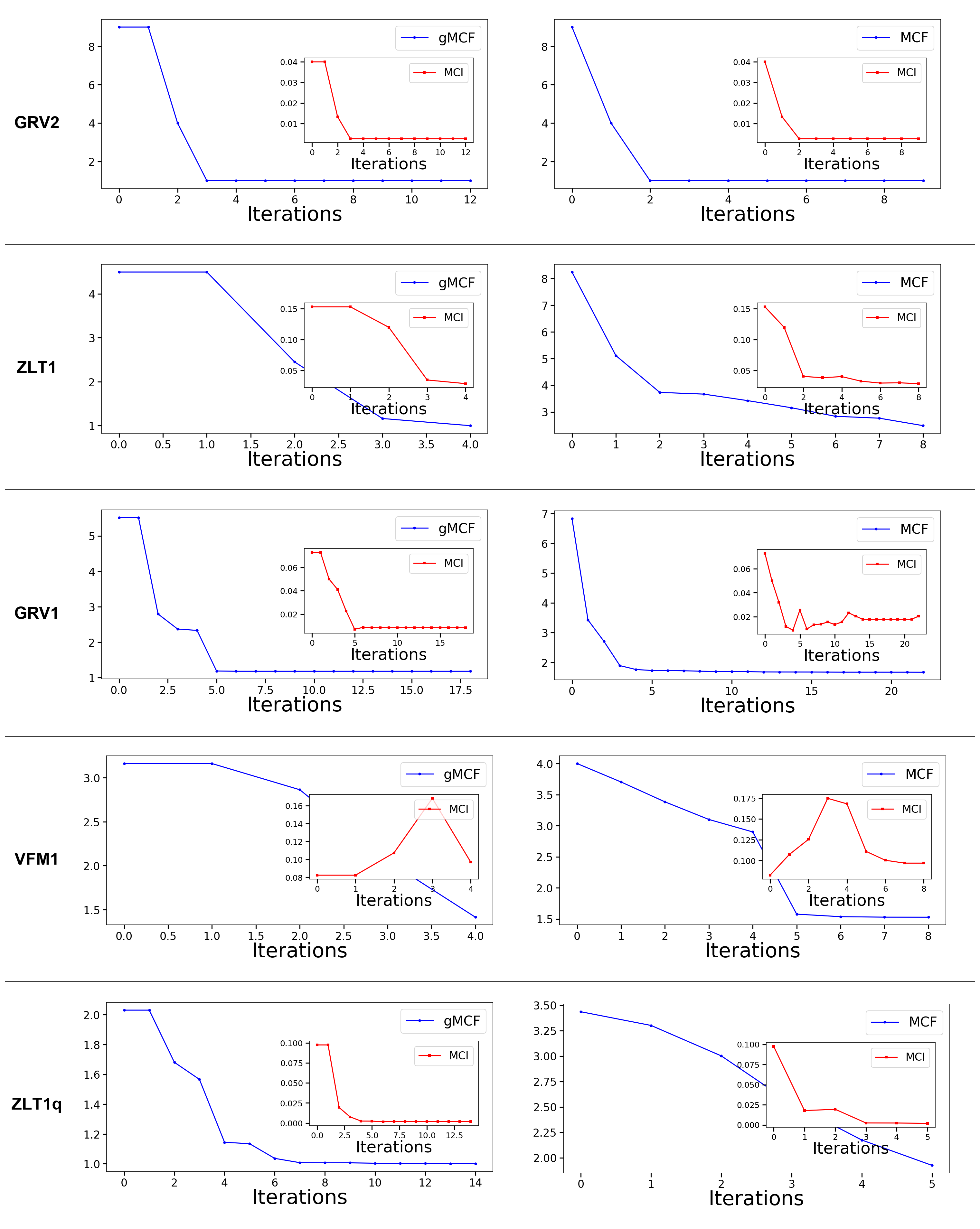}
    \caption{\tcb{Values of~$\MCFg$ (left column) and~$\MCF$ (right column) over the MADS iterations in the main step of Algorithm~\ref{alg:snee}, for problems GRV2, ZLT1, GRV1, VFM1, and ZLT1$q$ (one row per problem). Each panel's inset shows the corresponding~$\MCI$~\eqref{eq:metric} over the same iterations.}}
    \label{fig:Convergence-Unconstrained-MADS}
\end{figure}

\tcb{Figure~\ref{fig:WeightObjDecision-MADS} presents the weight vector, objective, and decision spaces for problems~GRV2, ZLT1, GRV1, and~VFM1 (one row per problem) when applying the MADS algorithm in the main step of Algorithm~\ref{alg:snee}. The} initial weight vector~$\lambda_0$ is set arbitrarily as follows: for problem~GRV2, $\lambda_0 = (0.9,0.1)$; for problems~ZLT1 and~GRV1, $\lambda_0 = (0.8,0.1,0.1)$; and for problems~VFM1 and~ZLT1$q$, $\lambda_0 = (0.4,0.2,0.4)$ and $\lambda_0 = (0.6,0.1,0.1,0.1,0.1)$, respectively. 
\tcb{In the weight vector space, the cyan square marks~$\lambda_0$, the stars mark the iterates, the black asterisk marks the knee weight vector obtained by minimizing~$\MCFg$ (``gMCF final''), and the red cross marks the knee weight vector obtained by minimizing~$\MCF$ (``MCF final''); the same color coding is used in the objective and decision spaces for the corresponding nondominated points and Pareto optimal solutions.}

Note that the knee nondominated points found by \tcb{the~MADS algorithm} lie in the intermediate regions of the Pareto fronts \tcb{for problems~GRV2 and~ZLT1, and on the upper edge curve of the Pareto front for problems~GRV1 and~VFM1}.
This outcome aligns with the concept of edge-knee nondominated points proposed in~\cite{KDeb_SGupta_2010} for problems with~2 objectives, but differs from the knee nondominated points found by the normal boundary intersection method~\citep{IDas_JEDennis_1998,IDas_1999}, which tends to identify points in the intermediate region of a Pareto front (when it is convex), depending on where the distance between a point on the front and the convex hull of the \tcb{single-objective (optimal) nondominated points} is maximized. 
If the front lacks a bulge, the point found by our \textit{snee} approach will not be in the intermediate region, as there is no clear knee nondominated point according to the verbal definition \tcb{in terms of the trade-off between improvement in one objective and deterioration in another}.

\tcb{Figure~\ref{fig:Convergence-Unconstrained-MADS} shows, for problems~GRV2, ZLT1, GRV1, VFM1, and~ZLT1$q$, the values of~$\MCFg$~\eqref{eq:maximal_change_function} (left column) and~$\MCF$~\eqref{eq:MCF} (right column) over the MADS iterations, with the corresponding~$\MCI$ shown in the inset of each panel.
One can observe from each figure that the~\tcb{MCI}~\eqref{eq:metric} is approximately correlated with both~MCF and~$\MCFg$, except for problem~VFM1. Note that the neighborhoods are required only for computing the~\tcb{MCI} and not for the~$\MCF$ or~$\MCFg$.}
For the purpose of evaluating the~\tcb{MCI}, we compute the local most-changing Pareto sub-fronts in the objective space. We visualize these sub-fronts along with the corresponding neighborhoods of weight vectors and Pareto optimal solutions in the weight vector and decision spaces, respectively. The neighborhoods of weight vectors shown in the figures are the ellipsoidal ones, centered at the \tcb{final} iterate returned by the~\tcb{MADS algorithm when minimizing the~gMCF}. For such neighborhoods, the size~$\alpha$ is set to~0.1 for all problems, except for~GRV2, where using a fixed neighborhood size results in ellipsoidal neighborhoods that either have no elements or contain the entire simplex set at certain iterations. Therefore, for~GRV2, we determine the neighborhood size using the adaptive rule described in Section~\ref{sec:neighborhoods} with~$\gamma=0.4$, which \tcb{empirically results in approximately~20\%} of the weight vectors in the fine-scale discretization of~$\Lambda$ \tcb{being} included in the neighborhood at each iteration.
The same adaptive rule will be used for the problems in Subsection~\ref{subsec:results_constr}. Remark~\ref{rem:MCFvsMCM} below explains that minimizing the~$\MCF$ \tcb{or~$\MCFg$} provides a more reliable approach than minimizing the~$\MCI$.

\begin{remark} \label{rem:MCFvsMCM}
Given the approximate correlation \tcb{of both~$\MCF$ and~$\MCFg$ with~$\MCI$ in most problems (except~VFM1)}, one might think that knee solutions could be obtained by minimizing the~$\MCI$ over~$\lambda_c$ instead of the~$\MCF$ over~$\lambda$. However, such an approach  presents challenges because~$\MCI=\MCI(\mathcal{E}_{\alpha}(\lambda_c))$ depends on the neighborhood size~$\alpha$. With fixed right-hand sides, the~$\MCI$ may achieve a small value for a neighborhood~$\mathcal{E}(\lambda_c)$ with a small size (particularly in regions of the weight vector space where the norm of~$\nabla_{\lambda}(F(x(\lambda_c)))$ is low), even if its center~$\lambda_c$ is not a knee solution. Therefore, minimizing the~$\MCF$ \tcb{or~$\MCFg$} provides a more reliable approach.
\hfill$\square$
\end{remark}

\subsection{Results for most-changing Pareto sub-fronts using Pareto sensitivity}\label{subsec:results_unconstr}

We now perform numerical experiments to assess which neighborhood is more effective at identifying the greatest variation across the objectives on the Pareto front through the corresponding Pareto sub-front. We consider the neighborhood achieving the highest value of the~\tcb{MCI}~\eqref{eq:metric} to be the most effective.
We again consider the problems from Table~\ref{tab:test_prob} in Subsection~\ref{subsec:moo_test_problems}.  In this section, we do not include results for problem~GRV2, where~$q = 2$, because both the ellipsoidal neighborhood~$\mathcal{E}_{\alpha}(\lambda_c)$ in~\eqref{eq:Ninverse} and neighborhood~$\mathcal{E}_{\beta}(\lambda_c)$ in~\eqref{eq:N} yield similar Pareto sub-fronts in the bi-objective case, providing little insight into our understanding of which one is more effective.
We use the weight vectors~$\lambda_0$ from Subsection~\ref{subsec:results_unc} as the neighborhood centers~$\lambda_c$ for problems~ZLT1, GRV1, VFM1, and~ZLT1$q$. 
We compare both neighborhoods~\eqref{eq:Ninverse} and~\eqref{eq:N} against a ball of radius~$r$, defined as follows
\begin{equation}\label{eq:ball} 
\mathcal{B}_r (\lambda_c) \; = \; \{ \lambda \in \mathbb{R}^q \; | \;
\| \lambda - \lambda_c \| \leq r \}.
\end{equation}

Note from~\eqref{eq:Ninverse} and~\eqref{eq:ball} that~$\mathcal{B}_r (\lambda_c)$ is equal to~$\mathcal{E}_{\alpha}(\lambda_c)$ when~$r = \alpha$ and the Jacobian of~$F(x(\cdot))$ is equal to the identity matrix. When using the ellipsoidal neighborhood~$\mathcal{E}_{\alpha}(\lambda_c)$, we set~$\alpha$ to~0.10, which results in neighborhoods that include roughly~10\%-25\% of the weight vectors in the fine-scale discretization of~$\Lambda$, \tcb{denoted by}~$\{\lambda_1, \ldots, \lambda_m\} \subset \Lambda$. In practice, the value of~$\alpha$ should be chosen by a decision maker. When using~$\mathcal{B}_r (\lambda_c)$ and~$\mathcal{E}_{\beta}(\lambda_c)$, we set~$r$ and~$\beta$ to values that produce neighborhoods containing approximately the same percentage of weight vectors as~$\mathcal{E}_{\alpha}(\lambda_c)$.

\begin{figure}[t]
    \centering
    \includegraphics[width=0.95\textwidth]{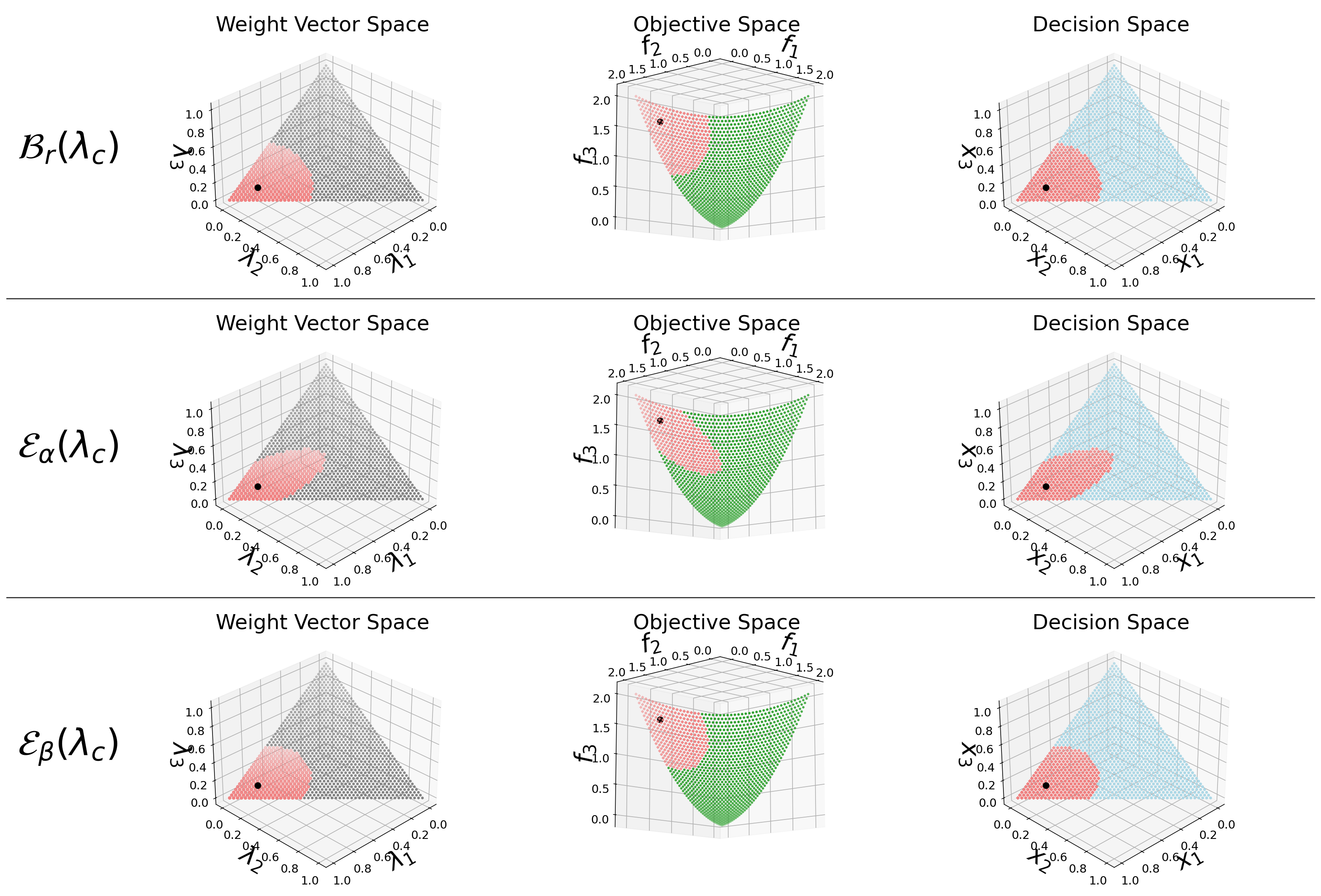}
    \caption{\tcb{Results for problem~ZLT1. Rows correspond to neighborhoods~$\mathcal{B}_r(\lambda_c)$, $\mathcal{E}_\alpha(\lambda_c)$, and $\mathcal{E}_\beta(\lambda_c)$; columns show the weight vector, objective, and decision spaces. The ellipsoidal neighborhood~$\mathcal{E}_\alpha(\lambda_c)$ is the most effective at capturing the greatest variation across the objectives.}}
    \label{fig:ZLT1}
\end{figure}

\begin{figure}[t]
    \centering
    \includegraphics[width=0.95\textwidth]{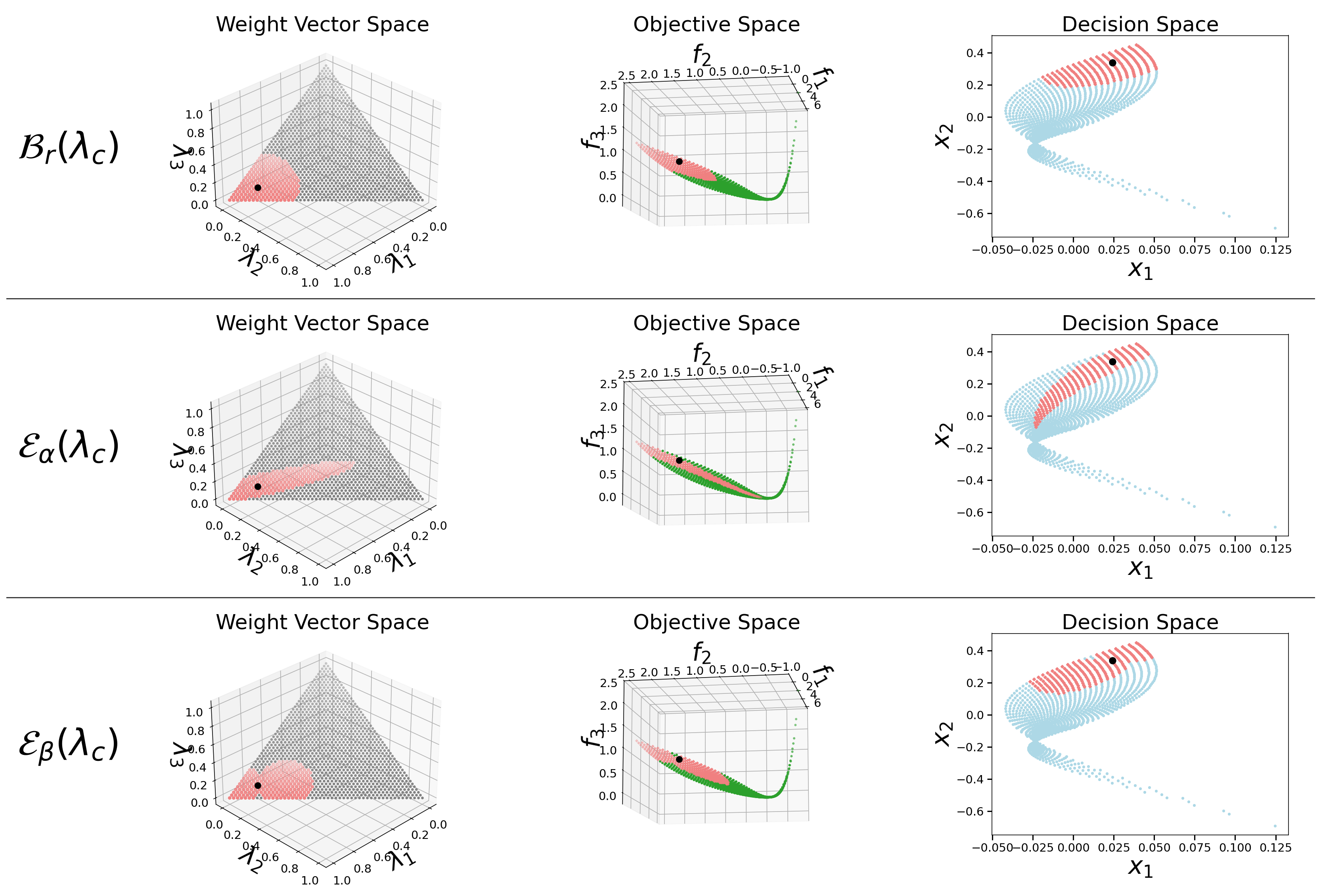}
    \caption{\tcb{Results for problem~GRV1 (for~$\bar{n}$ in Table~\ref{tab:test_prob} equal to~2). See the caption of Figure~\ref{fig:ZLT1} for details.}}
    \label{fig:GRV1}
\end{figure}

\begin{figure}[t]
    \centering
    \includegraphics[width=0.95\textwidth]{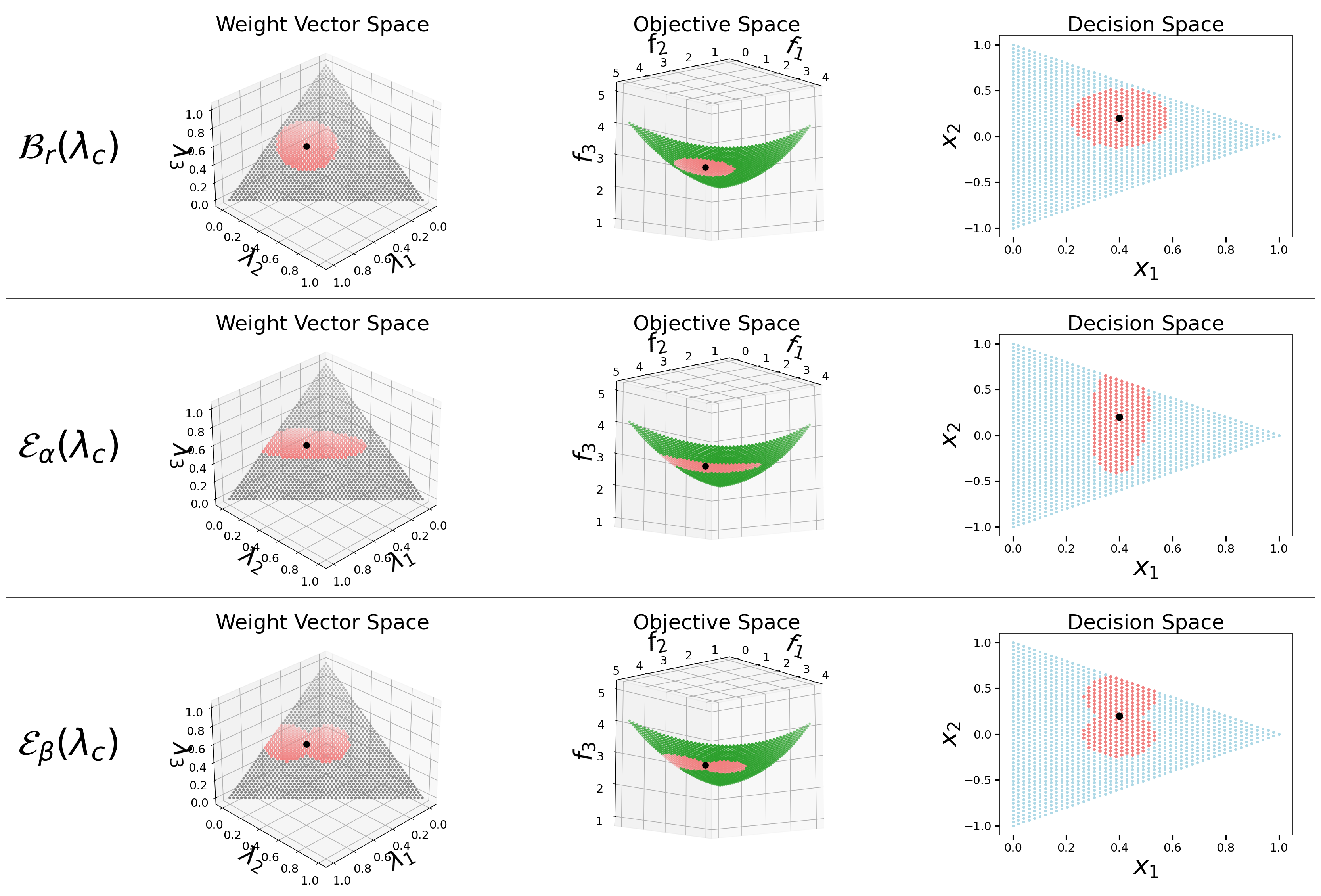}
    \caption{\tcb{Results for problem~VFM1. See the caption of Figure~\ref{fig:ZLT1} for details.}}
    \label{fig:VFM1}
\end{figure}

Figures~\ref{fig:ZLT1}--\ref{fig:VFM1} illustrate the results for problems~ZLT1, GRV1, and~VFM1. In each plot, the black dot represents the center of the corresponding neighborhood or sub-front.
Such figures show that~$\mathcal{E}_{\alpha}(\lambda_c)$ and~$\mathcal{E}_{\beta}(\lambda_c)$ lead to sub-fronts where nondominated points are distributed along directions of maximum change of the Pareto front. 
The values of the~\tcb{MCI}~\eqref{eq:metric} for problems~ZLT1, GRV1, VFM1, and~ZLT1$q$ are included in Table~\ref{tab:comparison_neighborhoods} below for each type of neighborhood, along with the values of~$n$ and~$q$ used in the experiments and the fraction of the weight vectors in~$\Lambda_m$ that are contained within each neighborhood. 
Similar to Subsection~\ref{subsec:results_unc}, to obtain the Pareto optimal solutions associated with each weight vector in a neighborhood, we minimized the corresponding weighted-sum function using the~BFGS algorithm implementation available in the Python SciPy library~\citep{RFletcher_1987,PVirtanen_etal_2020_Scipy}, with default parameters.
When computing the~\tcb{MCI} for the ball~\eqref{eq:ball}, we replace~$\mathcal{E}(\lambda_c)$ in~\eqref{eq:metric} with~$\mathcal{B}_r (\lambda_c)$.

    \begin{table}
    \centering
    \begin{tabular}{c|c|c|c|c|c}
         Problem & $n$ & $q$ & Neighborhood & $\MCI$ & (\# Neigh.~Points)$/|\Lambda_m|$\\ \hline
         \multirow{3}{*}{ZLT1} & \multirow{3}{*}{3} & \multirow{3}{*}{3} & $\mathcal{B}_r (\lambda_c)$, $r=0.40$ & $0.0895$ & 0.2392\\[0.5ex]
         &  &  & $\mathcal{E}_\alpha (\lambda_c)$, $\alpha=0.10$ & $0.1529$ & 0.2329\\[0.5ex]
         &  &  & $\mathcal{E}_\beta (\lambda_c)$, $\beta=7$\phantom{.11} & $0.0837$ & 0.2251\\ \hline 
         \multirow{3}{*}{GRV1} & \multirow{3}{*}{2} & \multirow{3}{*}{3} & $\mathcal{B}_r (\lambda_c)$, $r=0.30$ & 0.0136 & 0.1702\\[0.5ex]
          &  &  & $\mathcal{E}_\alpha (\lambda_c)$, $\alpha=0.10$ & 0.1078 & 0.1639\\[0.5ex]
          &  &  & $\mathcal{E}_\beta (\lambda_c)$, $\beta=10$\phantom{.1} & 0.0252 & 0.1749\\ \hline 
         \multirow{3}{*}{VFM1} & \multirow{3}{*}{2} & \multirow{3}{*}{3} & $\mathcal{B}_r (\lambda_c)$, $r=0.23$ & 0.0260 & 0.1788\\[0.5ex]
          &  &  & $\mathcal{E}_\alpha (\lambda_c)$, $\alpha=0.10$ & 0.0824 & 0.1804\\[0.5ex]
          &  &  & $\mathcal{E}_\beta (\lambda_c)$, $\beta=13$\phantom{.1} & 0.0546 & 0.1859 \\ \hline 
         \multirow{3}{*}{ZLT1$q$} & \multirow{3}{*}{5} & \multirow{3}{*}{5} & $\mathcal{B}_r (\lambda_c)$, $r=0.28$ & 0.0111 & 0.1026\\[0.5ex]
          &  &  & $\mathcal{E}_\alpha (\lambda_c)$, $\alpha=0.10$ & 0.0974 & 0.0948\\[0.5ex]
          &  &  & $\mathcal{E}_\beta (\lambda_c)$, $\beta=8.5$\phantom{1} & 0.0176 & 0.1007\\  
    \end{tabular}
    \caption{Comparison of~$\mathcal{B}_{r}(\lambda_c)$, $\mathcal{E}_{\alpha}(\lambda_c)$, and~$\mathcal{E}_{\beta}(\lambda_c)$ from the numerical experiments in Subsection~\ref{subsec:results_unconstr}}.\label{tab:comparison_neighborhoods}
    \end{table}

 From the~\tcb{MCI} column of Table~\ref{tab:comparison_neighborhoods}, one can observe that~$\mathcal{E}_{\alpha}(\lambda_c)$ outperforms~$\mathcal{E}_{\beta}(\lambda_c)$ and~$\mathcal{B}_{r}(\lambda_c)$ on all of the problems, even though all neighborhoods are comparable in size (see the last column of Table~\ref{tab:comparison_neighborhoods}). Additionally, $\mathcal{E}_{\beta}(\lambda_c)$ yields (nearly all the times) better results than~$\mathcal{B}_{r}(\lambda_c)$.

\section{The constrained case}\label{sec:constr}

Throughout this section, we focus on the case where~$X$ is not $\mathbb{R}^n$.

\subsection{Pareto sensitivity calculation in the constrained case}\label{subsec:pareto_sensitivity_calculation_constr} 

In this subsection, we extend the~\textit{snee} approach in Subsection~\ref{subsec:pareto_sensitivity_calculation} to handle constrained~MOO problems. 
Such an extension only requires a recalculation of the sensitivity~$\nabla x(\lambda)$ under the presence of constraints.
Specifically, instead of deriving~$\nabla x(\lambda)$ from the first-order necessary conditions for unconstrained optimality~\eqref{eq:1storderKKT} in Subsection~\ref{subsec:pareto_sensitivity_calculation}, we derive it from the corresponding conditions under the presence of constraints, typically referred to as the first-order~KKT conditions, in particular from the equality part of such conditions (the so-called~KKT system). 

Let us consider the following constrained~MOO problem 
\begin{equation}\label{eq:prob_constr}
\min_{x \in X} F(x), \ \text{ with } \ X = \{ x \in \mathbb{R}^n ~|~ c_j(x) \leq 0, \; j \in I, \text{ and } c_j(x) = 0, \; j \in E\},
\end{equation}
where~$I$ and~$E$ are finite index sets used to denote inequality and equality constraints, respectively.
We start with the derivation of~$\nabla x(\lambda)$ using the first-order~KKT system associated with problem~\eqref{eq:prob_constr}.
Denoting~$c_I(x) = (c_j(x), \, j \in I)$ and~$c_E(x) = (c_j(x), \, j \in E)$, and given~$\lambda \in \Lambda$, the Lagrangian function of problem~\eqref{eq:prob_constr} is defined as~$\mathcal{L} (x,z_I,z_E) = \sum_{i=1}^{q} \lambda_i f_i(x) + c_I(x)^\top z_I + c_E(x)^\top z_E$, where~$z_I$ and $z_E$ are vectors of Lagrange multipliers.
In Assumption~\ref{ass:LL_assumptions_constr} below, we assume that the constraint functions are twice continuously differentiable and there exists a solution~$x(\lambda)$ satisfying the~KKT conditions for problem~\eqref{eq:prob_constr} with associated vectors of multipliers~$(z_I(\lambda), z_E(\lambda))$. Specifically, the first-order~KKT system for problem~\eqref{eq:prob_constr} at~$x(\lambda)$ is given by~\citep[Theorem 3.1.5]{KMiettinen_1999}
\begin{equation}\label{eq:KKT_syst}
    \begin{cases}
        \sum_{i=1}^q \lambda_i \nabla_x f_i (x(\lambda)) + \nabla_x c_I(x(\lambda)) \, z_I(\lambda) + \nabla_x c_E(x(\lambda)) \, z_E(\lambda) = 0,\\
        z_{I}(\lambda) \circ c_{I}(x(\lambda)) = 0,\\
        c_E(x(\lambda)) = 0,
    \end{cases}
\end{equation}
where $\circ$ denotes the element-wise multiplication of two vectors. We do not include the non-negativity of the multipliers in~$z_I(\lambda)$ and the satisfaction of the inequality constraints in~\eqref{eq:KKT_syst} because they are not necessary for the derivation below. 

We also assume a well-known set of conditions in constrained optimization theory, 
as detailed in Assumption~\ref{ass:LL_assumptions_constr} below, in order to secure the required matrix invertibility in the derivative calculation. Specifically, we require the linear independence constraint qualification~(LICQ), strict complementarity slackness~(SCS), and the second-order sufficient condition~(SOSC), which are formally stated in~\citet[Chapter~12, Definitions~12.4-12.5 and Theorem~12.6, respectively]{JNocedal_SJWright_2006} and summarized next.
Given Lagrange multipliers~$(z_I(\lambda), z_E(\lambda))$ corresponding to a solution~$x(\lambda)$ of the~KKT conditions, we define the cone of critical directions as:
\begin{equation}\label{eq:critical_cone}
Z(\lambda) \; = \; \left\{\begin{matrix}
    \phantom{d^x \ne 0 ~:~  } \hspace{0.0cm}\nabla_x c_i(x(\lambda))^\top d^x \le 0, \ \forall i \in I(\lambda) \phantom{\text{ with } (z_I(\lambda))_i > 0}\\
    d^x \ne 0 ~:~ \nabla_x c_i(x(\lambda))^\top d^x = 0, \ \forall i \in I(\lambda) \text{ with } (z_I(\lambda))_i > 0\\ 
    \phantom{d^x \ne 0 : } \nabla_x c_i(x(\lambda))^\top d^x = 0, \ \forall i \in E \phantom{\text{ with } (z_I(\lambda))_i > 0}
    \end{matrix} \right\},
\end{equation}
where~$I(\lambda)$ denotes the index set of inequality constraints that are active at~$x(\lambda)$. LICQ requires that the gradients~$\nabla_x c_i(x(\lambda))$, for all~$i \in I(\lambda)\cup E$,
are linearly independent. SCS requires that for any multipliers~$(z_I(\lambda),z_E(\lambda))$ satisfying the~KKT conditions at~$x(\lambda)$, the strict positivity~$z_I(\lambda)_i > 0$ holds for every~$i \in I(\lambda)$. Finally, SOSC requires that for any multipliers~$(z_I(\lambda), z_E(\lambda))$ satisfying the~KKT conditions at~$x(\lambda)$ and for any~$d^x \in Z(\lambda)$, with~$Z(\lambda)$ as defined in~\eqref{eq:critical_cone}, the curvature condition~$(d^x)^\top \nabla^2_{xx} \mathcal{L}(x(\lambda),z_I(\lambda),z_E(\lambda)) d^x > 0$ is satisfied.

\begin{assumption}[Existence of Pareto optimal solutions (LL constrained case)]\label{ass:LL_assumptions_constr}
\phantom{.}\linebreak \tcb{Let~$X$, as defined in problem~\eqref{eq:prob_constr}, be nonempty.} The constraint functions~$c_j$, with~$j \in I \cup E$, are twice continuously differentiable. There exists a solution~$x(\lambda)$ satisfying the~KKT conditions with associated multipliers~$(z_I(\lambda), z_E(\lambda))$ such that the~LICQ, SCS, and~SOSC are satisfied. 
\end{assumption}

Under Assumption~\ref{ass:LL_assumptions_constr}, based on~\citet[p.~291, Theorem~2.1]{AFiacco_1976}, the vectors of multipliers~$z_I(\lambda)$ and~$z_E(\lambda)$ associated with~$x(\lambda)$ are unique, and the vector \linebreak function~$w(\lambda) = (x(\lambda), z_I(\lambda), z_E(\lambda))^{\top}$ is once continuously differentiable for any given~$\lambda$.
Let us now introduce a vector function~$K$ such that the~KKT system~\eqref{eq:KKT_syst} can be written as~$K(w(\lambda)) = 0$. 
Applying the chain rule to such an equation, we obtain~$\nabla_w K^{\top} \nabla w^{\top} = - \nabla_{\lambda} K^{\top}$, with 
\begin{equation}\label{eq:jacobians}
\nabla_{\lambda} K^{\top} = \begin{pmatrix} \nabla_{x\lambda}^2 \mathcal{L} \\  0 \\ 0  \end{pmatrix}
\; \text{ and } \;
\nabla_w K^{\top} = \begin{pmatrix} \nabla_{xx}^2 \mathcal{L} && \nabla_x c_I && \nabla_x c_E \\
z_I \circ \nabla_x c_I^{\top} && C_{I} && 0 \\
\nabla_x c_E^{\top} && 0 && 0 \end{pmatrix},
\end{equation}
where~$\nabla_{x\lambda} \mathcal{L}$ and~$\nabla_{xx} \mathcal{L}$ are evaluated at~$w(\lambda)$, the Jacobian matrices~$\nabla_x c_I^{\top}$ and~$\nabla_x c_E^{\top}$ are evaluated at~$x(\lambda)$, $C_I$ is a diagonal matrix with elements defined by~$c_I(x(\lambda))$, and~$z_I \circ \nabla_\lambda c_I^{\top}$ represents a matrix formed by element-wise multiplication of the entries of~$z_I$ with the corresponding rows of~$\nabla_x c_I^{\top}$. 

Under Assumption~\ref{ass:LL_assumptions_constr}, it is well known that the Jacobian~$\nabla_w K^{\top}$ is non-singular at~$w(\lambda)$, see~\citet[p.~294, Eq.~(2.4)]{AFiacco_1976} or~\citet[p.~225, Eq.~(5.2)]{AVFiacco_YIshizuka_1990}. Therefore, we have
\begin{equation}\label{jacobian_pi_11-19-2024}
    \nabla w \; = \; \begin{pmatrix} \nabla x, \nabla z_I, \nabla z_E \end{pmatrix} \; = \; - \nabla_{\lambda} K \, \nabla_w K^{-1}.
\end{equation}
Following~\citet[Subsection 2.2]{TGiovannelli_GKent_LNVicente_2022}, we can now introduce a matrix~$L = \begin{pmatrix} \mathbf{I}_n & \mathbf{0}\end{pmatrix}^{\top}$ to extract the columns of~\eqref{jacobian_pi_11-19-2024} that correspond to the~$\nabla x(\lambda)$ term, where~$\mathbf{I}_n$ denotes an identity matrix of size~$n$ and~$\mathbf{0}$ represents a null matrix of dimensions~$n \times \big( \vert I \vert + \vert E \vert \big)$, resulting in
\begin{equation}\label{nabla_y}
    \nabla x(\lambda) \; = \; - \nabla_\lambda K \nabla_w K^{-1} L.
\end{equation}

By applying the chain rule as we did to derive~\eqref{eq:adjoint}, but now using the matrix~$\nabla x(\lambda)$ from~\eqref{nabla_y} instead of~\eqref{eq:04}, we obtain that the transpose of the Jacobian matrix of~$F(x(\cdot))$ at~$\lambda$ is given by
\begin{equation}\label{eq:adjoint_2}
\begin{alignedat}{2}
    \nabla_{\lambda}(F(x(\lambda))) \; &= \; \nabla x(\lambda) \nabla_x F(x(\lambda)) \\ \; &= \;  - \nabla_\lambda K \nabla_w K^{-1} L (\nabla_x f_1, \ldots, \nabla_x f_q).
\end{alignedat}
\end{equation}
In~\eqref{eq:adjoint_2}, all gradients~$\nabla_x f_i$ are evaluated at $x(\lambda)$, with~$i \in \{1, \ldots, q\}$, and~$\nabla_\lambda K$ and~$\nabla_w K$ are evaluated at~$w(\lambda)$. 

We note that classical multi-objective optimization textbooks often invoke the Kuhn--Tucker constraint qualification~(KTCQ), see~\citet[p.~57, Definition~2.50]{MEhrgott_2005} and~\citet[p.~38, Definition~3.1.3]{KMiettinen_1999}, originally introduced in~\citet{HWKuhn_AWTucker_1951}. The KTCQ requires that for every direction in the linearized feasible cone, there exists a smooth feasible arc tangent to that direction, ensuring that the feasible set admits a local differentiable description at the solution. In those references, the~KTCQ is used either to establish necessary and sufficient conditions for properly efficient solutions or to guarantee that the KKT system~\eqref{eq:KKT_syst} holds as a necessary optimality condition when applying the weighted-sum method. We remark that our Assumption~\ref{ass:LL_assumptions_constr} already subsumes the KTCQ. As shown in~\citet[p.~645, Fig.~4]{DWPeterson_1973}, LICQ implies KTCQ, and therefore Assumption~\ref{ass:LL_assumptions_constr} ensures the validity of the~KTCQ invoked in those references.

\tcb{The example below illustrates two features of Pareto sensitivity in the presence of constraints. First, the variation of a Pareto optimal solution with respect to the weight vector is restricted to directions tangent to the feasible set~$X$. Second, the resulting objective sensitivities depend on the weight vector only through such tangent directions, via the constrained sensitivity~$\nabla x(\lambda)$. In this example, the~$\MCF$ and~$\MCFg$ criteria happen to attain the same value at every~$\lambda$ and are therefore minimized by the same~$\lambda_*$.}

\begin{example}\label{ex:constrained}
Consider the constrained problem with $n=2$ and $q=2$ given by $f_1(x) = (x_1^2+x_2^2)/2$ and $f_2(x) = ((x_1-2)^2+x_2^2)/2$, subject to the single equality constraint $c(x) = x_1+x_2-1 = 0$. Both objectives are strongly convex and the constraint is linear, so LICQ holds at every feasible point and Assumption~\ref{ass:LL_assumptions_constr} is satisfied.
To compute $\nabla_\lambda(f_i(x(\lambda)))$, $i \in \{1,2\}$, as required by $\eqref{eq:MCF}$ and $\eqref{eq:maximal_change_function}$, we must differentiate $x(\lambda)$ with respect to $\lambda_1$ and $\lambda_2$. Writing the Lagrangian as $\mathcal{L}(x,z) = \lambda_1 f_1(x) + \lambda_2 f_2(x) + z\,c(x)$, the KKT system~\eqref{eq:KKT_syst} reads
\[
\begin{cases}
\lambda_1 x_1 + \lambda_2(x_1-2) + z = 0,\\
\lambda_1 x_2 + \lambda_2 x_2 + z = 0,\\
x_1+x_2-1=0.
\end{cases}
\]
Solving for $x_1,x_2$ and writing $s = \lambda_1+\lambda_2$ gives~$x(\lambda) = \big( (\lambda_1+3\lambda_2)/(2s),\; (\lambda_1-\lambda_2)/(2s) \big)^\top$ and~$z(\lambda) = (\lambda_2-\lambda_1)/2$.
Let $t\in(0,1)$ parametrize~$\Lambda$ via $\lambda_1=t$ and $\lambda_2=1-t$ (so $s=1$). We have $x(\lambda) = ((3-2t)/2,\; (2t-1)/2)^\top$ and $z(\lambda) = (1-2t)/2$. The equality constraint~$c(x)=0$ requires~$x(\lambda)$ to remain on the line~$x_1+x_2=1$ for every~$\lambda$. Thus, as~$t$ varies,~$x(\lambda)$ can change only in a direction tangent to the constraint. Differentiating with respect to~$t$ gives~$
dx(\lambda)/dt = \left(-1,\;1\right)^\top$,
which is a tangent direction because~$\nabla_x c(x(\lambda))^\top (dx(\lambda)/dt) = (1,1)(-1,1)^{\top} =0$.
Differentiating $x(\lambda)$ with respect to $\lambda_1,\lambda_2$ and evaluating at $\lambda_1=t$, $\lambda_2=1-t$ gives
\[
\nabla x(\lambda) \; = \; \begin{pmatrix} t-1 & 1-t \\ t & -t \end{pmatrix}.
\]
Note that both rows of~$\nabla x(\lambda)$, $(t-1)(1,-1)^\top$ and $t(1,-1)^\top$, are each proportional to $(1,-1)$.
The objective gradients at $x(\lambda)$ are~$\nabla_x f_1(x(\lambda)) = \left((3-2t)/2,\; (2t-1)/2\right)^{\!\top}$, $\nabla_x f_2(x(\lambda)) = \left((-1-2t)/2,\; (2t-1)/2\right)^{\!\top}$, and therefore, by~\eqref{eq:adjoint_2}, $g_1~:=~\nabla_{\lambda}(f_1(x(\lambda))) = \nabla x(\lambda)\,\nabla_x f_1 = -2(1-t)(1-t,-t)^{\top}$ and~$g_2~:=~\nabla_{\lambda}(f_2(x(\lambda))) = \nabla x(\lambda)\,\nabla_x f_2 = 2t\,(1-t,-t)^{\top}$. Thus,~$\MCFg(\lambda) \;=\; \max\{(1-t)/t,\; t/(1-t)\}$,
which is minimized at $t_*=1/2$, giving $\lambda_*=(1/2,1/2)^{\top}$ and knee solution $x(\lambda_*)=(1,0)^\top$ with $\MCFg(\lambda_*)=1$.

For $t\in(0,1)$, the tangent cone is $T_\Lambda(\lambda)=\{d:d_1+d_2=0\}$, with unit directions $d^+=(1,-1)/\sqrt2=-d^-$. Since $g_1^\top d^+=-\sqrt2(1-t)$ and $g_2^\top d^+=\sqrt2\,t$, $d^+$ is the unique unit direction with $g_1^\top d\le-\varepsilon$ and $g_2^\top d\ge\varepsilon$ whenever $\varepsilon\le\sqrt2\min\{t,1-t\}$, so $D_\varepsilon^{12}(\lambda)=\{d^+\}$. By symmetry, $D_\varepsilon^{21}(\lambda)=\{d^-\}$. Substituting into~\eqref{eq:MCF} gives
\[
\MCF_\varepsilon(\lambda) = \max\Big\{\frac{-g_1^\top d^+}{g_2^\top d^+},\; \frac{-g_2^\top d^-}{g_1^\top d^-}\Big\} = \max\{(1-t)/t,\;t/(1-t)\} = \MCFg(\lambda)
\]
for every $\varepsilon\in(0,\sqrt2\min\{t,1-t\}]$. 
Over the values of~$t$ satisfying this condition, both criteria are thus minimized at $t_*=1/2$, giving the same knee solution $x(\lambda_*)=(1,0)^\top$ with $\MCF_\varepsilon(\lambda_*)=\MCFg(\lambda_*)=1$, provided that~$0<\varepsilon\le1/\sqrt2$.
\hfill$\square$
\end{example}

\subsection{MOO constrained test problems}

Table~\ref{tab:test_prob_constr} specifies the constrained test problems considered in the experiments for Subsection~\ref{subsec:results_constr}, along with the number of variables and objectives (i.e.,~$n$ and~$q$, respectively). Problems~DAS1 and~DO2DK are well-known in the~MOO literature. They have a convex Pareto front, despite some of the objective functions being non-convex.

\begin{table}
    \footnotesize
    \centering
    \begin{tabular}{ c|c|c|c|l } 
    Problem & $n$ & $q$ & Ref. & \multicolumn{1}{c}{Objective Functions and Constraints} \\[4pt]
    \hline
    \multirow{6}{*}{DAS1} & \multirow{6}{*}{$5$} & \multirow{6}{*}{$2$} & \multirow{6}{*}{\cite{IDas_JEDennis_1998}} & $f_1(x_1,\ldots,x_5) = x_1^2 + x_2^2 + x_3^2 + x_4^2 + x_5^2$ \\
 &&&& $f_2(x_1,\ldots,x_5) = 3x_1 + 2x_2 - x_3/3 + 0.01 (x_4 - x_5)^3$\\
    &&&& $c_1(x_1,\ldots,x_5) = x_1^2 + x_2^2 + x_3^2 + x_4^2 + x_5^2 - 10 \le 0$ \\
    &&&& $c_2(x_1,\ldots,x_5) = x_1 + 2 x_2 - x_3 - 0.5 x_4 + x_5 - 2 = 0$ \\
    &&&& $c_3(x_1,\ldots,x_5) = 4x_1 - 2 x_2 + 0.8 x_3 + 0.6 x_4 + 0.5 x_5^2 = 0$
    \\[4pt]
    \hline
    \multirow{6}{*}{DO2DK} & \multirow{6}{*}{$30$} & \multirow{6}{*}{$2$} & \multirow{6}{*}{\cite{JBranke_etal_2004,KDeb_LThiele_etal_2002,KDeb_1999_MOO}} & $f_1(x) =g_1(x) g_2(x_1)\left(\sin \left(\pi x_1 / 2+\pi\right)+1\right)$ \\
 &&&& $f_2(x) = g_1(x) g_2(x_1)\left(\cos \left(\pi x_1 / 2+\pi\right)+1\right)$\\
    & & & & $g_1(x) =1+\frac{9}{n-1} \sum_{i=2}^n x_i$\\
    &&&& $g_2(x_1)  = 5+10\left(x_1-0.5\right)^2+\cos \left(2 \pi x_1\right) \cdot 2^{\frac{1}{2}}$ \\
    &&&& $c_j(x) = -x_j \le 0, \quad j \in \{1,2, \ldots, n\}$ \\
    &&&& $c_{n+j}(x) = x_j - r \le 0, \quad j \in \{1,2, \ldots, n\}$
    \\[4pt]
    \hline
    \multirow{5}{*}{VFM1\textit{constr}} & \multirow{5}{*}{$2$} & \multirow{5}{*}{$3$} & & $f_1(x_1, x_2) = x_1^2 + (x_2 - 1)^2$ \\
 &&&& $f_2(x_1, x_2) = x_1^2 + (x_2 + 1)^2 + 1$\\
 &&&& $f_3(x_1, x_2) = (x_1 - 1)^2 + x_2^2 + 2$\\
    &&&& $c_1(x_1, x_2) = x_1^2 + x_2^2 - 0.8 \le 0$ \\
    &&&& $c_2(x_1, x_2) = (x_1 - 1)^2 + x_2^2 - 1\le 0$ \\[4pt]    
    \end{tabular}
    \caption{Constrained test problems ($r$ in~DO2DK is an arbitrary positive scalar used to tighten the feasible set).}\label{tab:test_prob_constr}
\end{table}

\subsection{Results for knee solutions using Pareto sensitivity in the constrained case}
\label{subsec:results_constr}

\tcb{To solve either problem~\eqref{prob:point_based_formulation_exact} or~\eqref{prob:point_based_formulation}} in the main step of Algorithm~\ref{alg:snee}, where~$\nabla_{\lambda} (f_i(x(\lambda)))$ and~$\nabla_{\lambda} (f_{j}(x(\lambda)))$ in the~$\MCF$ \tcb{and~$\MCFg$} now correspond to the columns of matrix~$\nabla_{\lambda}(F(x(\lambda)))$ in~\eqref{eq:adjoint_2} instead of~\eqref{eq:adjoint}, \tcb{we employ the MADS algorithm~\citep{CAudet_SLDigabel_VRochonMontplaisir_CTribes_2022}, as in Subsection~\ref{subsec:results_unc}, run with default parameters. When solving problem~\eqref{prob:point_based_formulation_exact}, the parameter~$\varepsilon$ in the~$\MCF$~\eqref{eq:MCF} is set adaptively as described in Subsection~\ref{subsec:results_unc}, with~$\eta = 0.5$ for all problems in Table~\ref{tab:test_prob_constr}, except~$\eta = 0.3$ for DO2DK~($r=0.5$). 
The same conditioning guard on~$g_{\min}(\lambda)$ and~$g_{\max}(\lambda)/g_{\min}(\lambda)$ described in Subsection~\ref{subsec:results_unc} is applied here as well, with the constrained-case gradients~$\nabla_{\lambda}(f_k(x(\lambda)))$ computed via~\eqref{eq:adjoint_2}.
We also solved problem~\eqref{prob:point_based_formulation_exact} and~\eqref{prob:point_based_formulation} with the NM and DIRECT algorithms introduced in Section~\ref{subsec:point_based}. Since the resulting knee solutions are nearly identical to those obtained with MADS for problems~DAS1 and~DO2DK (but differ for problem~VFM1\textit{constr}, as discussed at the end of this subsection), we report the full NM and DIRECT results in the Appendix rather than here.}


We considered two configurations for problem~DO2DK by selecting the constraint right-hand side~$r$ from~$\{0.5,1\}$, where~$1$ corresponds to a larger feasible set and~$0.5$ results in a tighter feasible set. 
\tcb{The} initial weight vector~$\lambda_0$ in Algorithm~\ref{alg:snee} is set arbitrarily as~$(0.4,0.6)$, $(0.2,0.8)$, and~$(0.4,0.2,0.4)$ for problems~DAS1, DO2DK, and~VFM1\textit{constr}, respectively. 
To perform~\textbf{Step~A} and obtain~$x(\lambda)$ for a given~$\lambda$, we solved the weighted-sum problem~\eqref{prob:multiobj_weightedsum} using the~SLSQP algorithm~\citep{DKraft_1988}, which is designed for solving constrained optimization problems through sequential quadratic programming. We ran the~SLSQP algorithm implementation available in the Python SciPy library~\citep{PVirtanen_etal_2020_Scipy}, with default parameters.

\begin{figure}[t]
    \centering
    \includegraphics[width=0.95\textwidth]{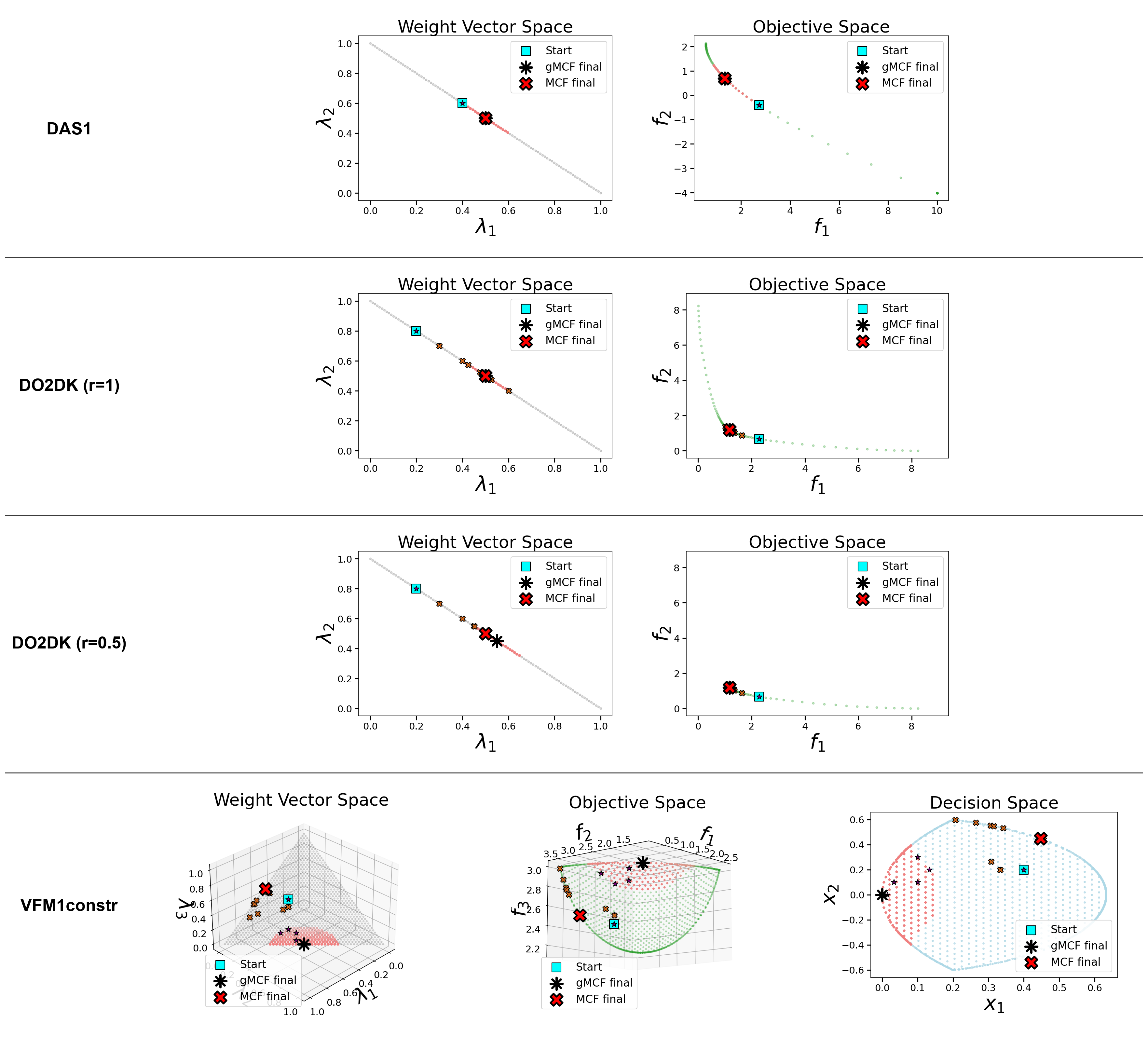}
    \caption{\tcb{Weight vector, objective, and decision spaces obtained when applying the MADS algorithm in the main step of Algorithm~\ref{alg:snee}, for problems DAS1, DO2DK ($r=1$), DO2DK ($r=0.5$), and VFM1\textit{constr} (one row per problem/configuration). See the caption of Figure~\ref{fig:WeightObjDecision-MADS} for the marker legend.}}
    \label{fig:WeightObjDecision-Constrained-MADS}
\end{figure}

\begin{figure}[t]
    \centering
    \includegraphics[width=0.95\textwidth]{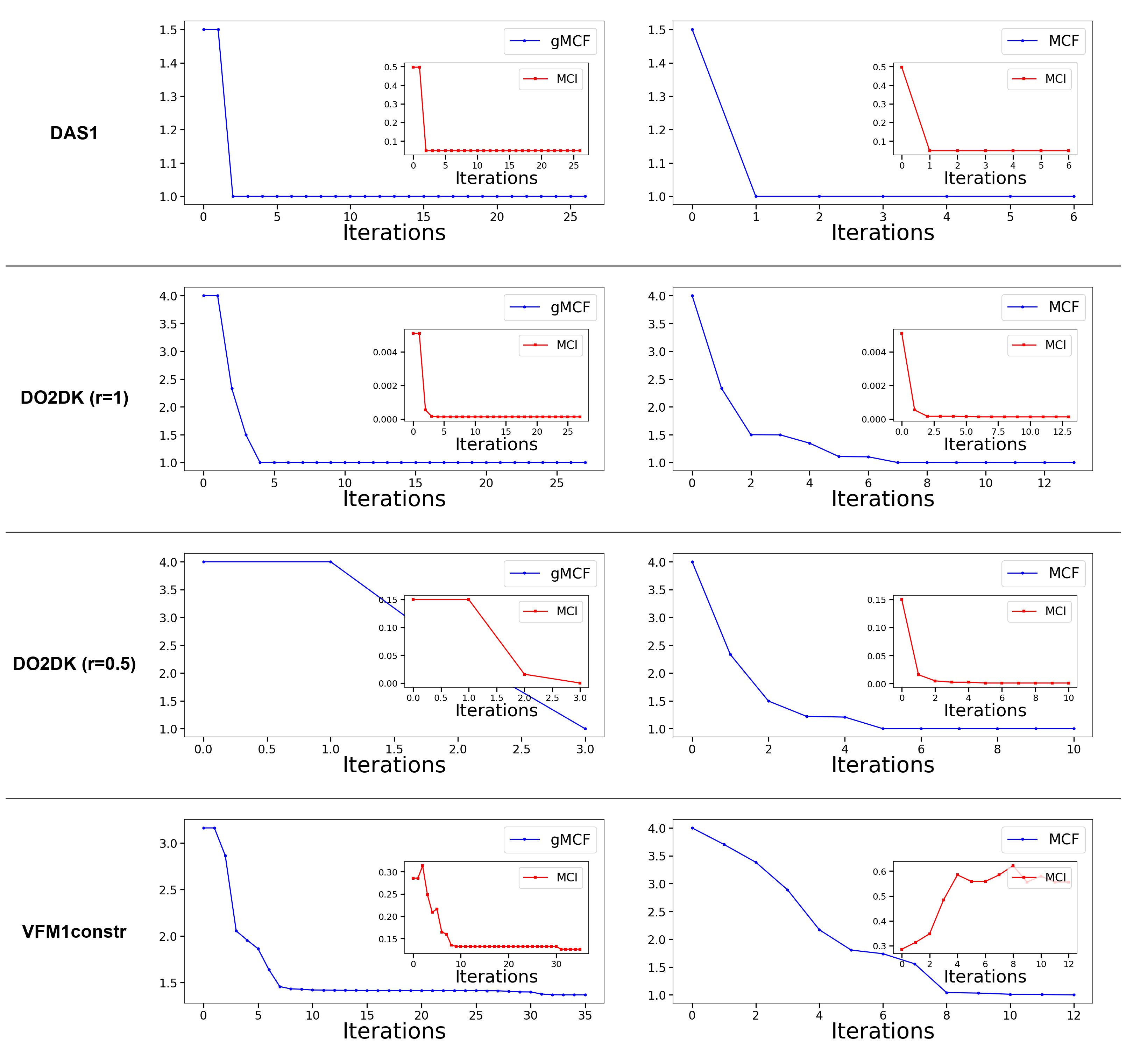}
    \caption{\tcb{Values of~$\MCFg$ (left column) and~$\MCF$ (right column) over the MADS iterations in the main step of Algorithm~\ref{alg:snee}, for problems DAS1, DO2DK ($r=1$), DO2DK ($r=0.5$), and VFM1\textit{constr}. Each panel's inset shows the corresponding~$\MCI$ over the same iterations.}}
    \label{fig:Convergence-Constrained-MADS}
\end{figure}

\tcb{Figure~\ref{fig:WeightObjDecision-Constrained-MADS} presents the weight vector, objective, and decision spaces for problems~DAS1, DO2DK ($r=1$), DO2DK ($r=0.5$), and~VFM1\textit{constr} when applying the MADS algorithm in the main step of Algorithm~\ref{alg:snee} (see the caption of Figure~\ref{fig:WeightObjDecision-MADS} for the marker legend).}
Figure~\ref{fig:WeightObjDecision-Constrained-MADS} \tcb{shows} that the knee nondominated points for problem~DAS1 and for problem~DO2DK when~$r = 1$ (i.e., the feasible set is relatively large) lie in the intermediate region of the Pareto front. In contrast, \tcb{for} DO2DK with~$r = 0.5$ and for~VFM1\textit{constr}, our \textit{snee} approach leads to knee nondominated points that correspond to points on the boundary of the Pareto front.
Similar to Subsection~\ref{subsec:results_unc}, the reason why we do not observe a knee nondominated point in the intermediate region for~DO2DK with~$r = 0.5$ and~VFM1\textit{constr} is due to the lack of a bulge in the Pareto front within the interior of the feasible set. Instead, the normal boundary intersection method would have returned a point in the intermediate region, regardless of the presence of a bulge.

\tcb{Figure~\ref{fig:Convergence-Constrained-MADS} shows, for problems~DAS1, DO2DK ($r=1$), DO2DK ($r=0.5$), and~VFM1\textit{constr}, the values of~$\MCFg$ and~$\MCF$ over the MADS iterations, with the corresponding~$\MCI$ shown in the inset of each panel. One can observe that the~MCI is approximately correlated with both~$\MCF$ and~$\MCFg$, except for problem~VFM1\textit{constr}.}

Similar to Subsection~\ref{subsec:results_unc}, the neighborhoods of weight vectors shown in \tcb{Figure~\ref{fig:WeightObjDecision-Constrained-MADS}} are the ellipsoidal ones, centered at the \tcb{final iterate returned by the MADS algorithm when minimizing the~gMCF}. For all the problems, using a fixed neighborhood size~$\alpha$ results in ellipsoidal neighborhoods that either have no elements or contain the entire simplex set at certain iterations. Therefore, we adopt the same adaptive rule used for~GRV2 in Subsection~\ref{subsec:results_unc}. Again, the neighborhoods are required only for computing the~\tcb{MCI} and not for the~$\MCF$ \tcb{or~$\MCFg$}.

\tcb{For problem~VFM1\textit{constr}, Figure~\ref{fig:WeightObjDecision-Constrained-MADS} shows that the knee nondominated point obtained by minimizing~$\MCFg$ (``gMCF final'') differs from the one obtained by minimizing~$\MCF$ (``MCF final''). Specifically, the gMCF-based knee lies near the left edge curve of the Pareto front, while the MCF-based knee lies near the upper edge curve. Analogous behavior is observed with the NM and DIRECT algorithms in the Appendix, where both knee nondominated points found by~NM lie on the left edge curve of the Pareto front (Figure~\ref{fig:WeightObjDecision-Constrained-NM}), whereas the ones found by~DIRECT lie on the upper edge curve (Figure~\ref{fig:WeightObjDecision-Constrained-DIRECT}). Such differences arise because the~$\MCF$ (and, similarly, the~$\MCFg$) is non-convex for this problem, so different criteria and different solvers can converge to different local minimizers along an edge curve of the Pareto front, rather than to a single global minimizer. It becomes evident that knee solutions can have a local or global nature.}

Interestingly, the results for problem~VFM1\textit{constr} shown in \tcb{Figure~\ref{fig:Convergence-Constrained-MADS}} suggest that a minimizer of the~$\MCF$ may not correspond to a minimizer of the~\tcb{MCI}. This discrepancy \tcb{for~$\MCF$} arises because the ellipsoidal neighborhood becomes degenerate at the weight vector corresponding to the knee solution of problem~\eqref{prob:point_based_formulation_exact} (and nearby points). The degeneracy occurs because the matrix~$\nabla_{\lambda} (F(x(\lambda)))$ \tcb{at MCF's minimizer} is highly ill-conditioned, with two of its three eigenvalues close to zero, resulting in an elongated ellipsoid that inflates the value of the~MCI. This observation confirms that minimizing the~MCI instead of~$\MCF$ is not a reliable approach for finding knee solutions, as previously noted in Remark~\ref{rem:MCFvsMCM} of Subsection~\ref{subsec:results_unc}.

\section{Concluding remarks and future work}\label{sec:future_work} 

In this paper, we showed how to use first-order rates of variation of Pareto optimal solutions to better understand the tradeoffs \tcb{along} a Pareto front in multi-objective optimization.
We have used such sensitivity rates to provide an answer to the open question of how to rigorously compute knee solutions \tcb{according to the verbal definition}, i.e., solutions corresponding to nondominated points where a small improvement in any objective leads to a large deterioration in at least one other objective. 
\tcb{We then introduced two optimization formulations: one that minimizes the MCF directly, given in~\eqref{prob:point_based_formulation_exact}, and another that minimizes its computationally simpler gradient-based proxy, the gMCF, given in~\eqref{prob:point_based_formulation}.}
The \tcb{cornerstone} of our approach is the ability to compute the gradient~$\nabla_{\lambda}( f_i(x(\lambda)))$ of each individual function~$f_i(x(\cdot))$ through the rate of variation~$\nabla x(\lambda)$ of the Pareto optimal solution~$x(\lambda)$, where~$\lambda$ denotes a vector of objective weights. 
\tcb{We used these formulations} to compute knee solutions \tcb{in both unconstrained and constrained settings}.
We called our approach {\it snee} to emphasize the use of Pareto sensitivity in the calculation of the derivatives of~$F(x(\lambda))$.

We also showed how to compute the most-changing Pareto sub-fronts around a nondominated point. Such sub-fronts are obtained from neighborhoods of weight vectors constructed using the pseudo-inverse of the Jacobian matrix of the vector function~$F(x(\cdot))$ and include nondominated points that are distributed along directions of maximum change. We considered the neighborhood achieving the highest value of the~\tcb{MCI} to be the most effective, and we evaluated the proposed neighborhoods numerically rather than analytically, as such a theoretical study would be out of the scope of this paper.

The techniques used in our approach are restricted to scalarized methods. In the current paper, we explored the weighted-sum method, which requires convexity of the objective functions in~$F$ for a complete coverage of the Pareto front.
Our approach can be extended to the~$\varepsilon$-constrained method, which converts a multi-objective problem into a single-objective constrained problem where one objective is optimized subject to constraints requiring the other objectives to be below varying thresholds. 
Given a positive threshold~$\varepsilon$ and denoting the corresponding Pareto optimal solution as~$x(\varepsilon)$, one can consider a neighborhood of~$\varepsilon$ to determine a Pareto neighborhood around~$x(\varepsilon)$.
\tcb{Knee solutions could then be computed using our {\it snee} approach by reformulating problems~\eqref{prob:point_based_formulation_exact} and~\eqref{prob:point_based_formulation}} in terms of~$\varepsilon$.
Although the~$\varepsilon$-constrained method does not explicitly require convexity for a complete coverage of the Pareto front, it requires a global solution of a constrained optimization problem. Therefore, our future work will focus on exploring techniques that do not rely on scalarized methods.

\section*{Acknowledgments}
This work is partially supported by the U.S. Air Force Office of Scientific Research~(AFOSR) award~FA9550-23-1-0217 and the U.S. Office of Naval Research~(ONR) award~N000142412656. We thank the two anonymous referees and the associated editor for their careful reading and constructive comments, which helped improve the paper.








  



\newpage

\begin{appendices}

\section{Results obtained with the Nelder--Mead and DIRECT Algorithms}
\label{app:nm-direct}

\tcb{This appendix reports the results obtained with the Nelder--Mead~(NM) and DIRECT algorithms, both run with default parameters, for the same test problems, initial weight vectors, and neighborhood settings described in Sections~\ref{sec:unconstr} and~\ref{sec:constr}.}

\section{Unconstrained problems}

\begin{figure}[h]
    \centering
    \includegraphics[width=0.95\textwidth]{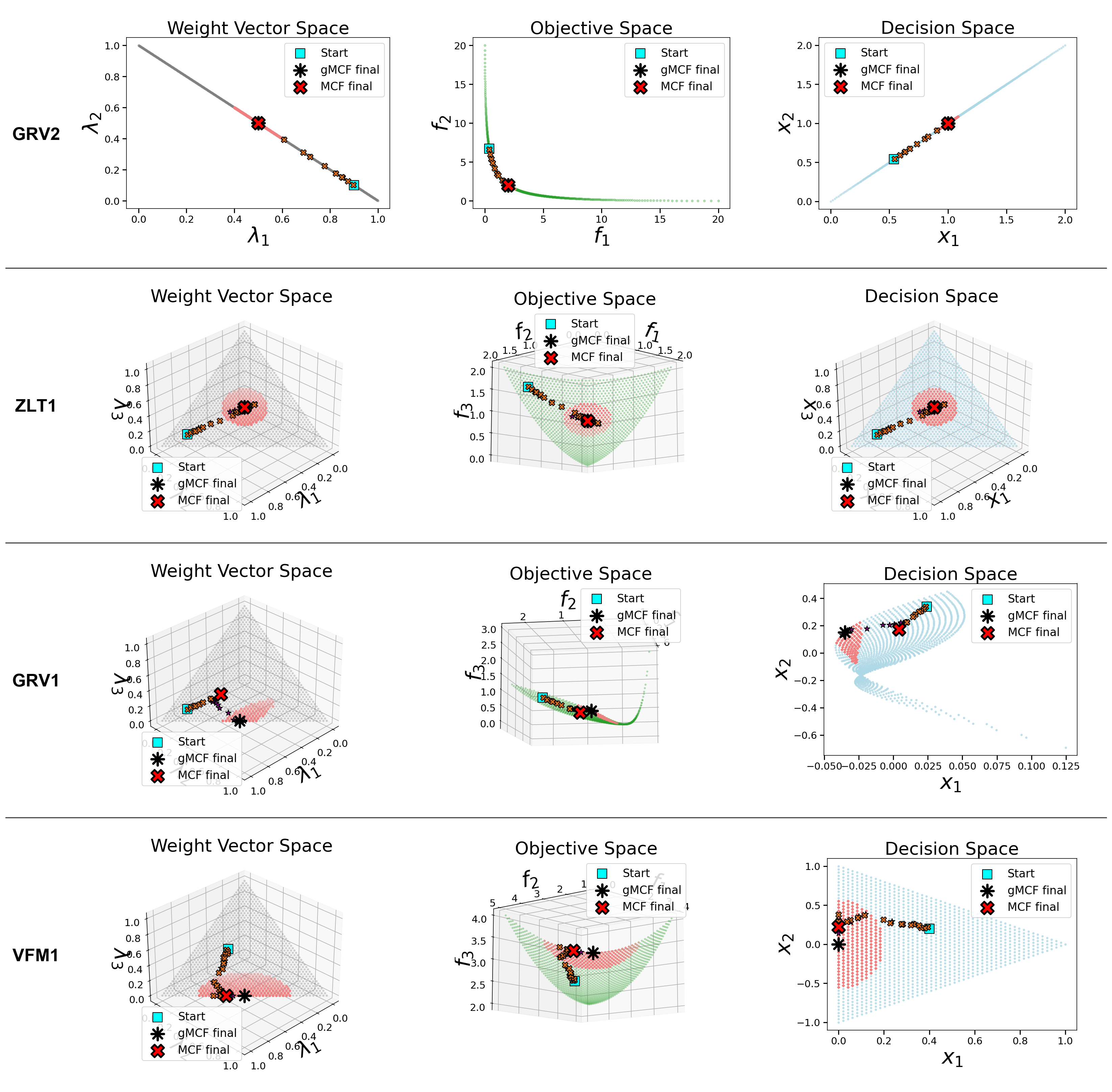}
    \caption{\tcb{Weight vector, objective, and decision spaces obtained when applying the NM algorithm, for problems GRV2, ZLT1, GRV1, and VFM1. See the caption of Figure~\ref{fig:WeightObjDecision-MADS} for the marker legend.}}
    \label{fig:WeightObjDecision-NM}
\end{figure}

\begin{figure}[h]
    \centering
    \includegraphics[width=0.95\textwidth]{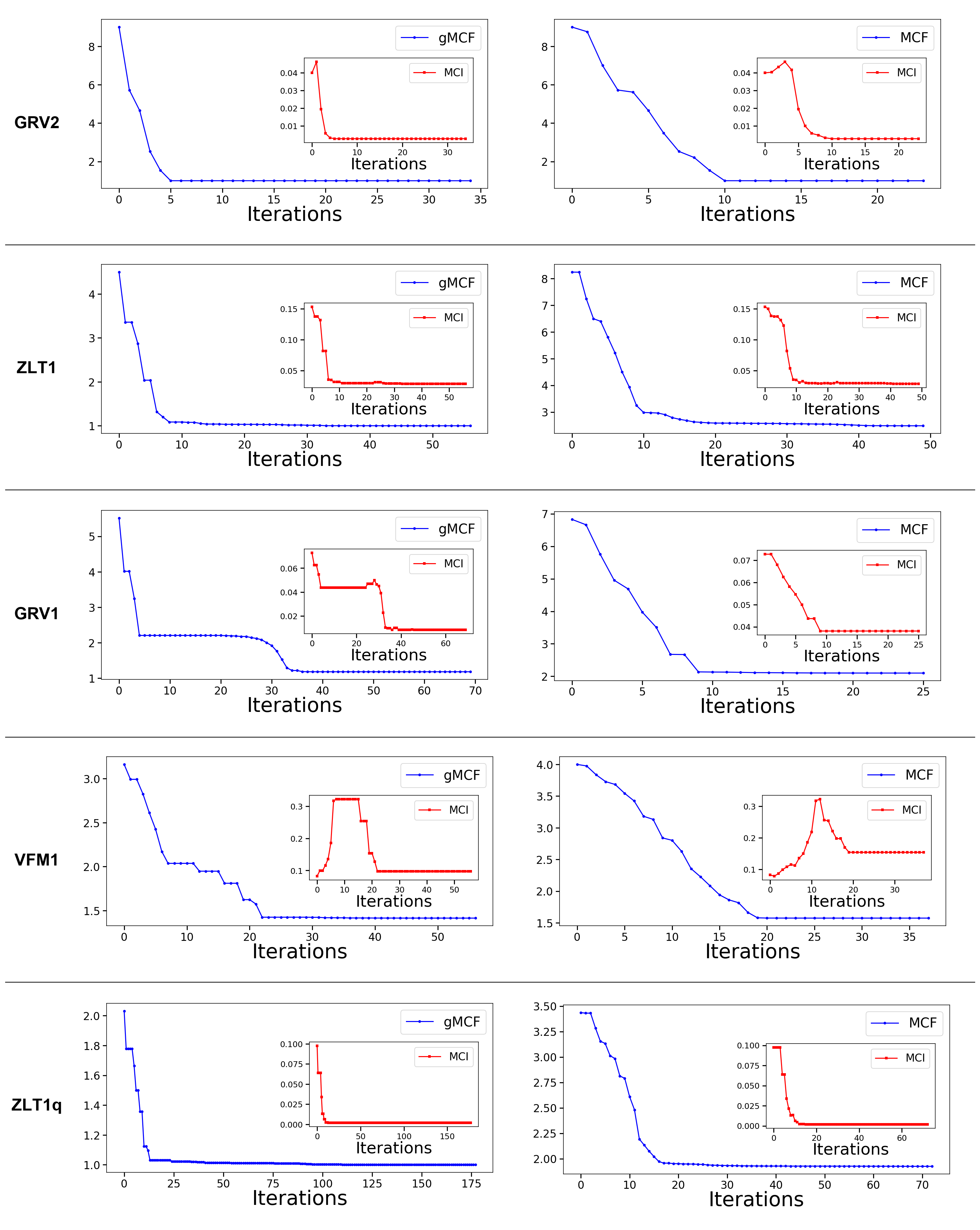}
    \caption{\tcb{Values of~$\MCFg$ and~$\MCF$ over the NM iterations, for problems GRV2, ZLT1, GRV1, VFM1, and ZLT1$q$.}}
    \label{fig:Convergence-Unconstrained-NM}
\end{figure}

\begin{figure}[h]
    \centering
    \includegraphics[width=0.95\textwidth]{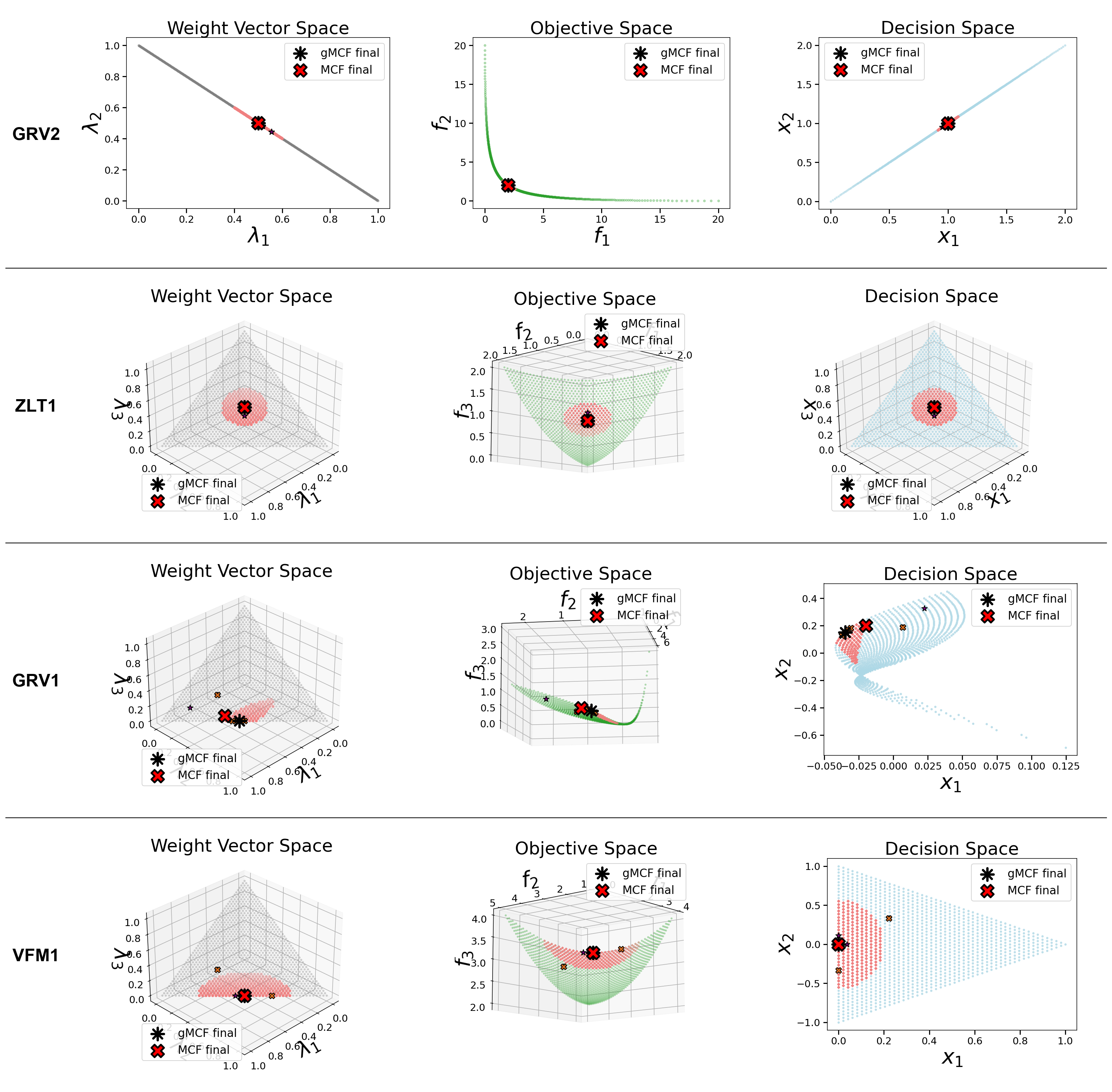}
    \caption{\tcb{Weight vector, objective, and decision spaces obtained when applying the DIRECT algorithm, for problems GRV2, ZLT1, GRV1, and VFM1. See the caption of Figure~\ref{fig:WeightObjDecision-MADS} for the marker legend.}}
    \label{fig:WeightObjDecision-DIRECT}
\end{figure}

\begin{figure}[h]
    \centering
    \includegraphics[width=0.95\textwidth]{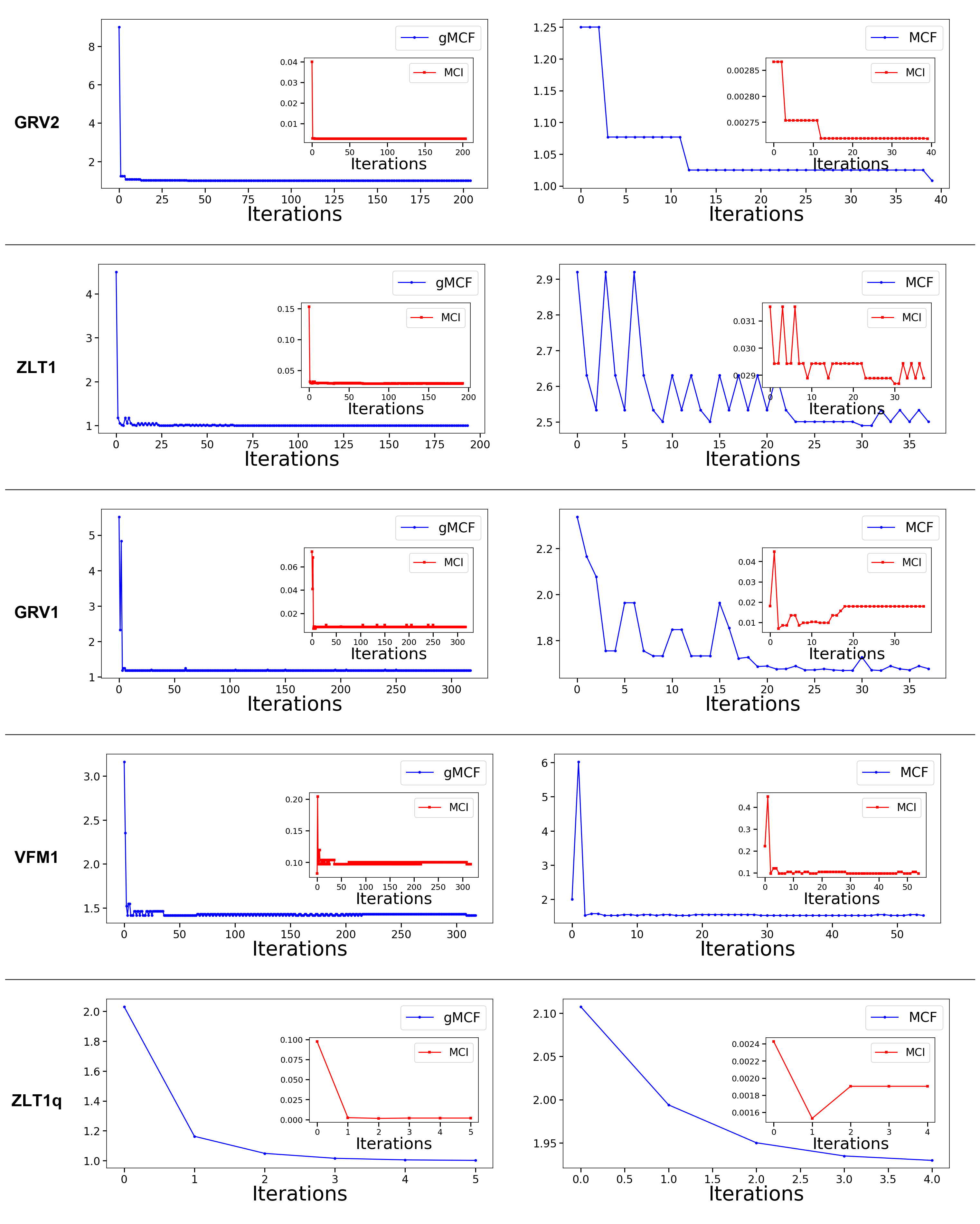}
    \caption{\tcb{Values of~$\MCFg$ and~$\MCF$ over the DIRECT iterations, for problems GRV2, ZLT1, GRV1, VFM1, and ZLT1$q$.}}
    \label{fig:Convergence-Unconstrained-DIRECT}
\end{figure}

\tcb{Figures~\ref{fig:WeightObjDecision-NM} and~\ref{fig:Convergence-Unconstrained-NM} show, respectively, the weight vector, objective, and decision spaces, and the values of~$\MCFg$ and~$\MCF$ over the iterations, obtained when applying the NM algorithm to problems GRV2, ZLT1, GRV1, VFM1, and ZLT1$q$. Figures~\ref{fig:WeightObjDecision-DIRECT} and~\ref{fig:Convergence-Unconstrained-DIRECT} report the analogous results for the DIRECT algorithm. The initial weight vector~$\lambda_0$ for each problem is the same as in Subsection~\ref{subsec:results_unc}. As in Figure~\ref{fig:WeightObjDecision-MADS}, the black asterisk marks the knee obtained by minimizing~$\MCFg$ (``gMCF final'') and the red cross marks the knee obtained by minimizing~$\MCF$ (``MCF final'').
Comparing Figures~\ref{fig:WeightObjDecision-NM} and~\ref{fig:WeightObjDecision-DIRECT} with Figure~\ref{fig:WeightObjDecision-MADS}, the knee nondominated points found by NM and DIRECT coincide with those found by MADS for all five problems, and the ``gMCF final'' and ``MCF final'' markers agree with one another to the same degree as under MADS. Similarly, the~$\MCFg$, $\MCF$, and~$\MCI$ values in Figures~\ref{fig:Convergence-Unconstrained-NM} and~\ref{fig:Convergence-Unconstrained-DIRECT} converge to final values similar to those in Figure~\ref{fig:Convergence-Unconstrained-MADS}.}

\subsection{Constrained problems}

\begin{figure}[h]
    \centering
    \includegraphics[width=0.95\textwidth]{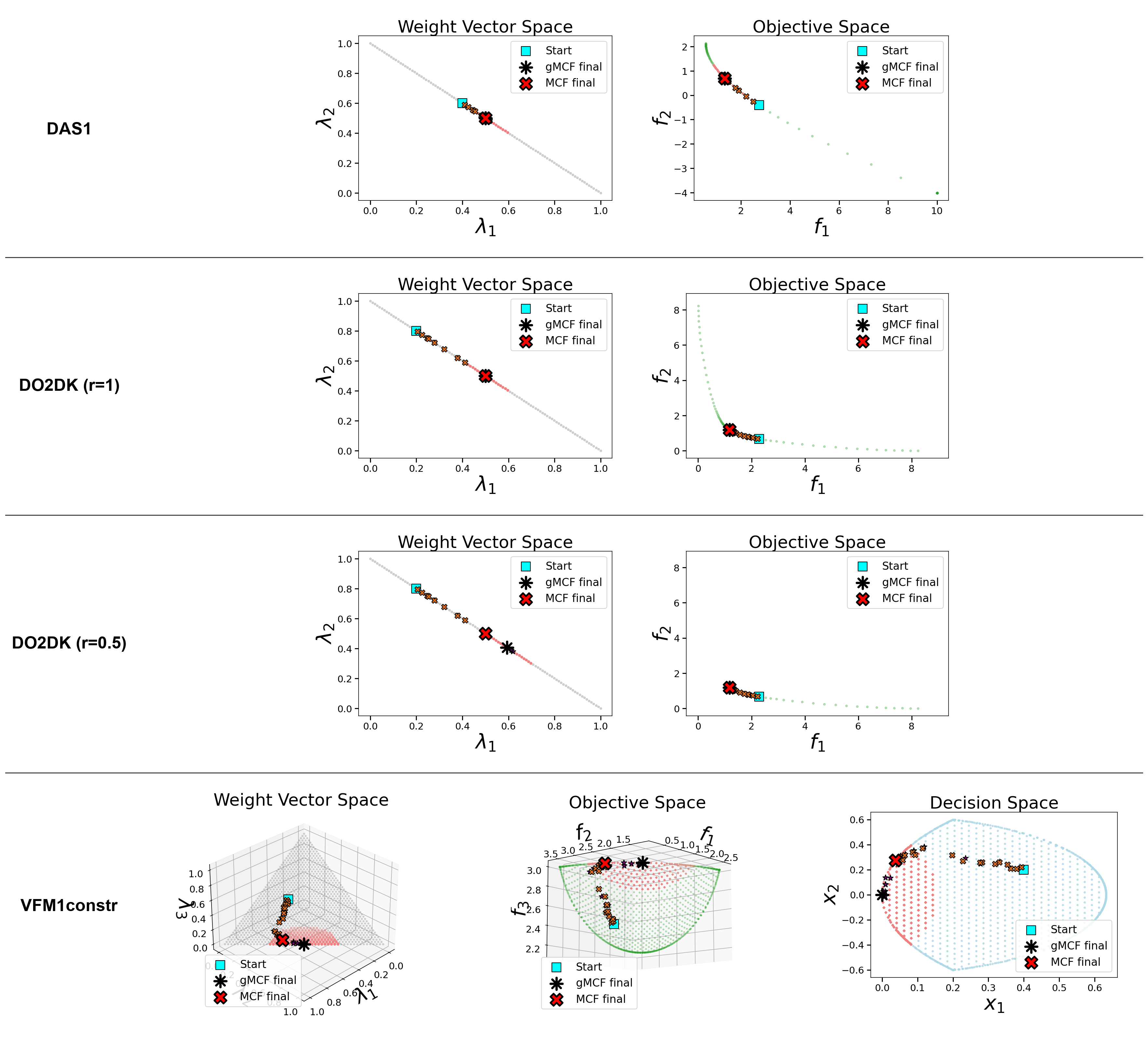}
    \caption{\tcb{Weight vector, objective, and decision spaces obtained when applying the NM algorithm, for problems DAS1, DO2DK ($r=1$), DO2DK ($r=0.5$), and VFM1\textit{constr}.}}
    \label{fig:WeightObjDecision-Constrained-NM}
\end{figure}

\begin{figure}[h]
    \centering
    \includegraphics[width=0.95\textwidth]{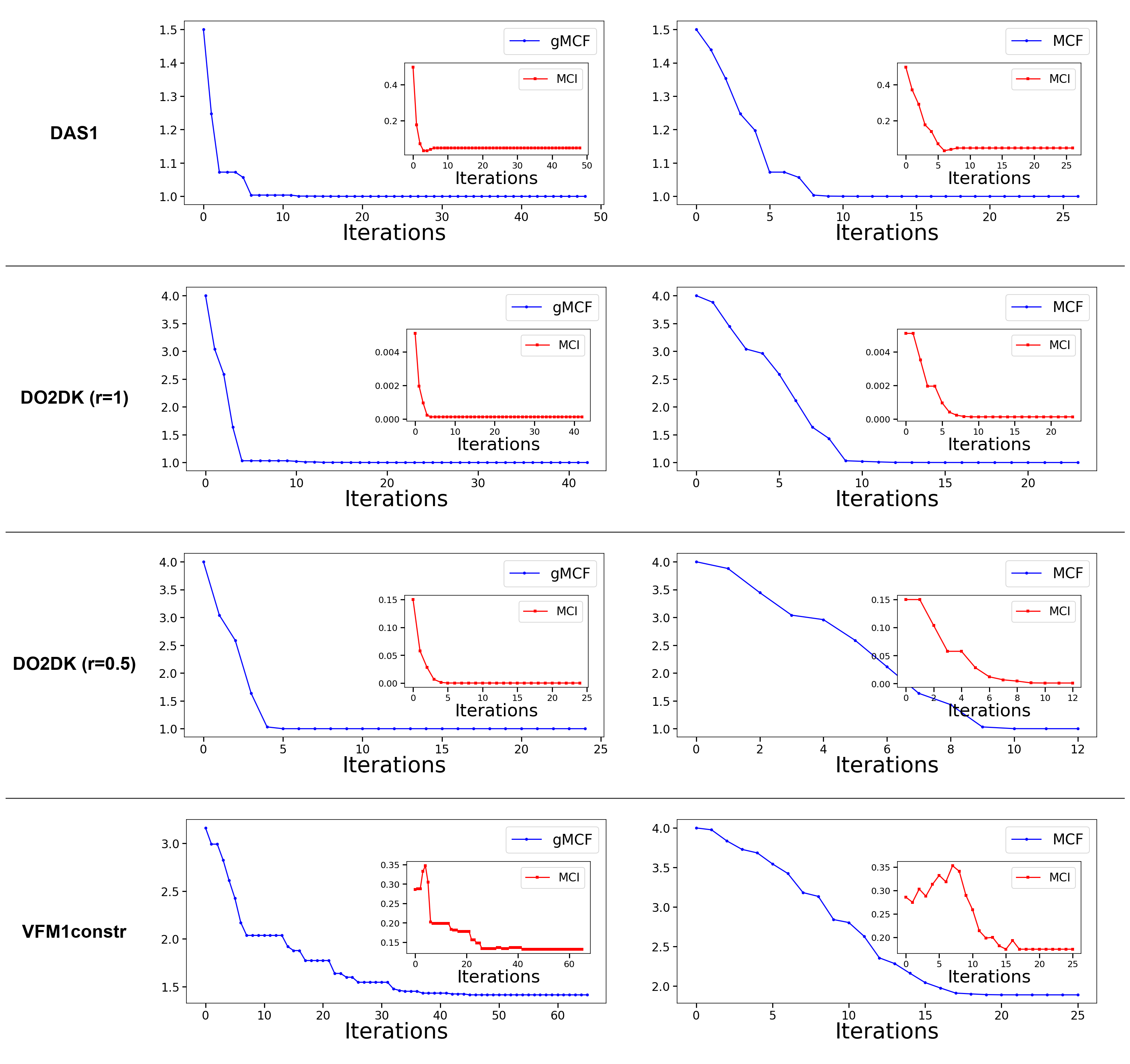}
    \caption{\tcb{Values of~$\MCFg$ and~$\MCF$ over the NM iterations, for problems DAS1, DO2DK ($r=1$), DO2DK ($r=0.5$), and VFM1\textit{constr}.}}
    \label{fig:Convergence-Constrained-NM}
\end{figure}

\begin{figure}[h]
    \centering
    \includegraphics[width=0.95\textwidth]{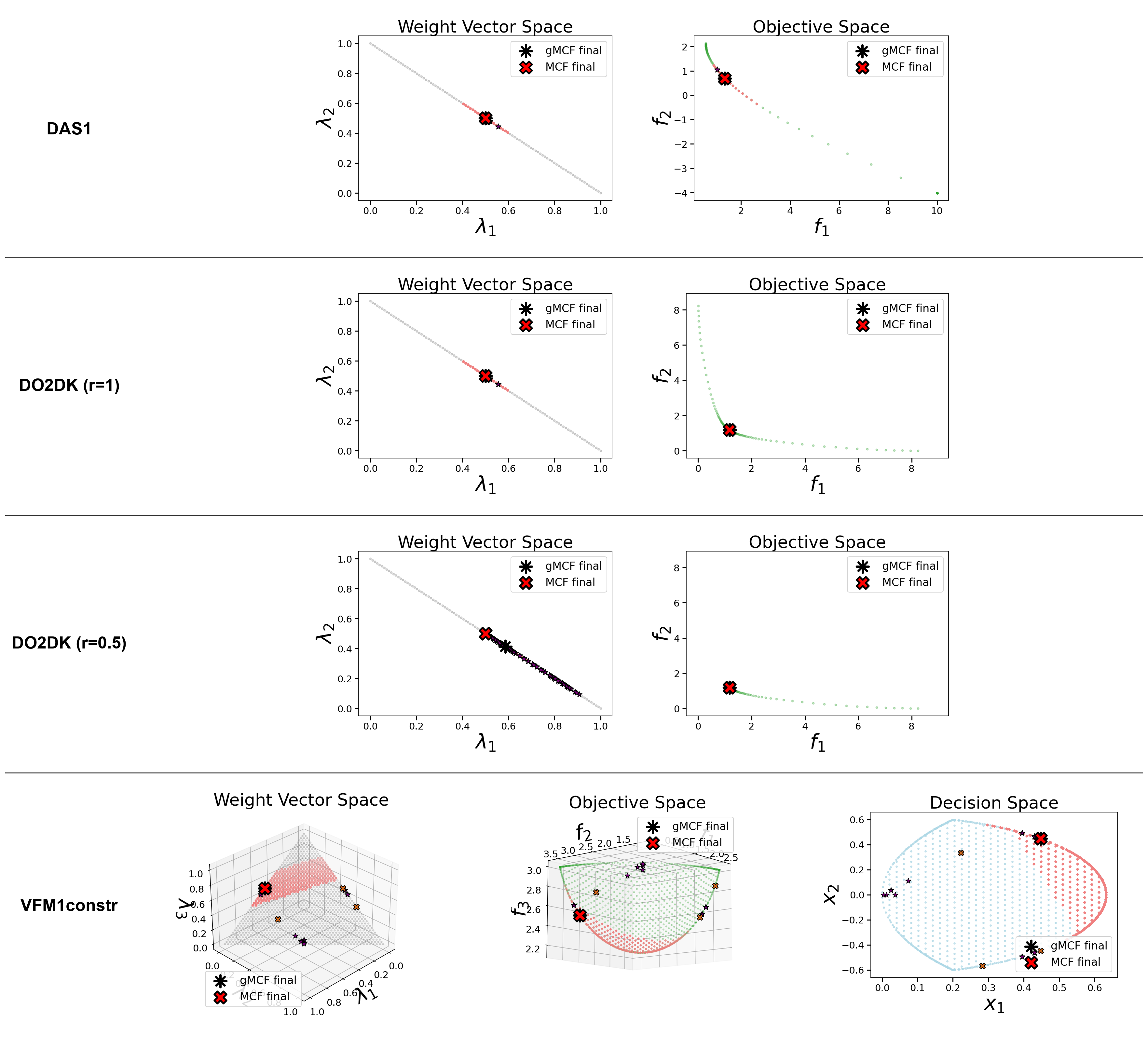}
    \caption{\tcb{Weight vector, objective, and decision spaces obtained when applying the DIRECT algorithm, for problems DAS1, DO2DK ($r=1$), DO2DK ($r=0.5$), and VFM1\textit{constr}.}}
    \label{fig:WeightObjDecision-Constrained-DIRECT}
\end{figure}

\begin{figure}[h]
    \centering
    \includegraphics[width=0.95\textwidth]{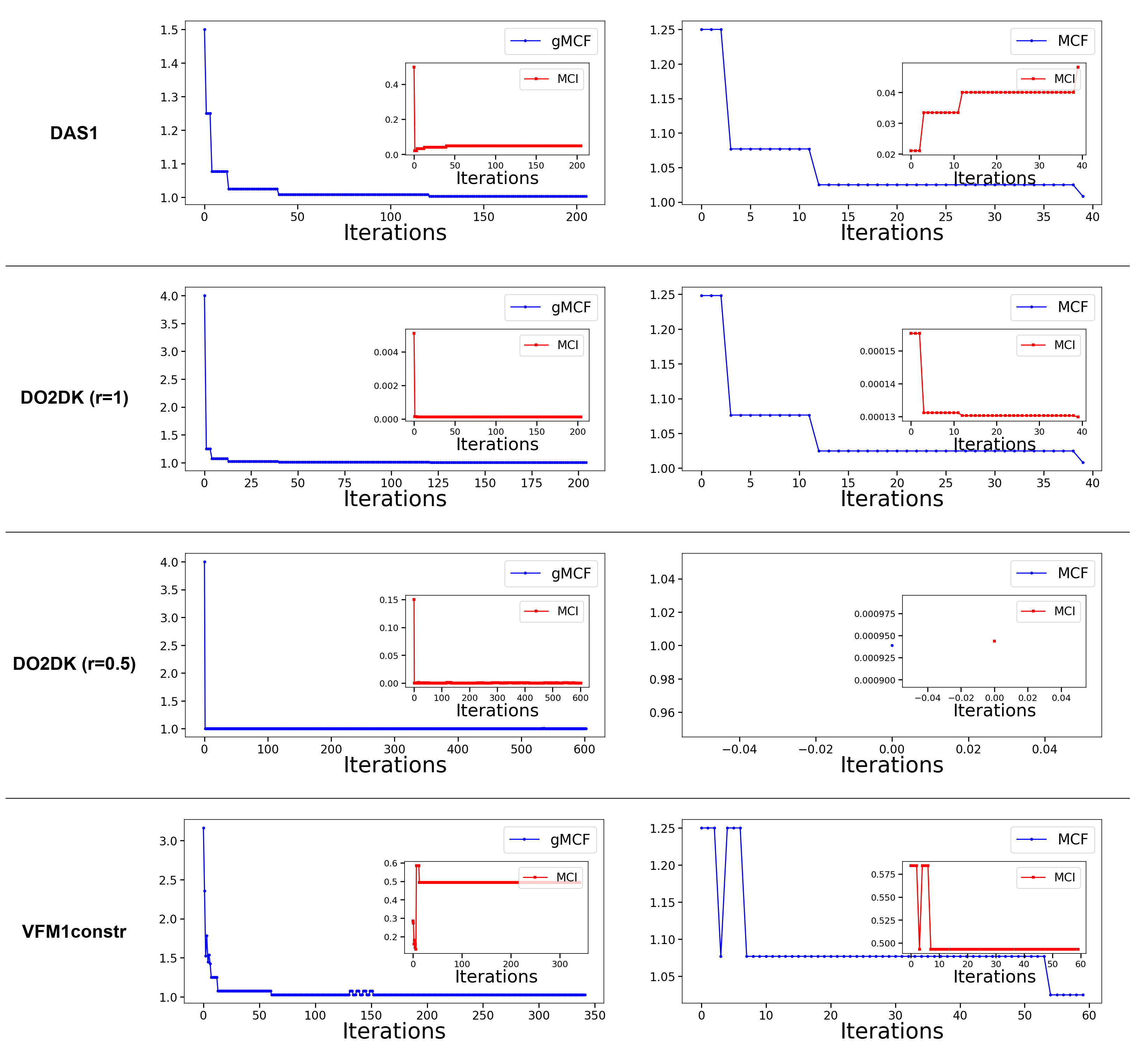}
    \caption{\tcb{Values of~$\MCFg$ and~$\MCF$ over the DIRECT iterations, for problems DAS1, DO2DK ($r=1$), DO2DK ($r=0.5$), and VFM1\textit{constr}.}}
    \label{fig:Convergence-Constrained-DIRECT}
\end{figure}

\tcb{Figures~\ref{fig:WeightObjDecision-Constrained-NM} and~\ref{fig:Convergence-Constrained-NM} report the NM results for problems DAS1, DO2DK ($r=1$), DO2DK ($r=0.5$), and VFM1\textit{constr}; Figures~\ref{fig:WeightObjDecision-Constrained-DIRECT} and~\ref{fig:Convergence-Constrained-DIRECT} report the analogous results for DIRECT. The initial weight vector~$\lambda_0$ and neighborhood settings for each problem/configuration are the same as in Subsection~\ref{subsec:results_constr}. For problems~DAS1 and~DO2DK, under both configurations of~$r$, the knee nondominated points found by NM and DIRECT coincide with those found by MADS in Figure~\ref{fig:WeightObjDecision-Constrained-MADS}, and the ``gMCF final'' and ``MCF final'' markers agree closely with one another, as they do under MADS.
For problem~VFM1\textit{constr}, however, the difference between the ``gMCF final'' and ``MCF final'' knees observed with MADS in Subsection~\ref{subsec:results_constr} is also present with~NM, and in fact NM and DIRECT identify different knees. The knees found by DIRECT lie on the upper edge curve of the Pareto front (Figure~\ref{fig:WeightObjDecision-Constrained-DIRECT}), whereas the knees found by NM lie on the left edge curve (Figure~\ref{fig:WeightObjDecision-Constrained-NM}). This is consistent with the discussion in Subsection~\ref{subsec:results_constr}. The corresponding~$\MCI$ values in Figure~\ref{fig:Convergence-Constrained-NM} and Figure~\ref{fig:Convergence-Constrained-DIRECT} confirm that minimizers of~$\MCF$ or~$\MCFg$ do not correspond to minimizers of~$\MCI$ for problem~VFM1\textit{constr}, for the same ill-conditioning discussed in Subsection~\ref{subsec:results_constr}. For problem~DO2DK~($r=0.5$), nearly every candidate weight vector encountered by DIRECT is discarded by the conditioning guard for the~$\MCF$ described in Subsection~\ref{subsec:results_constr}. Consequently, the~$\MCF$-vs-iterations plot for this problem in Figure~\ref{fig:Convergence-Constrained-DIRECT} reduces to a single point.
}

\end{appendices}

\end{document}